\documentclass[10pt,a4paper,reqno]{article}
\usepackage[numbers]{natbib}
\usepackage{graphicx, latexsym}
\usepackage{textcomp}
\usepackage{amsmath,amsfonts,amscd,amsthm,amssymb}
\usepackage{amscd}
\usepackage{epsfig}
\usepackage{epstopdf}
\usepackage{makeidx}
\usepackage{color} 
\usepackage{titlesec}
\usepackage{hyperref}
\usepackage{calc}
\usepackage{accents}
\usepackage{ulem}
\usepackage{empheq}
\usepackage{subcaption}
\usepackage{tikz}
\usepackage{titling} 
\usepackage{lipsum} 
\usepackage{authblk} 
\usepackage{wrapfig}

\makeindex
\newtheorem{theorem}{Theorem}[section]

\newtheorem{Remark}{Remark}[section]
\newtheorem{lemma}{Lemma}[section]

\numberwithin{equation}{section}

\newcommand{\R}{\mathbb{R}}                                                                
\newcommand{\ms}{\mathrm{s}}

\newcommand{\bX}{{\bf X}}
\newcommand{\bw}{{\bf w}}

\newcommand{\bP}{{\bf P}}
\newcommand{\bv}{{\bf v}}

\newcommand{\bV}{{\bf V}}

\newcommand{\bH}{{\bf H}}
\newcommand{\bL}{{\bf L}}

\newcommand{\f}{{\bf f}}
\newcommand{\bu}{{\bf u}}

\newcommand{\bn}{{\bf n}}
\newcommand{\bU}{{\bf U}}

\newcommand{\bK}{{\bf K}}
\newcommand{\xb}{\color{blue}}

\newcommand{\bz}{{\bf z}}
\newcommand{\bQ}{{\bf Q}}
\newcommand{\bO}{{\bf O}}
\newcommand{\bx}{{\bf x}}
\newcommand{\bZ}{{\bf Z}}
\newcommand{\bY}{{\bf Y}}
\newcommand{\bW}{{\bf W}}
\newcommand{\p}{ \boldsymbol{ \mathcal{P}}}
\newcommand{\Z}{ \boldsymbol{ \mathcal{X}}}

\newcommand{\bPi}{\mbox{\boldmath $\Pi$}}

\newcommand{\bOa}{\mbox{\boldmath $\textbf{O}_{fh}^{\bs}$}}
\newcommand{\bOb}{\mbox{\boldmath $\textbf{O}_{fh}^{\bu}$}}
\newcommand{\bOc}{\mbox{\boldmath $\textbf{O}_{ph}^{\bomega}$}}
\newcommand{\Od}{\mbox{ $O_{ph}^{\phi}$}}

\newcommand{\bTu}{\mbox{\boldmath ${\Theta}^{\bu}$}}

\newcommand{\bTsi}{{\mbox{\boldmath $\underline{\Theta}^{\bs i}$}}}
\newcommand{\bTui}{{\mbox{\boldmath ${\Theta}^{\bu i}$}}}
\newcommand{\bTqi}{{\mbox{\boldmath ${\Theta}^{\bomega i}$}}}
\newcommand{\Tpi}{{\mbox{ ${\Theta}^{\phi i}$}}}
\newcommand{\bTuI}{{\mbox{\boldmath ${\Theta}^{\bu (i-1)}$}}}

\newcommand{\bUpu}{\mbox{\boldmath ${\Upsilon}^{\bu}$}}

\newcommand{\Upp}{\mbox{ ${\Upsilon}^{\phi}$}}
\newcommand{\bups}{{\mbox{\boldmath$\underline{\Upsilon}^{\bs i}$}}}
\newcommand{\bupu}{{\mbox{\boldmath ${\Upsilon}^{\bu i}$}}}
\newcommand{\bupq}{{\mbox{\boldmath ${\Upsilon}^{\bomega i}$}}}
\newcommand{\upp}{{\mbox{ ${\Upsilon}^{\phi i}$}}}
\newcommand{\bUpsn}{{\mbox{\boldmath $\underline{\Upsilon}^{\bs n}$}}}
\newcommand{\bUpun}{{\mbox{\boldmath ${\Upsilon}^{\bu n}$}}}
\newcommand{\bUpqn}{{\mbox{\boldmath ${\Upsilon}^{\bomega n}$}}}
\newcommand{\Uppn}{{\mbox{ ${\Upsilon}^{\phi n}$}}}

\newcommand{\bUpui}{{\mbox{\boldmath ${\Upsilon}^{\bu (i-1))}$}}}

\newcommand{\Uppi}{{\mbox{ ${\Upsilon}^{\phi (i-1)}$}}}

\newcommand{\bUpsz}{{\mbox{\boldmath $\underline{\Upsilon}^{\bs 0}$}}}
\newcommand{\bUpuz}{{\mbox{\boldmath ${\Upsilon}^{\bu 0}$}}}

\newcommand{\Uppz}{{\mbox{ ${\Upsilon}^{\phi 0}$}}}

\newcommand{\bca}{{\mbox{\boldmath $\underline{\bzeta}^{\bs i}$}}}
\newcommand{\bcb}{{\mbox{\boldmath $\bzeta^{\bu i}$}}}
\newcommand{\bcc}{\mbox{\boldmath $\bzeta^{\bomega i}$}}
\newcommand{\cd}{{\mbox{ $\zeta^{\phi i}$}}}
\newcommand{\bcbI}{{\mbox{\boldmath $\bzeta^{\bu {(i-1)}}$}}}

\newcommand{\bKa}{\mbox{\boldmath $\underline{\bchi}^{\bs}$}}

\newcommand{\Kd}{\mbox{${\chi}^{\phi}$}}
\newcommand{\bka}{{\mbox{\boldmath$\underline{\chi}^{\bs i}$}}}
\newcommand{\bkb}{{\mbox{\boldmath ${\chi}^{\bu i}$}}}
\newcommand{\bkc}{\mbox{\boldmath ${\bchi}^{\bomega i}$}}
\newcommand{\kd}{\mbox{${\chi}^{\phi i}$}}
\newcommand{\bkA}{{\mbox{\boldmath$\underline{\chi}^{\bs 1}$}}}
\newcommand{\bkB}{{\mbox{\boldmath ${\chi}^{\bu 1}$}}}
\newcommand{\bkC}{\mbox{\boldmath ${\bchi}^{\bomega 1}$}}
\newcommand{\kD}{\mbox{${\chi}^{\phi 1}$}}

\newcommand{\bsq}{\boldsymbol{[}}
\newcommand{\bqs}{\boldsymbol{]}}
\newcommand{\bPa}{\mbox{\boldmath $\textbf{P}_{h}^{\bs}$}}
\newcommand{\bPb}{\mbox{\boldmath $\textbf{P}_{h}^{\bu}$}}
\newcommand{\bPc}{\mbox{\boldmath $\textbf{P}_{h}^{\bomega}$}}
\newcommand{\Pd}{\mbox{ $P_{h}^{\phi}$}}
\newcommand{\bDs}{\mbox{\boldmath ${\delta}^{\bs}$}}
\newcommand{\bDu}{\mbox{\boldmath ${\delta}^{\bu}$}}
\newcommand{\bDq}{\mbox{\boldmath ${\delta}^{\bomega}$}}
\newcommand{\Dp}{\mbox{ ${\delta}^{\phi}$}}

\newcommand{\bomega}{\mbox{\boldmath $\omega$}}
\newcommand{\bbeta}{\mbox{\boldmath $\beta$}}

\newcommand{\bnu}{\mbox{\boldmath $\nu$}}

\newcommand{\bs}{\mbox{\boldmath $\sigma$}}
\newcommand{\bSig}{\mbox{\boldmath $\Sigma$}}
\newcommand{\btau}{\mbox{\boldmath $\tau$}}

\newcommand{\bzeta}{\mbox{\boldmath $\zeta$}}
\newcommand{\bchi}{\mbox{\boldmath $\chi$}}

\newcommand{\bl}{\mbox{\boldmath $\lambda$}}

\newcommand{\e}{{\bf e}}

\DeclareMathOperator*{\esssup}{ess\,sup}

\newcommand{\xm}{\color{magenta}}

\newcommand{\vertiii}[1]{{\left\vert\kern-0.25ex\left\vert\kern-0.25ex\left\vert #1 
    \right\vert\kern-0.25ex\right\vert\kern-0.25ex\right\vert}}
    
     \newcommand{\correspondingauthor}{\textsuperscript{*}Corresponding author}

    \author[1]{\small Sweta Chakraborty}
\author[2]{\small Deepjyoti Goswami\thanks{\correspondingauthor}}
\author[3]{\small Saumya Bajpai}

\affil[1]{Department of Mathematical Sciences, Tezpur University, Tezpur, Sonitpur, Assam-784028, India. Email: msp21103@tezu.ac.in}
\affil[2]{Department of Mathematical Sciences, Tezpur University, Tezpur, Sonitpur, Assam-784028, India. Email: deepjyoti@tezu.ernet.in}
\affil[3]{School of Mathematics and Computer Science, Indian Institute of Technology Goa, Ponda, Goa-403401, India. Email: saumya@iitgoa.ac.in}

\date{}

\begin{document}
\thispagestyle{plain}

\title
{\large\bf A locally conservative LDG method for the evolutionary Navier-Stokes-Darcy coupled problem}
%
\maketitle
\noindent\textbf{Abstract:} This paper proposes and investigates a local discontinuous Galerkin (LDG) method for the evolutionary  Navier-Stokes-Darcy coupled problem representing the interaction between incompressible flows and porous media flows through the Beavers-Joseph-Saffman interface conditions. A fully discrete numerical scheme has been proposed and analyzed employing backward Euler temporal discretization and a locally conservative LDG scheme for the space discretization.  
The well-posedness of the numerical scheme has been obtained under a small data condition. Optimal {{\it a priori}} error estimates have been proven in the energy norm, as well as the $\bL^2$-norm for $t>0$. Numerical experiments ensure the efficiency of the theoretical results. 

 \vspace{0.25cm}
\vspace*{-2mm}
\section{Introduction}
\noindent The combination of incompressible flows and porous media flows is becoming increasingly popular nowadays \cite{COV26, D13} due to the abundance of applications in various fields of interest. Such complicated phenomena have numerous practical applications in areas like geosciences (modeling river and groundwater interactions, simulating the effect of flooding in dry areas), health sciences (filtration of blood via vessel walls, modeling blood flow and organs), industry (petroleum extraction, air or oil filters),  bioengineering (blood oxygenators, hemodialysis devices, etc.), and so on. The coupling is an integral part of many such technical applications related to filtration. 
It is well known that the Navier-Stokes-Darcy (NSD) model demonstrates these coupled flows effectively. 
The Navier-Stokes equations (NSEs) in the fluid region are typically coupled across an interface with the Darcy equations for the filtration velocity in the porous media to provide a model that is most accessible to large-scale computations. As a consequence, however, the coupled system of equations of various orders in different regions causes mathematical complications; for example, only the classical continuity conditions no longer work in the interface. Beavers and Joseph, in \cite{BJ67}, provided an additional interface condition describing the tangential interaction, based on their experimental study. Later on, Saffman in \cite{S71} observed that the Darcy tangential velocity is negligible in comparison to free flow slip velocity, leading to the simplified Beavers-Joseph-Saffman (BJS) interface condition. 
The first rigorous mathematical studies on coupling conditions at the interface can be ascribed to J\"{a}ger and Mikeli\'{c} (see \cite{JMN01} and the references therein). \vspace*{-3mm}

\noindent 
These coupled models, namely, the Stokes-Darcy and the steady NSD, have been studied extensively, more importantly, for various numerical approximations, including continuous and discontinuous Galerkin methods, see, for example,  \cite{COV26} and the references cited therein.
In contrast, literature for the time-dependent coupled NSD problem is limited. The well-posedness of the problem has been studied in \cite{CGR13}, whereas the global existence, uniqueness, and the large-time behaviour of the strong solution have been obtained under an additional smallness assumption on the initial data in \cite{CDG24}.
And recently, various finite elements have been studied in search of suitable numerical schemes to handle the coupling through different interface conditions; see \cite{LD22, WLS18} and the references therein.

\noindent
We would like to emphasize that although discontinuous Galerkin (DG) methods have gained considerable prominence over the past few decades due to their novel features, like the capability to handle complicated geometry, flexibility for adaptivity, locally varying polynomial order, higher-order accuracy, and local mass conservation, etc., the DG literature for this coupled problem is limited. 
For time-dependent NSEs and steady Darcy equations,  
primal DG methods are considered in \cite{CR09, CGPR17} along with the Crank-Nicolson and backward Euler schemes for time discretization, respectively.  
Recently, in \cite{CLR23}, a strongly conservative and pressure-robust hybridizable discontinuous Galerkin method for the transient NSEs coupled with the steady Darcy problem has been discussed. 

\noindent The LDG methods are one of the several DG methods with a wide range of applications, having properties of local conservativity and high degree of locality. Unlike the standard DG methods, which necessitate artificial stabilization to handle second-order terms, the LDG methods inherently transform these equations into first-order systems by the introduction of auxiliary variables. This leads to stability, flexibility in designing numerical fluxes, and accurate treatment of boundary and interface conditions.
 
\noindent The LDG methods have been introduced by Cockburn and Shu in \cite{CS98} as an extension to general convection-diffusion problems of the numerical scheme for compressible Navier-Stokes equations. 
These methods have then been extended to elliptic \cite{CCPS00}, Stokes \cite{CKSS02}, Oseen \cite{CKS04}, and steady NSEs \cite{CKS05}, and have been studied rigorously. For the time-dependent incompressible NSEs, in \cite{WDYW17}, the authors have developed a numerical scheme combining the characteristic method and the LDG method.
 
\noindent
The LDG method for the coupled Stokes-Darcy-transport problem has been presented in \cite{VD09}, where stability of the scheme and {\it a priori} error estimates have been analyzed. Dawson, in \cite{D06}, has studied a coupled ground water and surface water flow, where the model is approximated by the LDG method for steady ground water flow and a DG method for time dependent surface water flow, and {\it a priori} error estimates have been derived for this approach. \par
\noindent To the best of our knowledge, the analysis of the LDG method for the fully evolutionary NSD problem, where both problems are time dependent, remains unexplored, and our aim has been to fill this gap. At the very outset, the existence and uniqueness of the numerical scheme, along with the local conservation property, have been established. We note that the presence of the interface conditions disturbs the energy balance since the nonlinear convection of the NSEs is not completely eliminated (see (\ref{nonlinearity1})). This requires small data conditions for {\it a priori} error estimates; see \cite{CGPR17}.

\noindent 
Our analysis is based on discretizing the Oseen-Darcy equations, where the convective velocity is taken to be a projection of the approximate velocity following the approach of \cite{CKS05}. This approach leads our scheme to be locally conservative, as demonstrated in Remark \ref{locally conservative1}. Further, we have derived optimal error estimates for the velocity, pressure, and piezometric head.
%
This analysis 
 relies on the treatment of the interface along with the strong nonlinearity of the NSEs. 
Note that the use of standard projections yields optimal error bounds in the energy norm but fails in the case of $L^\infty(\bL^2)$-norm estimates, the reason being the appearance of the interface terms. The boundary integrals are estimated by means of trace inequalities, and the resulting negative power of $h$ is absorbed, in the DG setup, via the energy norms of broken spaces by including a jump term with negative power of $h$. However, in the case of NSD, in the definition of the energy norm (see (\ref{norm3})), the interface $\Gamma$ is not incorporated. As a result, the energy norm no longer compensates for $h^{-1/2}$ associated with the terms $\int_{\Gamma}$. This is a major bottleneck in our analysis, which also puts a restriction on analysing the nonlinear term effectively. For example,  
integration by parts of the nonlinear term is often performed to derive optimal error estimates, but now this would yield additional interface terms, which will only add to our difficulty. Our current work is devoted to overcoming these hurdles. \vspace*{1mm}

\noindent Regarding  $L^\infty(\bL^2)$ analysis, we have begun by introducing an $\bL^2$-projection $\bP_h$ onto a suitable LDG setup utilizing the approximation operator $\bPi_h$ and have derived the appropriate approximation properties. Next, we have constructed the modified Oseen operators in the LDG setup and established their optimal approximation properties by means of standard duality arguments. Then, temporal discretization has been done applying the backward Euler scheme, and fully discrete optimal $L^\infty(\textbf{L}^2)$ error estimates for velocity and $L^\infty({L}^2)$ error estimates for piezometric head have been obtained. Also, we have provided an optimal $L^{2}(L^2)$ and $L^{\infty}(L^2)$-pressure error estimate. Here also, we have overcome the technical challenges (see Remark \ref{pressure}).  
Our analysis has been performed for the two-dimensional case, but a similar analysis will go through for the three dimensional case as well, see Remark \ref{3D}. In fact, we would like to remark that all the estimates of our analysis are also valid for the steady NSD model, which are new results to the best of our knowledge. We end our discussion by providing numerical experiments to validate the theoretical results and the scheme's efficiency.

\noindent Below, we present a summary of the main findings from this article.
\begin{enumerate}
\item A fully discrete numerical scheme using the LDG method for a fully evolutionary NSD system has been proposed, and its existence, uniqueness, and local conservation property have been demonstrated. 
\item Modified Oseen operators have been established in an appropriate LDG setup, and their approximation properties have been explored, which plays a crucial role in tackling the interface terms while deriving the optimal error estimates in $L^\infty(\bL^2)$-norm for velocity and in $L^\infty(L^2)$-norm for piezometric head.
\item Optimal error estimates in the energy norm, for instance, $L^2(\bX_h)$ and $L^\infty(\bX_h)$-norm for fully discrete LDG approximate solutions to velocity, stress tensor, piezometric head, and Darcy velocity have been derived, where $\bX_h$-norm represents the appropriate broken energy norm.
\item Also, optimal error bounds in $L^\infty(\bL^2)$-norm for velocity and in $L^\infty(L^2)$-norm for piezometric head, along with a pressure error bound in $L^\infty(L^2)$-norm, have been obtained.
\item Numerical simulations have been done to validate our theoretical findings and the scheme's performance.
\end{enumerate}
\noindent The framework of this paper is as follows: Section \ref{model} describes the unsteady NSD model. Section \ref{s3} introduces the necessary functional spaces, notations, and key assumptions required for problem analysis, along with the mixed formulation. Section \ref{s4} presents the numerical scheme, its well-posedness, and discusses relevant trace inequalities. Section \ref{s05} contains standard approximation operators along with the establishment of modified Oseen operators. Section \ref{s5} mainly revolves around rigorous error estimates. Subsequently, Section \ref{s6} represents the numerical results. Finally, Section \ref{s7} summarizes the key points of the paper.
\vspace*{-2mm}
\section{The Fully Evolutionary Navier-Stokes-Darcy Model}
\label{model}
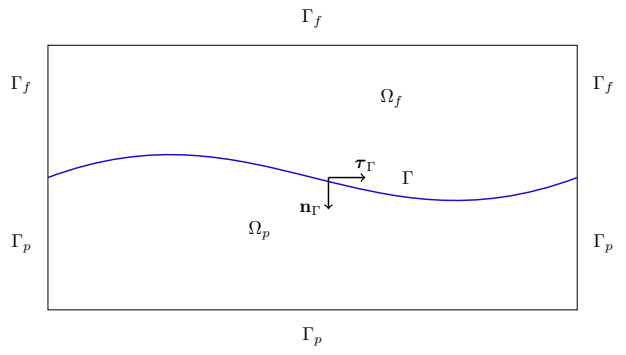
\begin{wrapfigure}{r}{0.25\textwidth}
    \centering
    \scalebox{0.35}{
    \begin{tikzpicture}
    \draw (0,0) rectangle (10,5);
 \draw [color={rgb,255:red,38; green,7; blue,211}, draw opacity=1, thick]
              (0,2.5) .. controls (4,4) and (6,1) .. (10,2.5)
              node[sloped,inner sep=0cm,above,pos=.235,
              anchor=south west,
              minimum height=(2.0)*0.3cm,minimum width=(3.0)*.3cm](N){};
 \coordinate (P1) at (4.7,2.5);
        \coordinate (P2) at (5.3,2.5);
\draw[->, thick, black] (P2) -- ++(0,-0.6) node[left] {$\bn_{\Gamma}$};
 \draw[->, thick, black] (P2) -- ++(0.7,0) node[above] {$\btau_{\Gamma}$};
 \node at (5, -0.5) {$\Gamma_{p}$};
        \node at (-0.5, 1.25) {$\Gamma_{p}$};
        \node at (10.5, 1.25) {$\Gamma_{p}$};
        \node at (5, 5.5) {$\Gamma_{f}$};
        \node at (-0.5, 4.25) {$\Gamma_{f}$};
        \node at (10.5, 4.25) {$\Gamma_{f}$};
        \node at (6.8, 2.5) {$\Gamma$};
        \node at (4, 1.5) {$\Omega_p$};
        \node at (6.5, 4) {$\Omega_f$};    
\end{tikzpicture}}
    \caption{Coupled domain of fluid flow and porous ~~media flow.} 
    \label{domain}
\end{wrapfigure}

Consider a coupled evolutionary NSD model in a bounded connected domain $\Omega=\Omega_f\cup \Omega_p\subset \mathbb{R}^d$ $(d=2,3)$, where $\Omega_f$  and $\Omega_p$ represent the fluid flow region and the porous media, respectively, satisfying $\Omega_f\cap\Omega_p=\emptyset$. 
We further assume that each subdomain is convex.
 Let $\partial\Omega$, $\partial\Omega_f$, and $\partial\Omega_p$ denote the boundaries of the domain $\Omega$ and subdomains $\Omega_f$ and $\Omega_p$, respectively.
 The interface between the two regions is given by $\Gamma=\bar{\Omega}_f\cap \bar{\Omega}_p$. 
Finally, the exterior boundaries of the fluid region and porous media are given by $\Gamma_f=\partial\Omega_f\cap\partial\Omega$ and $\Gamma_p=\partial\Omega_p\cap \partial\Omega$, respectively. {\xm See Figure \ref{domain} for a schematic diagram of two-dimensional geometry.}

\noindent The Navier-Stokes equations, governing the fluid flow in $\Omega_f$, are expressed as
\small{\begin{align} 
     \bu_t - \nabla\cdot(2\mu \textbf{D}(\bu)-p{\bf I})+ \bu \cdot \nabla \bu &= \f_f ~~\mbox{in}~ \Omega_f\times(0,T), \label{NSEDARCY1}\\
\nabla \cdot \bu &= 0 ~~\mbox{in}~  \Omega_f\times(0,T), \label{NSEDARCY2} \\
\bu(0)& = \bu_0 ~~\mbox{in}~ \Omega_f,\label{NSEDARCY3}
\end{align}}\normalsize
\noindent where $\bu_t$ denotes differentiation of $\bu$ in time. 
Here $\bu=(u_i)_{i=1}^d$ indicates the velocity vector, $p$ the pressure of the fluid in $\Omega_f$, the coefficient $\mu>0$ the kinetic viscosity, $\textbf{D}(\bu)=(1/2)(\nabla\bu+(\nabla\bu)^T)$ the deformation tensor with ${\bf I}$ being the identity tensor, and $T$ is the final time with $0<T<\infty$. The forcing term $\f_f$ is a given function in the domain $\Omega_f$. Other notations have the usual meaning as expressed below: \vspace*{-1mm}
\small{\[\bu\cdot\nabla\bu=\sum_{i=1}^d { u}_i\frac{\partial\bu}{\partial x_i}\quad \text{and}\quad \nabla\cdot \bu=\sum_{i=1}^d\frac{\partial u_i}{\partial x_i}.\]}
\noindent The porous media flow in $\Omega_p$ is governed by non-stationary Darcy equations, which are given by
\small{\begin{align}
\bomega &= -\bK\nabla \phi~~\mbox{in}~~ \Omega_p\times(0,T),\label{NSEDARCY4.0}\\
	S_0\phi_t+\nabla \cdot \bomega & = f_p ~~\text{in}~~\Omega_p\times(0,T), \label{NSEDARCY4}\\
	\phi(0) &= \phi_0~~ \mbox{in}~~ \Omega_p,\label{NSEDARCY5}
	\end{align}}
\noindent where 
$\phi$ represents the piezometric head, $\bomega$ the velocity of the porous media flow in $\Omega_p$, $S_0>0$ the specific mass storativity coefficient, and $f_p$ the forcing term in $\Omega_p$. $\bK\in\R^{d\times d}$ is the permeability tensor, which is presumed to be symmetric positive definite, with eigenvalues uniformly bounded above and away from zero. We impose homogeneous Dirichlet boundary conditions on $\Gamma_f$ and $\Gamma_p$:
\small{\begin{align}
	\bu ={\bf 0} \quad \mbox{on }\ \Gamma_f\times(0,T), \quad
	 \phi =0 \quad \mbox{on }\ {\Gamma_{p}}\times(0,T). \label{NSEDARCY9}
	\end{align}}
The two systems of equations, (\ref{NSEDARCY1})-(\ref{NSEDARCY3}) and (\ref{NSEDARCY4.0})-(\ref{NSEDARCY5}), describing the 
flow in each subdomain, are coupled by the following conditions across the interface $\Gamma$:
\small{\begin{align}
\bu\cdot\bn_{\Gamma} = \bomega\cdot\bn_{\Gamma} \quad &\mbox{on } \Gamma\times(0,T),\label{NSEDARCY10} \\
((-2\mu \textbf{D}(\bu)+p{\bf I})\bn_{\Gamma})\cdot\bn_{\Gamma} =  g\phi\quad &\mbox{on }\Gamma\times(0,T),\label{NSEDARCY11} \\
 \bu \cdot{\btau}^{j}_{\Gamma} = -2\mu G^j(\textbf{D}(\bu)\bn_{\Gamma})\cdot\btau^j_{\Gamma},\quad 1\leq j\leq d-1\quad &\mbox{on }~\Gamma\times(0,T), \label{NSEDARCY12}
\end{align}}
where $G^j=\frac{\mu\alpha}{(\bK\btau^j_{\Gamma},\btau^j_{\Gamma})^{1/2}}$, $\alpha>0$ is determined from experimental data, $\bn_{\Gamma}$ is the unit outward normal vector to $\Gamma$ pointing from $\Omega_f$ to $\Omega_p$, $\{\btau^j_{\Gamma}\}_{j=1}^{d-1}$ are linearly independent unit tangential vectors on $\Gamma$, and $g$ is the gravitational acceleration. The condition (\ref{NSEDARCY10}) presents the conservation of mass across the interface ${\Gamma}$. The second condition (\ref{NSEDARCY11}) guarantees the balance of forces across the interface, and the third one (\ref{NSEDARCY12}) is known as the BJS condition, which depicts that the tangential component of the velocity is proportional to the shear stress. For more details, we refer to \cite{LD22}.
\vspace*{-3mm}
\section{Preliminaries}\label{s3} 
For our subsequent use, we denote the $\mathbb{R}^d$-valued function spaces by boldface letters such as $\bH^1(\mathfrak{D})=(H^1(\mathfrak{D}))^d$ and $\bL^2{(\mathfrak{D})}=(L^2(\mathfrak{D}))^d$, as well as $\mathbb{R}^{d^2}$-valued function spaces such as ${ \underline{\bH}}^1(\mathfrak{D})=(H^1(\mathfrak{D}))^{d^2}$ and ${ \underline{\bL}}^2(\mathfrak{D})=(L^2(\mathfrak{D}))^{d^2}$, etc. We use the space of square integrable Lebesgue measurable functions $L^2(\mathfrak{D})$ with inner product ${(\cdot,\cdot)}$ and classical Sobolev spaces $W^{m,r}(\mathfrak{D})$ \cite{A75}, which are defined for any nonnegative integer $m$ and ${r}\geq 1$ as follows:
{\small\[W^{m,r}(\mathfrak{D})=\{v\in L^r(\mathfrak{D}):\, \partial^\ms v\in L^r(\mathfrak{D}) \,\,\forall\, |\ms|\leq m\}, \]}
where $\ms$ represents a  multi-index 
and $\partial^\ms v$ indicates the partial derivatives of $v$ of order $|\ms|$. The associated norm on $W^{m,r}(\mathfrak{D})$ is denoted by $\|\cdot\|_{m,r,\mathfrak{D}}$ and the seminorm by $|\cdot|_{m,r,\mathfrak{D}}$. For $r=2$, $W^{m,2}(\mathfrak{D})$ is a Hilbert space $H^m(\mathfrak{D})$, and its corresponding norm and semi-norm are denoted by $\|\cdot\|_{m,\mathfrak{D}}$ and $|\cdot|_{m,\mathfrak{D}}$, respectively.
By $L^p(0,T;\mathcal{Y})$ we denote the space of measurable functions $v$ on $(0,T)$ with values in a Banach space $\mathcal{Y}$ satisfying
{\small\[\int_{0}^{T} \|v(t)\|_{\mathcal{Y}}^{p} dt < \infty ~\mbox{for}~ 1<{ p} < \infty ~~\mbox{and}~~\esssup_{0\leq t\leq T } \|v(t)\|_\mathcal{Y}<\infty \,\mbox{for}\, { p}=\infty.\]}
\noindent To apply the LDG method, we consider the vectors $\bv,\, \bw \in \bL^2(\mathfrak{D})$, and the tensors $ \,\underline{\bs}, \,\underline{\btau} \in \underline{\bL}^2(\mathfrak{D})$, such that $ \bv=(v_i)_{i=1}^d,\, \bw=(w_i)_{i=1}^d, \,\text{and} \: \underline{\bs}=[\sigma_{ij}]_{1\leq i,\, j\leq d},\, \underline{\btau}=[\tau_{ij}]_{1\leq i,\, j\leq d}$. We use the standard notations $(\nabla\bu)_{ij}=\partial_ju_i \,\,~~ \mbox{and}~~ (\nabla\cdot\underline{\bs})_i =\sum_{j=1}^{d}\partial_j\sigma_{ij}.$ 
We also denote by $\bu\otimes\bv$ the matrix whose $ij^\text{th}$ component is $u_iv_j$ and write 
\vspace{-2mm}
{\scriptsize\[ \underline{\bs}: \underline{\btau}:=\sum_{i,j=1}^{d}{ \sigma}_{ij}\tau_{ij} \,\,\,\mbox{and}~~ \bv\cdot\underline{\bs}\cdot\bw:=\sum_{i,j=1}^{d} v_i\sigma_{ij}w_j=\underline{\bs}:\bv\otimes\bw.\]}\normalsize
\noindent To formulate (\ref{NSEDARCY1})-(\ref{NSEDARCY3}) in a mixed setup, we set $\underline{\bs} =\textbf{D}(\bu)$. In a similar way, we will also proceed with the second and third interface conditions, given by (\ref{NSEDARCY11})-(\ref{NSEDARCY12}). Next, we introduce the functional spaces required for the weak solution $(\underline{\bs},\bu,p,\bomega,\phi)$, which are
{\scriptsize\begin{align*}
&\bX = \{\bv\in { \bH^1(\Omega_f)}: \bv = 0 ~\text{on}~{\Gamma}_f\}{\xb ,} \quad
\underline{\bY} =\underline{\bL}^2(\Omega_f),\quad \textbf{W}=\bL^2(\Omega_p),\quad \bV = \{\bv\in \bX: \nabla\cdot\bv = 0 ~\text{in}~{\Omega}_f\}, \\
&  M_p = \{\psi\in H^1(\Omega_p): \psi = 0 ~\text{on}~\Gamma_p\}, \,\, M_f = L_0^2(\Omega_f)=\Big\{q\in L^2(\Omega_f):\int_{\Omega_f}q\,d\Omega_f=0\Big\}
.
\end{align*}}%
\noindent Now we are ready to state the mixed variational formulation: Find  $(\underline{\bs},\bu,p,\bomega,\phi)\in \underline{\bY}\times\bX\times M_f\times\textbf{W}\times M_p$, such that for all $(\underline{\btau},\bv, q,\bw,\psi)\in \underline{\bY}\times\bX\times M_f\times\textbf{W}\times M_p$ and for $t>0$,
{\scriptsize
\begin{align*}
(\underline{\bs},\underline{\btau})_{\Omega_f} &= (\textbf{D}(\bu),\underline{\btau}), \\
(\bu_t,\bv)_{\Omega_f}+2\mu (\underline{\bs},\textbf{D}(\bv))_{\Omega_f}+ (\bu \cdot \nabla \bu, \bv)_{\Omega_f}   
- (p , \nabla \cdot \bv)_{\Omega_f}
&+\sum_{j=1}^{d-1}\frac{1}{G^j}(\bu\cdot \btau_{\Gamma}^j, \bv\cdot \btau_{\Gamma}^j)_{\Gamma}
+  g(\phi,\bv\cdot \bn_\Gamma)_{\Gamma}
=({\textbf{f}}_f,\bv)_{\Omega_f},\\
(\nabla \cdot\bu , q)_{\Omega_f} &= 0,\\
\bK^{-1}(\bomega,{\bw})_{\Omega_p}&=-(\nabla\phi,{\bw})_{\Omega_p},\\
S_0 g(\phi_t,\psi)_{\Omega_p}-  g( \bomega , \nabla\psi)_{\Omega_p} - g(\bu\cdot\bn_\Gamma , \psi)_{\Gamma} &= g(f_p,\psi)_{\Omega_p} , \\
(\bu(0),\bv)_{\Omega_f} = (\bu_0,\bv)_{\Omega_f},\,\,
(\phi(0),\psi)_{\Omega_p} &= (\phi_0,\psi)_{\Omega_p}.
\end{align*}}\normalsize
\noindent Under the assumptions $\bu_0\in \bL^2(\Omega_f), ~\phi_0\in L^2(\Omega_p), ~\textbf{f}_f\in L^2(0,T;\bL^2(\Omega_f))$, $f_p\in L^2(0,T;L^2(\Omega_p))$, $\bK\in L^\infty(\Omega_p)^{d^2}$, one can show the existence and uniqueness of the weak solution for the coupled NSD system (see \cite[Theorem 2]{DO17}).

\noindent Next, we recall some standard inequalities such as Korn,  Poincar{\`e}, and Sobolev's inequalities \cite[pp. 2057]{GR09}:
there exist constants $d_1$, $d_2$, $d_3>0$ that depend only on $\Omega_f$, such that $\forall{\bz}\in \bX$,
\begin{align}
| \bz|_{H^1(\Omega_f)}  \leq d_1\|\textbf{D}(\bz)\|_{L^2(\Omega_f)},\:
\| \bz\|_{L^p(\Omega_f)}  \leq d_2|\bz|_{H^1(\Omega_f)},(p=2,3,4), \: \| \bz\|_{L^r(\Gamma)}\leq d_3 |\bz|_{H^1(\Omega_f)},\,(r=2,4).\label{ctsineq}
\end{align}
\section{Numerical Scheme}\label{s4}
\noindent In order to discretize the coupled problem, we assume a regular family of triangulations $\mathcal{E}_{h}$ of $\bar{\Omega}$, consisting of triangles or tetrahedra {$E$}, depending on the dimension. Here, $h$ represents the maximum element diameter, and $h_E$ indicates the diameter of each element $E$.  
We restrict $\mathcal{E}_{h}$ to the subdomains $\Omega_f$ and $\Omega_p$ to obtain two regular families of triangulations $\mathcal{E}_{ih}$ for $i=f,p$, provided that the interface $\Gamma$ consists of a finite union of triangle ({tetrahedra}) edges. 
By regular, we mean that there exists a constant  $\beta>0$, independent of $h$, such that
\small{
\begin{equation*}
\frac{h_E}{\varrho_E}\leq \beta, \quad \forall E\in \mathcal{E}_h,\label{regular}
\end{equation*}}
\noindent where $\varrho_E$ denotes the diameter of the circle inscribed in $E$. We denote by $\Gamma_{ih}$, the set of all edges of $\mathcal{E}_{ih}$ interior to $\Omega_i$ for $i=f,\,p$. With each interior edge $e$, we associate a unit normal vector $\bn_e$. For a boundary edge $e\in \Gamma_f\cup\Gamma_p$, the unit vector $\bn_e$ coincides with the outward normal vector. For $e\in\Gamma$, we set $\bn_e=\bn_{\Gamma}$ such that  $\bn_e$ is the outward normal from $\Omega_f$ to $\Omega_p$.

\noindent Now, we introduce notations related to traces. Let $E^+$ and $E^-$ be two adjacent elements of $\mathcal{E}_{ih}$. Let $\bx$ be an arbitrary point of the interior edge $e=\partial E^+\cap\partial E^-\in \Gamma_{ih}$, and $\bn^+$ and $\bn^-$ be the corresponding outward normal vectors at that point.
Let $(\underline{\bs},\bu,p,\bomega,\phi)$ 
be a smooth function inside each element $E^{\pm}$. Also, let $(\underline{\bs}^{\pm},\bu^{\pm},p^{\pm},\bomega^{\pm},\phi^{\pm})$ denote the traces of $(\underline{\bs},\bu,p,\bomega,\phi)$ on $e$ from the interior of $E^{\pm}$. Then the averages $\{\cdot\}$ and jumps $[\cdot]$ are defined at $\bx\in e$ as 
{\scriptsize\begin{align*}
 \{\underline{\bs}\} &:=\frac{1}{2}(\underline{\bs}^++\underline{\bs}^-),\,\,  \{\bomega\}:= \frac{1}{2}(\bomega^+ +\bomega^-),\,\,
\{p\} := \frac{1}{2}(p^+ +p^-), \,\, 
 \{\bu\}:=\frac{1}{2}(\bu^+ +\bu^-), \,\,
 \{\phi\}:=\frac{1}{2}(\phi^+ + \phi^-), \\
{ \bsq\underline{\bs}\bqs} &:=\underline{\bs}^+\cdot\bn^+ +\underline{\bs}^-\cdot\bn^-,  \,\, [\bomega] := \bomega^+\cdot\bn^+ +\bomega^-\cdot\bn^-,\,\,
 \bsq p\bqs := p^+\bn^+ + p^-\bn^-,  [\bu]:=\bu^+\cdot\bn^+ + \bu^-\cdot\bn^-, \,\, \bsq \phi\bqs:=\phi^+\bn^+ +\phi^-\bn^-.
\end{align*}}\normalsize
For an edge {$e\in\Gamma_f\cup\Gamma_p$}, the averages of $\underline{\bs},\,\bu,\,p,\,\bomega,\,\phi$ on $e$ coincide with the value of $\underline{\bs},\,\bu,\,p,\,\bomega,\,\phi$ on $e$, and jumps are $\bsq p\bqs=p\cdot\bn,\,\,[\bu]=\bu\cdot\bn,\,\,[\underline{\bs}] =\underline{\bs}\cdot\bn$, with $\bn$ denoting the outward unit normal vector on the boundary $\Gamma_f$ or $\Gamma_p$. The jumps of $\bsq p \bqs$, $[\underline{\bs}]$ are vectors, and $[\bu]$ is a scalar. We now define jump of the velocity $\bu$, which is a matrix, given by
\small{\[\underline{\bsq\bu\bqs}:=\bu^+\otimes\bn^+ + \bu^-\otimes\bn^-.\]} Also, for {$e\in\Gamma_f$}, 
 the jump is $\underline{\bsq\bu\bqs}:=\bu\otimes\bn.$
 
\noindent We next introduce 
the following discontinuous spaces to approximate the stress tensor, NS velocity, pressure, velocity in porous media, and the piezometric head: 
{\small\begin{align*}
& \underline{\bY}_h = \{ \underline{\bs}_h\in  \underline{\bL}^2(\Omega_f):\, \underline{\bs}_h|_{E}\in (\mathbb{P}_{k_1-1}(E))^{d^2},\quad\forall E\in \mathcal{E}_{fh}\},\\
& \bX_h = \{\bv_h\in \bL^2(\Omega_f):\, \bv_h|_{E}\in (\mathbb{P}_{k_1}(E))^d{,}\quad\forall E\in \mathcal{E}_{fh}\},\\
& M_{fh} = \{q_h\in L^2(\Omega_f):\, q_h|_{E}{\,\in \mathbb{P}_{k_1-1}}(E){,}\quad\forall E\in \mathcal{E}_{fh}\},\\
&  M_{ph} = \{\psi_h\in L^2(\Omega_p):\, \psi_h|_{E}{\,\in \mathbb{P}_{k_2}}(E){,}\quad\forall E\in \mathcal{E}_{ph}\},\\
&  \textbf{W}_h = \{\bw_h\in \bL^2(\Omega_p):\, \bw_h|_{E}\in (\mathbb{P}_{k_2-1}(E))^d{,}\quad \forall E\in \mathcal{E}_{ph}\},
 \end{align*}}%
\noindent where $\mathbb{P}_k(E)$ represents the polynomial spaces of degree $k$ on $E$ with two positive integers $k_1\geq 1$ and $k_2\geq 1$.

\noindent To develop an LDG scheme for the coupled NSD problem, our idea is to begin with discretization of the Oseen equations in place of the NSEs, where the nonlinear term $(\bu\cdot\nabla)\bu$ will be replaced by the linear convective term $(\Z\cdot\nabla)\bu$. The convective velocity $\Z$ would be a suitable projection of the approximate velocity $\bu_h$ to the space $\textbf{J}(\mathcal{E}_{fh})$ to get back the NSE setup, where
\begin{equation*}
\scalebox{0.95}{${J}(\mathcal{E}_{fh})= \{\bv\in \bL^2(\Omega_f):\nabla\cdot\bv=0 ~~\mbox{and}~~ \bv|_{E}\in (H^1(E))^d \:\forall E\in \mathcal{E}_{fh} \}.$}
\end{equation*} 
\noindent For time discretization, we use the backward Euler method. Let $\{t_i\}_{i=0}^N$ be a uniform partition of the time interval $[0,\,T]$ with the time step $\Delta t>0$ and {$t_N=T$} as the final time. Therefore, $t_i=i\Delta t$, $0\leq i\leq N$. Further, for a given sequence $\{\nu^i\}_{i\geq 0}$, we define ${ \partial_t}\nu^i$ as follows:
\begin{equation*}
\scalebox{0.95}{$\partial_t\nu^i=\frac{1}{\Delta t}(\nu^i-\nu^{i-1}),\,\,i\geq 1.$}
\end{equation*} 
\noindent 
 We now define the elementwise formulation for the approximate solution $(\underline{\bSig}_h^i,\bU_h^i,P_h^i,\bQ_h^i,\Phi_h^i)\in \underline{\bY}_h\times \bX_h\times M_{fh}\times \bW_h\times M_{ph}$ such that for all $E\in\mathcal{E}_{fh}$ and $E\in\mathcal{E}_{ph}$,
  {\scriptsize
\begin{align}
    \int_{E} \underline{\bSig}^i_h:\underline{\btau}_h d\bx =  \int_{E} \bU_h^i\cdot\nabla\cdot{\underline{\btau}_h}  d\bx + \int_{\partial E} \hat{\bU}^i_{h,\sigma}&\cdot\underline{\btau}_h\cdot\bn_E\,ds,\label{elementwise1}\\
     \int_{E} \partial_t\bU_{h}^i\cdot\bv_h\, d\bx+2\mu\int_{E} \underline{\bSig}^i_h:\textbf{D}(\bv_h) \, d\bx -2\mu\int_{\partial E} \hat{\underline{\bSig}}^i_h :(\bv_h\otimes\bn_E) ds- \int_{E} P^i_h&\nabla\cdot\bv_h  d\bx
   \nonumber\\
   + \int_{\partial E} \hat{P}^i_h\bv_h\cdot\bn_E ds
     - \int_{E}\bU^i_h\nabla\cdot(\bv_h\otimes\Z^i)\,d\bx + \int_{\partial E}\Z^i\cdot\bn_E\hat{\bU}^i_{h,\Z}\cdot\bv_h\,ds 
&=\int_{E} \f_f^i\bv_h  d\bx,\label{elementwise2}\\
     -\int_{E} \bU^i_h\cdot\nabla q_h\, d\bx +  \int_{\partial E} \hat{\bU}^i_{h,p}\cdot\bn_E q_h\,ds &= 0,\label{elementwise3}\\
     \bK^{-1} \int_{E} \bQ^i_h\cdot{\bw}_h\, d\bx =  \int_{E} \Phi^i_h\nabla\cdot{\bw}_h\, d\bx - \int_{\partial E} &\hat{\Phi}^i_{h}{\bw}_h\cdot\bn_E \,ds,\label{elementwise4}\\
     gS_0 \int_{E} \partial_t\Phi^i_{h}\cdot\psi_h \, d\bx - g\int_{E}\bQ^i_h\cdot\nabla\psi_h\, d\bx + g\int_{\partial E} \psi_h\hat{\bQ}^i_h\cdot\bn_E\,ds &= g\int_{E} f_p^i\cdot\psi_h \, dx, \label{elementwise5}
     \end{align}}
     for all  $(\underline{\btau}_h,\bv_h,q_h,{\bw}_h,\psi_h)\in \underline{\bY}_h\times \bX_h\times M_{fh}\times {\bW}_h\times M_{ph}$.\par
     \vspace{0.25cm}
     \noindent Here 
     $\bn_E$ denotes the unit outward normal on the boundary $\partial E$ of the element $E$. Also, it is worth mentioning that $ \hat{\underline{\bSig}}^i_h$, $\hat{\bU}^i_{h,\sigma}$, $\hat{\bU}^i_{h,\Z}$, $\hat{\bU}^i_{h,p}$, $\hat{P}^i_h$, $\hat{\bQ}^i_h$, and $\hat{\Phi}^i_{h}$ are the numerical fluxes that will be defined in the forthcoming discussion. 
     \begin{Remark}\label{locally conservative}
 Observe that if one chooses the test function $\bv_h=1$ in (\ref{elementwise2}) on $E$ belonging to $\mathcal{E}_{fh}$, excluding the elements having the edge lying on $\Gamma$ and zero otherwise, for each element we obtain
\small{\begin{align*}
     \int_E \partial_t\bU_{h}^i \, d\bx + \int_{\partial E} (-2\mu \underline{\hat{\bSig}}_h^i \cdot\bn_E + (\Z^i\cdot\bn_E)\hat{\bU}_{h,\bw }^i + \hat{P}_h^i\bn_E)\,ds  = \int_E\textbf{f}^i_f\, d\bx.
     \end{align*}}
     In addition, if we take the test function $\psi_h=1$ in (\ref{elementwise5}) on $E$  belonging to $\mathcal{E}_{ph}$, excluding the elements having the edge lying on $\Gamma$ and zero
otherwise, for each element we also have
   \small{ \begin{align*}
    gS_0 \int_{E} \partial_t\Phi_{h}^i \, d\bx + g\int_{\partial E} \hat{\bQ}_h^i\cdot\bn_E\,ds &= g\int_{E} f^i_p \, dx.
\end{align*}}
\end{Remark}
\noindent The structure of the LDG method will be completed after defining the numerical fluxes. It is worth mentioning that the fluxes have to be properly designed to ensure the stability of the method. We begin by defining the numerical fluxes {$\hat{\underline{\bSig}}_h^i$, $\hat{\bU}^i_{h,\sigma}$, $\hat{\bQ}^i_h$, and $\hat{\Phi}_h^i$} associated with the Laplacian (see \cite{CCPS00}). On a face $e$ inside the domain $\Omega_f$, we take
{\scriptsize\begin{align*}
\hat{\underline{\bSig}}_h^i &= \{\underline{\bSig}^i_{ h}\}- C_{11}\underline{\bsq\bU_h^i\bqs} - \bsq\underline{\bSig}_h^i\bqs\otimes\textbf{C}_{12},\quad
        \hat{\bU}^i_{h,\sigma} = \{\bU^i_h\} + \underline{\bsq \bU^i_h\bqs}\cdot\textbf{C}_{12},   \\     
        \hat{\bQ}^i_h &= \{\bQ^i_h\} + C_{11}\bsq\Phi^i_h\bqs-\textbf{C}_{12}[\bQ^i_h],\quad
        \hat{\Phi}_h^i= \{\Phi^i_h\} + \textbf{C}_{12}\cdot\bsq\Phi^i_h\bqs,
     \end{align*}}
and if e lies on $\partial\Omega_f$,
{\scriptsize\begin{equation*}
  \hat{\underline{\bSig}}_h^i = {\underline{\bSig}}_h^{i+}  - C_{11}{\bU}_h^{i+}\otimes\bn^+,\quad
       \hat{\bU}^i_{h,\sigma}= \textbf{0},  \quad
       \hat{\bQ}^i_h = \bQ_h^{i+}+C_{11}\Phi_h^{i+}\bn^+,\quad
       \hat{\Phi}^i_h= 0. 
     \end{equation*} }
    When $e$ lies on interface $\Gamma${,}   
    {\scriptsize\begin{equation*}
       \hat{\underline{\bSig}}_h^i = {\underline{\bSig}_h}^{i+} ,\quad
       \hat{\bU}^i_{h,\sigma} = \bU_h^{i+},\quad
       \hat{\bU}^i_{h,\Z }=\bU_h^{i+}, \quad
       \hat{\bQ}_h^i= \bQ_h^{i-},\quad
       \hat{\Phi}_h^i= \Phi_h^{i-}.
     \end{equation*} }\normalsize
\noindent The parameters $C_{11}>0$ and $\textbf{C}_{12}\in\mathbb{R}^d$ are defined for each edge $e$ and have to be chosen properly (see \cite{CCPS00}).

\noindent The numerical fluxes associated with the incompressibility constraint are given by $\hat{\bU}^i_{h,p}$ and $\hat{P}^i_h$ (for instance, \cite{CKS04}). If the face $e$ is inside $\Omega_f${,} we take
\small{\begin{equation*}
{\bU}^i_{h,p} = \{\bU^i_h\} + D_{11}\bsq P^i_h\bqs + \textbf{D}_{12}[\bU^i_h],\,\, \,\,   \hat{P^i_h} = \{P^i_h\}-\textbf{D}_{12}\cdot\bsq P^i_h\bqs,
\end{equation*} }
and if e lies on $\partial\Omega_f${,}
\small{\begin{equation*}
\hat{\bU}^i_{h,p} = \textbf{0}
,\,\, \hat{P}^i_h = P_h^{i+}.
\end{equation*} }
When $e$ lies on the interface $\Gamma$,
\small{\begin{equation*}\hat{\bU}^i_{h,p} = \bU_h^{i+},\,\, \hat{P}^i_h = P_h^{i+}.
\end{equation*} }
\noindent Here, $D_{11}$ and $\textbf{D}_{12}\in\mathbb{R}^d$ are the parameters depending on $\bx\in e$ (see \cite{CKS04}). In our analysis, we consider $D_{11}=0$ (see \cite[pp. 580]{CKS04}). 
For the convective numerical flux, we take the standard upwind flux \cite{CKS04} defined as
\small{\begin{equation*}
 \hat{\bU}^i _{h ,\Z}(\bx) = 
{ \displaystyle\lim_{\epsilon \to 0}} \bU^i_h (\bx - \epsilon \Z(\bx)). 
\end{equation*} }%
\noindent This completes the definition of the LDG method for the Oseen-Darcy problem.
For the sake of completeness of the definition of the LDG method for the NSD equations (\ref{NSEDARCY1})-{(\ref{NSEDARCY12})}, it only remains to introduce what we refer to as the post-processing operator $\p$ and to plug $\Z^i=\p\bU^i_h$ into the approximation (\ref{elementwise1})-(\ref{elementwise5}).
For a piecewise smooth velocity field $\bu$, the operator $\p$ (see \cite{CKS05}) is defined as
\small{\begin{equation*}
\p\bu|_E=\p_E(\bu|_E,\hat{\bu}_{p}),\,\, \forall E\in\mathcal{E}_{fh},
\end{equation*} }
\noindent where {$\hat{\bu}_{p}$} is the flux related to the incompressibility constraint. For further details, one can refer to \cite[Section 3.3]{CKS05}.

\noindent In preparation for the forthcoming analysis, we define the following bilinear forms:
 $a_{NS}: \underline{\bY}_h\times \underline{\bY}_h\rightarrow \mathbb{R}$,  $b_{NS}:\textbf{X}_h\times \underline{\bY}_h\rightarrow \mathbb{R}$, $c_{NS}:\textbf{X}_h\times \textbf{X}_h\rightarrow \mathbb{R}$, $d_{NS}:\textbf{X}_h\times \textbf{X}_h\rightarrow  \mathbb{R}$,   $a_{D}:\textbf{W}_h\times \textbf{W}_h\rightarrow \mathbb{R}$, $b_{D}:{M_{ph}}\times \textbf{W}_h\rightarrow \mathbb{R}$ and $c_{D}:{M_{ph}}\times {M_{ph}}\rightarrow \mathbb{R}$, such that 
{\small\begin{align*}
a_{NS}( \underline{\bs}_h, \underline{\btau}_h)& = \int_{\Omega_f} \underline{\bs}_h:\underline{\btau}_h  d\bx,\\
b_{NS}(\bu_h, \underline{\btau}_h)& = 2\mu\sum_{E\in \mathcal{E}_{fh}} \int_{E} \bu_h\cdot\nabla\cdot \underline{\btau}_h  d\bx-2\mu\sum_{e\in \Gamma_{fh}}\int_{e}(\{\bu_h\}+\underline{\bsq\bu_h\bqs}\cdot\textbf{C}_{12})\cdot\bsq \underline{\btau}_h\bqs\, ds-2\mu\int_{\Gamma}\bu_h\cdot  \underline{\btau}_h\cdot\bn\, ds,\\
c_{NS}(\bu_h,\bv_h)& = 2\mu \sum_{e\in \Gamma_{fh}}\int_{e} \frac{C_{11}}{|e|}\underline{\bsq\bu_h\bqs}:\underline{\bsq\bv_h\bqs}\,ds +2\mu  \sum_{e\in \Gamma_f}\int_{e}\frac{C_{11}}{|e|} (\bu_h\otimes\bn):(\bv_h\otimes\bn)\,ds ,\\
d_{NS}(p_h,\bv_h)& = -\sum_{E\in \mathcal{E}_{fh}} \int_{E}p_h\nabla\cdot\bv_h\, d\bx + \sum_{e\in \Gamma_{fh}}\int_{e}(\{p_h\}-\textbf{D}_{12}\cdot\bsq p_h\bqs)[\bv_h]\, d\bx + \int_{\Gamma_f} p_h\bv_h\cdot\bn \, ds,\\
a_{D}(\bomega_h,\bw_h) &=\bK^{-1}\int_{\Omega_p} \bomega_h\cdot \bw_h  d\bx,\\
 b_{D}(\phi_h,\bw_h) &= g\sum_{E\in \mathcal{E}_{ph}}\int_{E}\phi_h\nabla\cdot\bw_h\, d\bx - g\sum_{e\in \Gamma_{fh}}\int_{e}(\{\phi_h\}+\textbf{C}_{12}\cdot\bsq\phi_h\bqs)[\bw_h]\,d\bx-g\int_{\Gamma}\phi_h\bw_h\cdot\bn\,ds,\\
 c_{D}(\phi_h,\psi_h) &= g\sum_{e\in \Gamma_{fh}}\int_{e} \frac{C_{11}}{|e|}\bsq\phi_h\bqs\cdot\bsq\psi_h\bqs\,ds + g\sum_{e\in \Gamma_p}\int_{e} \frac{C_{11}}{|e|}\phi_h\psi_h\,ds.
\end{align*}}
\noindent Also, using integration by parts, the forms $b_{NS} $ and $b_{D}$ can be written as
{\scriptsize\begin{align*}
 b_{NS}(\bu_h,\underline{\btau}_h)& = -2\mu\sum_{E\in \mathcal{E}_{fh}} \int_{E} \textbf{D}(\bu_h):\underline{\btau}_h { d\bx} +2\mu \int_{\Gamma_h}(\{{ \underline{\btau}_h}\}-\bsq \underline{\btau}_h\bqs { \otimes}\textbf{C}_{12}):\underline{\bsq{\bu_h}\bqs}\, ds \nonumber+ 2\mu\int_{\Gamma_f}\underline{\btau}_h:(\bu_h\otimes\bn)\, ds, \\
 b_{D}(\phi_h,\bw_h) &= -g\sum_{E\in \mathcal{E}_{ph}}\int_{E}\bw_h\cdot\nabla \phi_h\,{ d\bx} +g \int_{\Gamma_h}(\{\bw_h\}-\textbf{C}_{12}\cdot[\bw_h])\cdot\bsq\phi_h\bqs\, ds+g\int_{\Gamma_p}\phi_h\bw_h\cdot\bn\,ds.
\end{align*}} 
\noindent Next, we define the seminorm (\cite[pp. 328]{CKSS02}) and energy norm in discrete spaces $\underline{\bY}_h$, $\bX_h$, $\bW_h$, and ${M_{ph}}$ as follows:
{\small\begin{align}
&|( \underline{\bs}_h,\bv_h)|^2_{\mathcal{A}} =a_{NS}( \underline{\bs}_h\: \underline{\bs}_h)+c_{NS}(\bv_h,\bv_h){,}\,\,\,\,\forall   \underline{\bs}_h\in  \underline{\bY}_h,\, \forall \bv_h\in\bX_h,\label{norm1}\\
&|( \bomega_h,\psi_h)|^2_{\mathcal{B}} =a_{D}(\bomega_h,\bomega_h)+c_{D}(\psi_h,\psi_h){,}\,\,\,\,\forall  \bomega_h\in {\bW_h}, \forall \psi_h\in M_{ph},\label{norm2}\\
&\|\bv_h\|^2_{\bX_h}=\bigg(2\sum_{E\in\mathcal{E}_{fh}}\|\textbf{D} (\bv_h)\|^2_{L^2(E)}+c_{NS}(\bv_h,\bv_h)\bigg)^{1/2},\,\forall\bv_h\in \bX_h,\label{norm3}\\
&\|\psi_h\|^2_{M_{ph}}=\bigg(\sum_{E\in\mathcal{E}_{ph}}\|\bK^{1/2} \nabla \psi_h\|^2_{L^2(E)}+c_{D}(\bv_h,\bv_h)\bigg)^{1/2}{,}\,\forall \psi_h\in {M_{ph}}.\label{norm4}
\end{align}}
\noindent For $M_{fh}$, we take the $L^2$-induced norm.
In order to define the DG discretization of the trilinear term, we write the form $O_{NS}(\cdot,\cdot,\cdot)$ (see \cite{CKS05}) as
\begin{align}
O_{NS}(\Z,\bu_h,\bv_h) = -\sum_{E\in \mathcal{E}_{fh}} \int_{E} \bu_h\cdot\nabla\cdot(\bv_h\otimes\Z){ d\bx} + \sum_{E\in \mathcal{E}_{fh}} \int_{\partial E \setminus \Gamma_f^-} \Z\cdot\bn_E\,\hat{\bu}_{h,\Z}\cdot\bv_h\, ds.\label{trilinear}
\end{align} 
\noindent where \small{\begin{equation*}
\Gamma_f^-=\{\bx\in\Gamma_f:\Z(\bx)\cdot\bn(\bx)<\textbf{0}\},\end{equation*} }
 represents the inflow part of $\Gamma_f$.
\noindent Now, we present the fully discrete LDG formulation: For $1\leq i\leq N$, given functions  $\underline{\bSig}_h^0$, $\bU_h^0$, $\bQ_h^0$, and $\Phi_h^0$, find $({\underline{\bSig}_h^i},\bU_h^i,P_h^i,\bQ_h^i,\Phi_h^i)\in \underline{\bY}_h\times\bX_h\times M_{fh}\times{\bW}_h\times M_{ph}$, such that{,} for all $(\underline{\btau}_h,\bv_h,q_h,\bw_h,\psi_h)\in \underline{\bY}_h\times \bX_h\times M_{fh}\times {\bW}_h\times  M_{ph}$,
\small{\begin{align}
  a_{NS}({ \underline{\bSig}_h^i},{ \underline{\btau}_h}) +  b_{NS}(\bU_h^i,{ \underline{\btau}_h}) &= 0,\label{FSLNSD1}\\ 
  ( \partial_t\bU_h^i,\bv_h)- b_{NS}(\bv_h,{ \underline{\bSig}_h^i}) +  c_{NS}(\bU_h^i,\bv_h) + d_{NS}(&\bv_h,P_h^i) \nonumber\\ 
  + O_{NS}(\mathfrak{P}_h\bU^{i-1}_h,\bU_h^i,\bv_h)
  + \sum_{j=1}^{d-1}\frac{1}{G^j}(\bU_h^i\cdot\btau_{\Gamma},\bv_h\cdot\btau_{\Gamma})_{{\Gamma}}
+  g(\Phi_h^i,\bv_h\cdot{\bn_{\Gamma}})_{{\Gamma}} & = (\textbf{f}_f^i,\bv_h)_{{\Omega}_f},\label{FSLNSD2}\\
d_{NS}(\bU_h^i,q_h)&=0,\label{FSLNSD3}\\ 
a_{D}(\bQ_h^i,\bw_h)-b_{D}(\Phi_h^i,\bw_h)&=0,\label{FSLNSD4}\\
gS_0(\partial_t\Phi_{h}^i,\psi_h)+  b_{D}(\psi_h,\bQ_h^i)+ c_{D}(\Phi_h^i,\psi_h)-g(\bU_h^i\cdot\bn_{\Gamma},\psi_h)_{{\Gamma}}&= g(f_p^i,\psi_h)_{{\Omega}_p}.\label{FSLNSD5}
\end{align}}
\noindent Note that $\bU_h^0$ and $\Phi_h^0$ are suitable approximations of $\bu_0$ and $\phi_0$, respectively. 
\begin{Remark}\label{locally conservative1}
Applying Remark \ref{locally conservative}, we observe that the scheme (\ref{FSLNSD1})-(\ref{FSLNSD5}) is locally conservative.
\end{Remark}
\noindent The next lemma gives the consistency of the scheme (\ref{FSLNSD1})-(\ref{FSLNSD5}). 
\begin{lemma}
Let $(\underline{\bs},\bu,p,\bomega,\phi)$ be a solution to (\ref{NSEDARCY1})-(\ref{NSEDARCY12}). Then $(\underline{\bs},\bu,p,\bomega,\phi)$ satisfies (\ref{FSLNSD1})-(\ref{FSLNSD5}).
\end{lemma}
\noindent Utilizing the idea that the numerical fluxes are consistent (see \cite[pp. 325]{CKSS02}), i.e.
\small{\begin{equation*} \hat{\underline{\bSig}}^i_h=\bs,\,\, \hat{\bU}^i_{h,\sigma}=\bu,\,\, \hat{\bU}^i_{h,\Z}=\bu,\,\, \hat{\bU}^i_{h,p}=\bu,\,\, \hat{P}^i_h=p,\,\, \hat{\bQ}^i_h=\bomega,\,\, \hat{\Phi}^i_{h}=\phi,
\end{equation*} }
and following similar footsteps as in \cite[Section 3.3]{ABCM02}, we can prove the above lemma.

\noindent We also define the discrete divergence-free space, which is analogous to $\bV$, as follows:
\small{\begin{equation*}\bZ_h=\{\bv_h\in\bX_h:\forall { q_h}\in M_{fh},\,d_{NS}({q_h},\bv_h)=0\}.
\end{equation*} }
\noindent Considering the problem (\ref{FSLNSD1})-(\ref{FSLNSD5})~in $\textbf{Z}_h$, we obtain the equivalent LDG formulation for the system as follows: For all $(\underline{\btau}_h,\bv_h,\bw_h,\psi_h)\in \underline{\bY}_h\times \bZ_h\times {\bW}_h\times  M_{ph}$, find $({\underline{\bSig}_h^i},\bU_h^i,\bQ_h^i,\Phi_h^i)_{i\geq 1}\in \underline{\bY}_h\times\bZ_h\times{\bW}_h\times M_{ph}$, such that
\small{\begin{align}
  a_{NS}({ \underline{\bSig}_h^i},{ \underline{\btau}_h}) +  b_{NS}(\bU_h^i,{ \underline{\btau}_h}) &= 0,\label{DFSLNSD1}\\ 
  ( \partial_t\bU_h^i,\bv_h)- b_{NS}(\bv_h,{ \underline{\bSig}_h^i}) +   c_{NS}(\bU_h^i,\bv_h)  
  + O_{NS}(\p\bU_h^{i-1},\bU_h^i,&\bv_h) 
  + \sum_{j=1}^{d-1}\frac{1}{G^j}(\bU_h^i\cdot\btau_{\Gamma}^j,\bv\cdot\btau_{\Gamma}^j)_{{\Gamma}}\nonumber\\
+  g(\Phi_h^i,\bv_h\cdot{\bn_{\Gamma}})_{{\Gamma}} & = (\textbf{f}_f^i,\bv_h)_{{\Omega}_f},\label{DFSLNSD2}\\
a_{D}(\bQ_h^i,\bw_h)-b_{D}(\Phi_h^i,\bw_h)&=0,\label{DFSLNSD3}\\
gS_0(\partial_t\Phi_{h},\psi_h)+  b_{D}(\psi_h,\bQ_h^i)+ c_{D}(\Phi_h^i,\psi_h)-g(\bU_h^i\cdot\bn_{\Gamma},\psi_h)_{{\Gamma}}&= g(f_p^i,\psi_h)_{{\Omega}_p},\label{DFSLNSD4}
\end{align}}
\noindent $\forall (\underline{\btau}_h,\bv_h,q_h,\bw_h,\psi_h)\in \underline{\bY}_h\times \bZ_h\times M_{fh}\times {\bW}_h\times  M_{ph}$.\par
\vspace*{2mm}
\noindent We now recall some approximation properties for the discrete spaces $\bX_h$, $M_{fh}$, and $M_{ph}$. We also note that henceforth the constant $C > 0$ that appears in the subsequent analysis is independent of $h$ and $E$.
\begin{lemma}\cite[Lemma 6.1]{CGPR17}
There exists an approximation operator ${\Pi}_h\in \mathcal{L}(\bX,\bX_h)$ satisfying
\small{\begin{align*}
d_{NS}({\Pi}_h(\bv)-\bv,\,q_h)=0\quad \forall \bv\in \bX, \,\,\forall q\in M_{fh},
\end{align*}}
and for all $E $ in ${\mathcal{E}}_{fh}$, for all $\bv$ in $\bX\cap W^{s,r}(E)^d,\quad 1\leq r \leq \infty, \,\, 1\leq s \leq k_1+1$,
\small{\begin{align*}
\|\bv-\Pi_h(\bv)\|_{L^r(E)} + h_E\|\nabla(\bv-\Pi_h(\bv))\|_{L^r(E)}&\leq Ch_E^s |\bv|_{W^{s,r}(\Delta_ E)},
\end{align*}}
where $\Delta_ E\subset \Omega_f$ is a macro-element used in the construction of $\Pi_h(\bv)$ in $E$. Also, for all $s, 1 \leq s \leq k_1 + 1$, $\Pi_h$ satisfies the following bound:
\small{\begin{align*}
\|\bv-\Pi_h(\bv)\|_{L^2(\Omega_f)} + h\|\bv-\Pi_h(\bv)\|_{\bX_h}&\leq Ch^s |\bv|_{s,\Omega_f}.
\end{align*}}
\end{lemma}
\noindent Moreover, there exists a local projection operator $r_h:L_0^2(\Omega_f)\rightarrow M_{fh}$ (see \cite[(3.10)-(3.11)]{GRW005}), defined for any $E\in \mathcal{E}_{fh}$ and $z_h\in \mathbb{P}_{r-1}(E)$,
\small{\begin{align}
\int_E z_h{(y-r_h y)}\,dx=0, \quad\forall y\in L_0^2{(\Omega_f)}.\label{projection1defn}
\end{align}}
\noindent Additionally, $r_h$ satisfies 
\small{\begin{align}
\| y-r_h( y)\|_{L^2(E)}\leq Ch^s_E | y|_{s,E},\quad \forall s\in [0,r],\,\, \forall  y\in H^s(\Omega_f)\cap L_0^2(\Omega_f).\label{projection1}
\end{align}}
\noindent To interpolate the piezometric head for the Darcy region, we use the Lagrange interpolant ${ P_h^\phi}$ in $\mathbb{P}_{k_2+1}$ (see \cite[pp. 38]{CGPR17}) satisfying the following approximation errors for all $r$, $1\leq r\leq \infty$ 
\small{\begin{align}
 E\in \mathcal{E}_{ph}, \, \forall\, s \in  [2,k_2+1], \, l=0,\,1,\, \forall \,\psi \in  W^{s,r}(E),\,|\psi-P_h^\phi\psi|_{W^{l,r}( E)}\leq Ch^{s-l} |\psi|_{W^{s,r}( E)}.\label{projection3}
\end{align}}
\noindent Now, we recall 
the $L^\mathfrak{p}$-estimate for functions in the broken-$H^1$ space $\bX_1$ $= \{\bv\in \bL^2(\Omega_f):\, \bv|_{E}\in \bH^1(E),\forall E\in\mathcal{E}_{fh}\}$ \cite[Lemma 6.2]{GRW005}: for each $\mathfrak{p}\in [2,\infty)$,
\small{\begin{align}
\|\bnu\|_{L^\mathfrak{p}(\Omega_i)} \leq C_\mathfrak{p}\|\bnu\|_{\bX_h},\label{brokenpoincare}
\end{align} }
\noindent presuming that $C_\mathfrak{p}$ is a positive constant.
In the case of $3D$, the result is valid for $2\leq\mathfrak{p}\leq 6$.
Next comes the standard trace inequalities. 
\begin{lemma}\cite[pp. $23$, $(2.1)$-$(2.2)$]{R08}, \cite[pp. $118$, $(7.7)$]{KJ98}\label{TL1}
For each $E\in \mathcal{E}_{fh}$, the following inequalities hold: \small{\begin{align}
&\|\nabla^i\bv\|_{0,e}  \leq C(h_E^{-1/2}\|\nabla^i\bv\|_{0,E}+h_E^{1/2}\|\nabla^{i+1} \bv\|_{0,E})\quad \forall e\in \partial E,\,\,\forall \bv\in \bH^{i+1}(E),\,i=0,1,\label{T1}\\
&\|\bv\|_{L^4(e)}  \leq C(h_E^{-1/4}\|\bv\|_{0,E}+h_E^{ (3-d)/4}\|\nabla \bv\|_{0,E})\quad \forall e\in \partial E,\,\, \bv\in \bH^1(E).\label{T3}
\end{align}}%
\end{lemma}
\noindent Similar to the trace inequalities defined in (\ref{T1})-(\ref{T3}), we have the following trace and inverse inequalities for discrete spaces as well.
\begin{lemma} \label{TL2} \cite[pp. 23, $(2.3)$-$(2.5)$]{R08},\cite[Section 1.4.3]{DDE11} For each $E\in \mathcal{E}_{fh}$ and $1\leq p \leq\infty $, the following inequalities hold:\vspace*{-5mm}
{\small\begin{align}
&{\|\nabla^i\bv_h\|_{0,e}}{\leq Ch_E^{-1/2}\|\nabla^i\bv_h\|_{0,E}}\quad {\forall e\in \partial E,} \,\,{\forall \bv_h\in\bX_h,\,i=0,1,}\label{T4}\\
&\|{\underline{\bs}_h}\|_{0,e}\leq Ch_E^{-1/2}\|{\underline{\bs}_h}\|_{0,E}\quad \forall e\in \partial E, \,\,\forall \underline{\bs}_h\in\underline{\bY}_h, \label{ST7}\\
&\|\bv_h\|_{L^p(e)}\leq C{h_E}^{-1/p}\|\bv_h\|_{L^p(E)}\,\,\forall e\in \partial E, \,\,\forall p\geq 1,\label{ST8}\\
&\|\nabla\bv_h\|_{0,E}\leq Ch_E^{-1}\|\bv_h\|_{0,E}\quad \forall \bv_h\in\bX_h,\label{T6}\\
&\|\bv_h\|_{L^p(E)}\leq Ch_E^{ d(1/p-1/2)}\|\bv_h\|_{0,E}\quad \forall \bv_h\in\bX_h.\label{T7}
\end{align}}
\end{lemma}
\noindent The above-mentioned inequalities also hold similarly for the discrete spaces $\bW_h$ and $M_{ph}$. 
The following estimates yield the 
usual Sobolev embeddings and trace inequalities to 
the discrete spaces (see \cite[(34)-(38)]{CGPR17}): 
\small{\begin{align}
&\text{For}\,\, 2\leq r\leq 6,\,\, \|\bv_h\|_{L^r(\Omega_f)} \leq C  \|\bv_h\|_{\bX_{h}},\quad \|\psi_h\|_{L^r(\Omega_p)} \leq C  \|\psi_h\|_{M_{ph}}, ~\forall \bv_h\in \bX_{h},\,\psi_h\in M_{ph}. \label{T8}\\
&\text{For}\,\, 2\leq r\leq 4,\,\, \|\bv_h\|_{L^r(\Gamma)} \leq C  \|\bv_h\|_{\bX_{h}},\,\,\,\quad\|\psi_h\|_{L^r({\Gamma})} \leq C  \|\psi_h\|_{M_{ph}}, ~\forall \bv_h\in \bX_{h},\, \psi_h\in M_{ph}.\label{T9}
\end{align}}
\noindent For $r=2$, the constants of (\ref{T8}) will be represented by $C_2$ and $C_1$, respectively.
Also, utilizing (\ref{T9}), the following term can be bounded:
\small{\begin{align}
 \frac{1}{2}(\bv_h\cdot\bn_{\gamma},\bw_h\cdot\bw_h)_{\Gamma}\leq C_0 \|\bv_h\|_{\bX_h} \|\bw_h\|^2_{\bX_h}\quad \forall \bw_h, \bv_h\in {\bX}_h.\label{nonliearity2}
\end{align}}
\noindent There is another crucial inequality, which is very helpful to handle the interface terms.
\begin{lemma}\cite[Theorem 5.2]{CGPR17})
 $\forall \psi_h\in M_{ph}, \bu_h\in\bZ_h $, the following estimate holds: 
\small{\begin{align}
|(\psi_h,\,\bu_h\cdot\bn_{\gamma})_{\Gamma}|\leq C\|\psi_h\|_{M_{ph}}\|\bu_h\|_{L^2(\Omega_f)}.\label{T12}
\end{align}}
\end{lemma}
\noindent The next lemma demonstrates the boundedness property of the bilinear forms $a_{NS}(\cdot, \cdot) + c_{NS}(\cdot, \cdot)$ and $a_{D}(\cdot, \cdot) + c_{D}(\cdot, \cdot)$.
\begin{lemma}\label{boundedness}
$\forall \underline{\bs}_h,\,\underline{\btau}_h\in \underline{\bY}_h$, $\bv_h,\,\bz_h\in\bX_h$, $\bomega_h,\,\bw_h \in \bW_h$, and $\psi_h, \,\theta_h \in M_{ph}$, the following inequalities hold:
\small{\begin{align*}
&|a_{NS}(\underline{\bs}_h,\underline{\btau}_h)+ c_{NS}(\bv_h,\bz_h)|\leq C |(\underline{\bs}_h,\bv_h )|_{\mathcal{A}} |(\underline{\btau}_h,\bz_h)|_{\mathcal{A}}, \,\,|a_{D}(\bomega_h,\bw_h)+ c_{D}(\psi_h,\theta_h)|\leq C |(\bomega_h,\psi_h)|_{\mathcal{B}} |(\bw_h,\theta_h)|_{\mathcal{B}}.
\end{align*}}
\end{lemma}
\noindent Applying Cauchy-Schwarz's and Young's inequalities, the estimates can be proved easily.\par
\noindent Also, we have the following inf-sup condition for the space $(\bX_h,M_{fh})$, which will be required for obtaining the pressure estimates.
\begin{lemma}\cite[Proposition 4.4]{CKS05}\label{infsup}
There exists a constant $\beta^{*}>0$, independent of $h$, such that
\small{\begin{align}
\inf_{q_h\in M_{fh}} \sup_{\bv_h\in {\bX}_h } \frac{b_{NS}(\bv_h,q_h)}{\|\bv_h\|_{{\bX}_h }\|q_h\|_{M_{fh}}}\geq \beta^{*}.
\end{align}}
\end{lemma}
\noindent A careful inspection reveals that the above inf-sup condition is also satisfied for the smaller space $\overline{\bX}_h$ (see \cite[pp. 1078]{CKS05} ) defined as
\[\overline{\bX}_h = \{\bv\in\bX_h:\nabla\cdot\bv\in \bL^2(\Omega_f)\}.\] 
\noindent Furthermore, we have the $L^4$-norm estimate for the elements of $\bX_h$.
\begin{lemma}\cite[Lemma 3.8]{GRW005S}
When $\Omega\subset \mathbb{R}^2$ is convex, the following estimate holds true: 
\begin{align}
\|\bv_h\|_{L^4(\Omega)}\leq C \|\bv_h\|^{1/2}_{L^2(\Omega)}\|\bv_h\|^{1/2}_{\bX_h}\quad \forall\bv_h\in\bX_h .\label{nonlinearity8}
\end{align}
\end{lemma}
\noindent Next, we observe that the linear operator $\p$ is bounded from $\bX_h$ to $\overline{\bX}_h$ with respect to the norm $\|\cdot\|_{\bX_h}$.
\begin{lemma}\cite[Proposition 4.6]{CKS05}
Let $\bv_h \in \bX_h$. Then, we have
\begin{align}
\|\p\bv_h\|_{\bX_h}\leq C_a \|\bv_h\|_{\bX_h}, \label{Postprocess1}
\end{align}
with a constant $C_a>0$ that is independent of $h$.
\end{lemma}
\noindent 
Since $\p$ produces functions in $\overline{\bX}_h$, we have $\p\bar{\bv}=\bar{\bv}$ and satisfy the following bound (see \cite[pp. 1091]{CKS05}):
\begin{align}
\inf_{\bar{\bv}\in \overline{\bX}_h}\|\bu-\bar{\bv}\|_{\bX_h} \leq C h^{k_1}\|\bu\|_{k_1+1,\Omega_f}.\label{Postprocess2}
\end{align}%
\noindent 
We now provide a relation between the $|\cdot|_{\mathcal{A}}${-} seminorm, $\|\cdot\|_{\bX_h}${-}norm and,  $|\cdot |_{\mathcal{B}}${-} seminorm and $\|\cdot\|_{M_{ph}}${-}norm.
\begin{lemma}\label{equivalence}
There are positive constants $K_\mathcal{A}$ and $K_\mathcal{B}$, independent of $h$, such that for all $\underline{\bs}_h\in \underline{\bY}_h$, $\bv_h\in\bX_h$, $\bw_h \in \bW_h$ and { $\psi_h\in M_{ph} $}, { we have}
\small{\begin{align}
&\|\bu_h\|^2_{\bX_h}\leq K_\mathcal{A}|({ \underline{\bs}_h},\bu_h)|^2_{\mathcal{A}} ,\quad |({ \underline{\bs}_h},\bu_h)|^2_{\mathcal{A}}\leq K_\mathcal{A}\|\bu_h\|^2_{\bX_h},\label{e1}\\
&\|\psi_h\|^2_{M_{ph}}\leq K_\mathcal{B}|(\bw_h,\psi_h)|^2_{\mathcal{B}} ,\quad |(\bw_h,\psi_h)|^2_{\mathcal{B}}\leq K_\mathcal{B}\|\psi_h\|^2_{M_{ph}}.\label{e2}
\end{align}}
\vspace*{-2mm}
\begin{proof}
From (\ref{DFSLNSD1}), we have 
\begin{align}
\sum_{E\in \mathcal{E}_{fh}} \int_{E} \textbf{D}(\bu_h):\underline{\btau}_h { d\bx} = \int_{\Omega_f} { \underline{\bs}_h}:{ \underline{\btau}_h}\, { d\bx}   + \int_{\Gamma_h}(\{{ \underline{\btau}_h}\}-[\underline{\btau}_h]\otimes\textbf{C}_{12}):\underline{[{\bu_h}]}\,ds + \int_{\Gamma_f}\underline{\btau}_h:(\bu_h\otimes\bn)\, ds. \label{semi1}
\end{align}
Using the Cauchy-Schwarz inequality and Lemma \ref{TL2}, 
we obtain
\begin{align*}
\sum_{E\in \mathcal{E}_{fh}} \int_{E} \textbf{D}(\bu_h):\underline{\btau}_h {d\bx} \leq \|{ \underline{\bs}_h}\|_{L^2(\Omega_f)}\|{ \underline{\btau}_h}\|_{L^2(\Omega_f)} + C\, c_{NS}(\bu_h,\bu_h)^{1/2}\|{ \underline{\btau}_h}\|_{L^2(\Omega_f)}.
\end{align*}
Now, applying the definition of the $\|\cdot\|_{\bX_h}$--norm and { $|\cdot|_{\mathcal{A}}$--seminorm,} we {can achieve} the first estimate of Lemma \ref{equivalence}.
Again, putting ${ \underline{\btau}_h}={ \underline{\bs}_h}$ in (\ref{DFSLNSD1}), using {the} Cauchy-Schwarz inequality, Lemma \ref{TL2}, and applying similar steps as above, we obtain the second estimate. Similar steps will also be followed for (\ref{e2}).
\end{proof}
\end{lemma}
\noindent Below, we provide some estimates of the form { $O_{NS}(\cdot,\cdot,\cdot)$}. 
\begin{lemma} \label{con}
The following inequalities for the trilinear term hold true:
\small{\begin{align}
&|O_{NS}(\bw,\bu,\bv_h)|\leq C \|\bw\|_{\bX_h}\|\bu\|_{\bX_h}\|\bv_h\|_{\bX_h},\,\bw,\bu\in\bX_1,\bv_h\in\bX_h,\label{trilinear1}\\
&|O_{NS}(\bz_h,\bu_h,\bv_h)|\leq C \|\bz_h\|_{\bX_h}\|\bu_h\|_{\bX_h}\|\bv_h\|_{\bX_h},\,\bz_h,\bu_h,\bv_h\in\bX_h, \label{continuity}\\
&|O_{NS}(\bw_h,\bu,\bv_h)|\leq C \|\bw_h\|_{L^2(\Omega_f)}\|\bu\|_{2,\Omega_f}\|\bv_h\|_{\bX_h},\,\bu\in W^{1,4}(\Omega_f),\bw_h,\bv_h\in\bX_h\label{trilinear2}.
\end{align}}
\end{lemma}
\begin{proof}
The trilinear term in (\ref{trilinear}) can be written as
\small{\begin{align}
O_{NS}(\bw,&\bu,\bv_h) = -\sum_{E\in \mathcal{E}_{fh}} \int_{E} (\nabla\cdot\bw)\bu\cdot\bv_h{ d\bx} -\sum_{E\in \mathcal{E}_{fh}} \int_{E} \bu\cdot(\bw\cdot\nabla)\bv_h{ d\bx} \nonumber\\
&+\bigg( \sum_{e\subset\Gamma_f^+} \int_{e}(\hat{\bu}^\bw\otimes\bw): (\bv_h\otimes\bn)\, ds + \sum_{e\subset\Gamma_{fh}} \int_{e}(\hat{\bu}^\bw\otimes\bw): \underline{[\bv_h]} ds\bigg) = T_1+T_2+T_3,\label{outflow}
\end{align}}
\noindent where \[\Gamma_f^+=\{\bx\in\Gamma_f:\bw(\bx)\cdot\bn(\bx)\geq {0}\},\] represents the outflow part of $\Gamma_f$.\par
\noindent Using the Cauchy-Schwarz inequality and (\ref{brokenpoincare}), the first two terms can be bounded as follows:
\small{\begin{align*}
|T_1|+|T_2|&\leq \sum_{E\in \mathcal{E}_{fh}} \|\nabla\bw\|_{L^2(E)}\|\bu\|_{L^4(E)}\|\bv_h\|_{L^4(E)} + \sum_{E\in \mathcal{E}_{fh}} \|\bu\|_{L^4(E)}\|\bw\|_{L^4(E)}\|\nabla\bv_h\|_{L^2(E)}\leq C\|\bw\|_{\bX_h}\|\bu\|_{\bX_h}\|\bv_h\|_{\bX_h}.
\end{align*}}
\noindent Inserting the definition of the upwind flux and applying Lemma \ref{TL1}, we obtain
\small{ \begin{align*}
 |T_3|&\leq \sum_{e\subset\Gamma_f^+\cup\Gamma_{fh}}
 \int_{e}\|\bw\|_{L^4(e)}\|\bu\|_{L^4(e)}\|[\bv_h]\|_{L^2(e)}\, ds\\
 &\leq C|e|^{1/2}
 \sum_{E\in \mathcal{E}_{fh}} h_E^{-1/2}(\|\bw\|_{L^2(E)}+h_E^{1/2}\|\bw\|_{1,E})(\|\bu\|_{L^2(E)}+h_E^{1/2}\|\bu\|_{1,E})\|\bv_h\|_{\bX_h}\leq C\|\bw\|_{\bX_h}\|\bu\|_{\bX_h}\|\bv_h\|_{\bX_h}.
 \end{align*}}
\noindent Combining the above bounds, we arrive at the estimate (\ref{trilinear1}).
Following similar steps, we can obtain (\ref{continuity}). For the third estimate, we can bound the first two terms in (\ref{outflow})  using a similar technique as in \cite[Proposition 4.1, $(4.1)$]{GRW005S}. For the third term, using Lemma \ref{TL2} we obtain
\begin{align*}
    T_3
    \leq C\|\bu\|_{L^\infty(\Omega_f)}\sum_{E\in \mathcal{E}_{fh}}\|\bw_h\|_{L^2(E)}\|\bv_h\|_{\bX_h}\leq C\|\bu\|_{2,\Omega_f}\|\bw_h\|_{L^2(\Omega_f)}\|\bv_h\|_{\bX_h}.
\end{align*}
Incorporating all the resulting bounds in (\ref{outflow}), we achieve (\ref{trilinear2}).
\end{proof}
\noindent Following the techniques utilized for { \cite[pp. 577, (2.9)]{CKS04}}, we have obtained the following lemma, but it does not directly follow the standard ellipticity property due to the involvement of interface terms. A similar result can be found in \cite[Lemma 4.1]{GR09}.
\begin{lemma}\label{non}
The term $O_{NS}(\cdot,\,\cdot\,,\cdot)$ satisfies the following property:
\begin{align}
O_{NS}(\Z,\bu_h,\bu_h)\geq  \frac{1}{2}  (\Z\cdot\bn ,\bu_h\cdot \bu_h)_{ \Gamma}\quad \forall { \bu_h  \in}{\bX}_h,\,\Z\in \textbf{J}(\mathcal{E}_{fh}). \label{nonlinearity1}
\end{align}
\end{lemma}
\noindent We now state the existence and uniqueness of the proposed scheme. Since this is a linear system in finite dimensions, it suffices to show the uniqueness of the solution. But it is worth mentioning here that although, at each time step, (\ref{DFSLNSD1})-(\ref{DFSLNSD4}) is a square linear system of finite dimension, proving the nonsingularity of its matrix without restricting the time step is not straightforward (see \cite[Section 5.2]{CGPR17}). 
\begin{lemma}\label{unique}
If $\bU_h^{i-1}$ satisfies 
\begin{align}
\|\bU_h^{i-1}\|_{\bX_h}\leq \frac{\mu }{C_a C_0K_\mathcal{A}}, \label{1}
\end{align}
then the solution  $(\underline{\bSig}^i_{h},\bU^i_{h},\bQ^i_{h},\Phi^i_{h})$  is unique. 
\end{lemma}
\begin{proof}
Assume that at the $i$-th time step (\ref{DFSLNSD1})-(\ref{DFSLNSD4}) have two solutions $(\underline{\bSig}^i_{1h},\bU^i_{1h},{\bQ}^i_{1h},\Phi^i_{1h})$ and $(\underline{\bSig}^i_{2h},\bU^i_{2h}$, $\bQ^i_{2h},\Phi^i_{2h})$, then their difference $(\underline{\bSig}^i_{Dh},\bU^i_{Dh},\bQ^i_{Dh},\Phi^i_{Dh})$ solves
{\scriptsize
\begin{align*}
  a_{NS}(\underline{\bSig}^i_{Dh},\underline{\btau_{h}}) +  b_{NS}(\bU^i_{Dh}, \underline{\btau}_h) &= 0,\\ 
  ( \partial_t\bU^i_{Dh},\bv_h)- b_{NS}(\bv_h,\underline{\bSig}^i_{Dh}) +  c_{NS}(\bU^i_{Dh},\bv_h) + O_{NS}(\p\bU^{i-1}_{h},&\bU^i_{Dh},\bv_h) \nonumber\\ 
  + \sum_{j=1}^{d-1}\frac{1}{G^j}(\bU^i_{Dh}\cdot\btau_{\Gamma},\bv_h\cdot\btau_{\Gamma})_{{\Gamma}}
+  g(\Phi^i_{Dh},\bv_h\cdot{\bn_{\Gamma}})_{{\Gamma}} & = 0,\\
a_{D}(\bQ^i_{Dh},\bw_h)-b_{D}(\Phi^i_{Dh},{\bw}_h)&=0,\\
gS_0(\partial_t\Phi^i_{Dh},\psi_h)+b_{D}(\psi_h,\bQ^i_{Dh})+ c_{D}(\Phi^i_{Dh},\psi_h)-g(\bU^i_{Dh}\cdot\bn_{\Gamma},\psi_h)_{{\Gamma}}&= 0,
\end{align*}}
\noindent $\forall (\underline{\btau}_h,\bv_h,q_h,\bw_h,\psi_h)\in \underline{\bY}_h\times \bZ_h\times M_{fh}\times {\bW}_h\times  M_{ph}$.

\noindent Now, putting $\underline{\btau}_h=2\mu\underline{\bSig}^i_{Dh}$, $\bv_h=\bU^i_{Dh}$, $\bw_h=g\bQ^i_{Dh}$, and $\psi_h=\Phi^i_{Dh}$, adding all equations and using (\ref{nonliearity2}) and (\ref{nonlinearity1}), we obtain 
\begin{equation*}\scalebox{0.9}{$
\frac{1}{\Delta t} \|\bU^i_{Dh}\|^2_{L^2(\Omega_f)} + \frac{gS_0}{\Delta t} \|\Phi^i_{Dh}\|^2_{L^2(\Omega_p)} + \mu |(\underline{\bSig}^i_{Dh},\bU^i_{Dh})|^2_{\mathcal{A}} + g |(\tilde{\bQ}^i_{Dh},\phi_h)|^2_{\mathcal{B}} + \sum_{j=1}^{d-1}\frac{1}{G^j}\|\bU^i_{Dh}\cdot\btau_{\Gamma}\|^2_{L^2({\Gamma})}
\leq C_0 \|\p\bU^{i-1}_{h}\|_{\bX_h}\|\bU^i_{Dh}\|^2_{\bX_h}.$}
\end{equation*}
\noindent Then applying (\ref{Postprocess1}) and (\ref{1}), we find $\underline{\bSig}^i_{Dh}=\textbf{0},~~\bU^i_{Dh}=\textbf{0},~~\bQ^i_{Dh}=\textbf{0},~~\Phi^i_{Dh}=0$.
\end{proof}
\noindent The existence and uniqueness of the solution of (\ref{DFSLNSD1})-(\ref{DFSLNSD4}) at each time step can be shown as in Section $5.2$ in \cite{CGPR17}. 
Below, in Lemma \ref{apriori}, we present the stability of the scheme. 
\begin{lemma} \label{apriori}
Let $0<\alpha<min\{\frac{\mu }{2C_2K_A},\frac{ g}{2C_1K_B}\}$. 
If for all $0\leq i\leq N-1$, 
\begin{equation}\scalebox{0.98}{$
\|\bU^{i}_h\|_{\bX_h}\leq\frac{\mu}{2C_aC_0K_A},$} \label{aprioricondition1}
\end{equation}
then the discrete solution $(\{\underline{\bSig}^i_{h}\},\{\bU^i_h\},\{\bQ^i_{h}\}, \{\Phi^i_h\})_{1\geq n\leq N}$ of (\ref{DFSLNSD1})-(\ref{DFSLNSD4}) satisfies 
{\small\begin{align*}
\|\bU^n_h\|^2_{L^2(\Omega_f)} +  gS_0\|\Phi^n_h\|^2_{L^2(\Omega_p)} + e^{-2\alpha t_n}\Delta t \sum_{i=1}^{n} e^{2\alpha t_i}|( \underline{\bSig}^i_h,\bU^i_h)|^2_{\mathcal{A}} + e^{-2\alpha t_n}\Delta t \sum_{i=1}^{n} e^{2\alpha t_i}|({\bQ}_h^i,\Phi_h^i)|^2_{\mathcal{B}}\leq C.
\end{align*}}
Here the constant $C=e^{2\alpha \Delta t}\big( \|\bU_h^0\|^2 + \|\Phi_h^0\|^2 \big) + \frac{C_2^2K_A}{2\mu }\Delta t \sum_{i=1}^{n} e^{2\alpha t_i} \|\f_f^i\|^2_{L^2(\Omega_f)} + \frac{C_1^2K_B}{2  g}\Delta t \sum_{i=1}^{n} e^{2\alpha t_i} \|f_p^i\|^2_{L^2(\Omega_p)}$. 
\end{lemma}
\section{Approximation Operators}\label{s05}
In order to derive optimal error estimates, the error analysis mainly revolves around weakly divergence-free spaces. Consequently, it becomes essential to define an approximation operator in $\bV\cap \bH^{k_1+1}(\Omega_f)$.
\begin{lemma}\label{iprojection}
There is an approximation $i_h:\bV\cap \bH^{k_1+1}(\Omega_f) \rightarrow \bZ_h$, such that
\begin{align}
\|\bz-i_h\bz\|_{L^2(\Omega_f)} + h\|\bz-i_h\bz\|_{\bX_h} \leq Ch^{k_1+1} |\bz|_{k_1+1,\Omega_f}, \forall\bz\in\bV\cap\bH^{k_1+1}(\Omega_f).
\end{align}
 Here, $i_h$ is the restriction of $\bPi_h$ on $\bV\cap \bH^{k_1+1}(\Omega_f)$.
\end{lemma}
\noindent Next, we define the following $\bL^2$ projections, which will be used in the subsequent analysis:
$\textbf{P}_h^{\bs}: \underline{\bL}^2(\Omega_f)\rightarrow \underline{\bY}_h$, $\bP_h^\bu:{\bL}^2(\Omega_f)\rightarrow\bZ_h$, and $\textbf{P}_h^{\bomega}: \underline{\bL}^2(\Omega_p)\rightarrow \bW_h$,
{\scriptsize\begin{align}
&\,(\underline{\bs}-\textbf{P}_h^{\bs}\underline{\bs},\underline{\btau}_h)  =0, \quad \underline{\bs} \in \underline{\bL}^2(\Omega_f),\,\,  \underline{\btau}_h \in  \underline{\bY}_h,\\
&\,(\bu-\textbf{P}_h^{\bu}\bu,\bv_h)  \,\,=0,\quad \bu\in \bL^2(\Omega_f),\,\, \bv_h\in \bZ_h,\\
&(\bomega-\textbf{P}_h^{\bomega}\bomega,\bw_h) =0,\quad \bw\in {\bL}^2(\Omega_p),\,\, \bw_h\in \bW_h.
\end{align}}
\begin{lemma}\label{allprojections}
The projections $\textbf{P}_h^{\bs}$, $\textbf{P}_h^{\bu}$, $\textbf{P}_h^{\bomega}$ satisfy 
{\scriptsize\begin{align}
&\|\bu-\textbf{P}_h^{\bu}\bu\|_{L^2(\Omega_f)}+h\|\bu -\textbf{P}_h^{\bu}\bu\|_{\bX_h}\leq Ch^{k_1+1} |\bu|_{k_1+1,\Omega_f}, \,\, \forall \bu\in \bV\cap\bH^{k_1+1},\label{projection2}\\
&\|\underline{\bs}-\textbf{P}_h^{\bs}\underline{\bs}\|_{L^2(\Omega_f)} 
\leq Ch^{k_1}|\underline{\bs}|_{k_1,\Omega_f}\,\, \forall \underline{\bs}\in \underline{\bY},\label{N1}\\
&\|\underline{\bs}-\textbf{P}_h^{\bs}\underline{\bs}\|_{L^2(e)} \leq Ch_E^{k_1/2}\|\underline{\bs}\|_{L^2(E)},\,\, \forall \underline{\bs}\in \underline{\bY},\label{N2}\\
&\|\bw-\textbf{P}_h^{\bomega}{\bw}\|_{L^2(\Omega_p)} 
\leq Ch^{k_2}|{\bw}|_{k_2,\,\Omega_p},\,\, \forall \bw\in \bW,\label{N3}\\
&\|\bw-\textbf{P}_h^{\bomega}\bw\|_{L^2(e)} \leq Ch^{k_2/2}|\bw|_{k_2,\Omega_p},\,\, \forall \bw\in \bW,\label{N4}\\
&|(\underline{\bs}-\bPa\underline{\bs},\bu-\bPb\bu)|_{\mathcal{A}} \leq Ch^{k_1}(|\underline{\bs}|_{k_1,\Omega_f}+|\bu|_{k_1+1,\Omega_f})\,\, \forall \underline{\bs}\in \underline{\bY},  \,\forall \bu\in \bV\cap\bH^{k_1+1},\label{N5}\\
&|(\bw-\textbf{P}_h^{\bomega}\bw,\phi-\Pd\phi)|^2_{\mathcal{B}}\leq Ch^{k_2}(|\bw|_{k_2,\Omega_p}+|\phi|_{k_2+1,\Omega_p})\,\, \forall \bw\in \bW,\, \forall\phi\in M_p.\label{N6}
\end{align}}
\end{lemma}
\noindent The proof of (\ref{projection2}) follows from the orthogonal property of $\textbf{P}_h^{\bu}$ and (\ref{iprojection}). For the proof, one can refer to \cite[Lemma 3.1]{CKSS02} and \cite[Lemma 3.1]{CCPS00}.
Also, an application of the trace inequality (\ref{T4}), the Lagrange interpolant, and (\ref{projection2}) yields the stability of the projection $\textbf{P}_h^{\bu}$ as follows:
\begin{align}
\|\textbf{P}_h^{\bu}\bv\|_{\bX_h}\leq C|\bv|_1,\quad\forall\bv\in\bV.\label{stability} 
\end{align}
\noindent We observe that higher-order polynomial approximation requires sufficient regularity for the solution $(\underline{\bs},\bu,p,$ $\bomega,\phi)$, and therefore we assume the following:

\noindent $(\textbf{A}\textbf{1}).$ The exact solution satisfies:\\
\small{\begin{tabular}{c c c }
	&$\bu \in L^{\infty}(0\,,T;\bH^{k_1+1}(\Omega_f))\cap L^{2}(0\,,T;\bH^{k_1+1}(\Omega_f) )$, $p \in L^{\infty}(0\,,T;\bH^{k_1}(\Omega_f))\cap L^{2}(0\,,T;\bH^{k_1}(\Omega_f) )$,     \\
	&$\underline{\bs} \in L^{\infty}(0,T;\underline{\bH}^{k_1}(\Omega_f))\cap L^{2}(0,T;\underline{\bH}^{k_1}(\Omega_f) )$, $\bomega \in L^{\infty}(0,T;\bH^{k_2}(\Omega_f))\cap L^{2}(0,T;\bH^{k_2}(\Omega_f) ),$     \\
	&$\phi  \in L^{\infty}(0\,,T;H^{k_2+1}(\Omega_p))\cap L^{2}(0\,,T;H^{k_2+1}(\Omega_p) )$,\,\, $p_t\in L^2(0,T;H^{k_1}(\Omega_f))$ ,        \\
   & $\bu_t \in L^{2}(0\,,T;\bH^{k_1+1}(\Omega_f) )\cap L^{\infty}(0\,,T;\bH^{k_1}(\Omega_f))$,\,\, $\bu_{tt} \in L^{2}(0\,,T;\bL^2(\Omega_f) )$,  \\
 &  $\phi_t \in L^{2}(0\,,T;H^{k_2+1}(\Omega_p) )\cap L^{\infty}(0\,,T;H^{k_2}(\Omega_p))$,\,\, $\phi_{tt}  \in L^{2}(0\,,T;L^2(\Omega_p) )$. \\
		\end{tabular}}
		
\noindent We now present below the modified Oseen operators, which are crucial for our $L^2$-error analysis. This and the rest of our analysis will be carried out for the two dimensional case; for the three dimensional case, please see Remark \ref{3D}.
\subsection{Modified Oseen Operators}
\noindent As mentioned in the 'Introduction', the NSD problem faces difficulty due to the presence of the interface terms, and the standard operators have not been adequate to tackle the interface terms in the case of $L^{\infty}(\textbf{L}^2)$ estimates. Therefore, we have defined the modified Oseen operators $\bOa$ for stress tensor, $\bOb$ for velocity, $\bOc$ for velocity in porous media, and $\Od$ for piezometric head, respectively, satisfying 
\small{\begin{align}
a_{NS}(\underline{\bs} - \bOa\underline{\bs},&\underline{\btau}_h) + b_{NS}(\bu - \bOb \bu, \underline{\btau}_h ) = 0,\,\, \forall\underline{\btau}_h \in \underline{\bY}_h, \label{Oseen1}\\
 - b_{NS}(\bv_h,\underline{\bs} - \bOa\underline{\bs}) +  c_{NS}(\bu - \bOb &\bu,\bv_h) + O_{NS}(\bu,\bu - \bOb \bu,\bv_h)   
 + O_{NS}(\bu - \bOb \bu,\bu,\bv_h) \nonumber\\
 + \frac{1}{G}((\bu - \bOb \bu)\cdot\btau_{\Gamma},\bv_h\cdot\btau_{\Gamma}&)_{{\Gamma}}
+  g((\phi - \Od \phi) ,\bv_h\cdot{\bn_{\Gamma}})_{{\Gamma}}  = -d_{NS}(\bv_h,p),\,\, \forall\bv_h \in \bZ_h,\label{Oseen2}\\
a_{D}(\bomega - \bOc \bomega,&\bw_h)-b_{D}(\phi - \Od \phi,\bw_h)=0,\,\, \forall \bw_h \in \bW_h,\label{Oseen3}\\
b_{D}(\psi_h,\bomega  -\bOc \bomega )+ c_{D}(\phi - \Od \phi ,&\psi_h)-g((\bu - \bOb \bu)\cdot\bn_{\Gamma},\psi_h)_{{\Gamma}}= 0,\,\, \forall \psi_h \in M_{ph},\label{Oseen4}
\end{align}}
given $(\underline{\bs},\bu,p,\bomega,\phi)$ is the weak solution of the NSD system. 

\noindent Being a linearized steady NSD problem, the existence and uniqueness of (\ref{Oseen1})-(\ref{Oseen4}) can be proved in similar ways as in Lemma \ref{unique}. The approximation properties of the operators $\underline{\bs}-\bOa\underline{\bs}$, $\bu-\bOb\bu$, $\bomega-\bOc\bomega$, and $\phi-\Od\phi$ follow next.
\begin{lemma} \label{Opproperty}
If $\bu$ satisfies
\small{\begin{align}
    \|\textbf{D}(\bu)\|_{L^2(\Omega_f)} \leq \frac{\mu}{4d_1d_2},\label{EO1}
\end{align}}
the terms $\underline{\bs}-\bOa\underline{\bs}$, $\bu-\bOb\bu$, $\bomega-\bOc\bomega$ and $\phi-\Od\phi$ satisfy 
\small{\begin{align}
& \|\bu - \bOb \bu\|^2_{L^2(\Omega_f)}+\|\phi-\Od\phi\|^2_{L^2(\Omega_p)} + h^2\Big(|(\underline{\bs} - \bOa\underline{\bs},\bu - \bOb \bu)|^2_{\mathcal{A}} 
 +|(\bomega - \bOc \bomega,\phi - \Od \phi)|^2_{\mathcal{B}}\Big)\nonumber\\
&\quad\quad\quad\quad\quad\quad\quad\quad
\leq C (h^{2k_1+2}(|\bu|^2_{k_1+1,\Omega_f}
+|p|^2_{k_1,\Omega_f})+h^{2k_2+2}|\phi|^2_{k_2+1,\Omega_p}),\label{Oseenprojection1}\\
&\|(\bu - \bOb \bu)_t\|^2_{L^2(\Omega_f)}+\|(\phi-\Od\phi)_t\|^2_{L^2(\Omega_p)}+h^2|((\underline{\bs} - \bOa\underline{\bs})_t,(\bu - \bOb \bu)_t)|^2_{\mathcal{A}} \nonumber\\
&\,\,+h^2|((\bomega - \bOc \bomega)_t,(\phi - \Od \phi)_t)|^2_{\mathcal{B}}
\leq C (h^{2k_1+2}(|\bu_t|^2_{k_1+1,\Omega_f}
+|p_t|^2_{k_1,\Omega_f})+h^{2k_2+2}|\phi_t|^2_{k_2+1,\Omega_p}),\label{Oseenprojection2}
\end{align}}
where $C = C(\mu,g,d_1,d_2, K_\mathcal{A}, K_\mathcal{B})$ is a positive constant independent of $h$.
\end{lemma}
\begin{proof}
Since \small{\begin{align*}
&\|\underline{\bs} - \bOa\underline{\bs}\|_{L^2(\Omega_f)} \leq \|\underline{\bs} - \bPa\underline{\bs}\|_{L^2(\Omega_f)} + \|\bPa\underline{\bs} - \bOa\underline{\bs}\|_{L^2(\Omega_f)},\\
&\|\bu - \bOb \bu\|_{\bX_h} \leq \|\bu-\bPb\bu\|_{\bX_h}+ \|\bPb \bu - \bOb \bu\|_{\bX_h},\\
&\|\bomega - \bOc \bomega\|_{L^2(\Omega_p)} \leq \|\bomega - \bPc \bomega\|_{L^2(\Omega_p)} + \|\bPc \bomega - \bOc \bomega\|_{L^2(\Omega_p)},\\
&\|\phi - \Od \phi\|_{M_{ph}} \leq \|\phi-\Pd \phi\|_{M_{ph}}+\|\Pd \phi-\Od\phi\|_{M_{ph}},
\end{align*}}
\noindent it is sufficient to estimate $\bPa\underline{\bs} - \bOa\underline{\bs}$, $\bPb \bu - \bOb \bu $, $\bPc \bomega - \bOc \bomega$ and $\Pd \phi-\Od\phi $.
We choose $\underline{\btau}_h=\bPa\underline{\bs} - \bOa\underline{\bs}=:\bDs$,
$\bv_h=\bPb \bu - \bOb \bu=:\bDu$, $\bw_h=\bPc \bomega - \bOc \bomega=: \bDq$, $\psi_h=\Pd \phi-\Od\phi=: \Dp$ in (\ref{Oseen1})-(\ref{Oseen4}) and notice that, since $d_{NS}(\bDu,r_hp)=0$,
\small{\begin{align}
&a_{NS}(\bDs,\bDs)+c_{NS}(\bDu,\bDu)+a_{D}(\bDq,\bDq)+c_{D}(\Dp,\Dp)+\frac{1}{G}(\bDu\cdot\btau_\Gamma,\bDu\cdot\btau_\Gamma)_\Gamma \leq -d_{NS}(\bDu,p-r_hp)\nonumber\\
&\quad-a_{NS}(\underline{\bs} - \bPa\underline{\bs},\bDs)-b_{NS}(\bu - \bPb \bu, \bDs )+b_{NS}(\bDu, \underline{\bs}- \bPa \underline{\bs} )-c_{NS}(\bu - \bPb \bu,\bDu)-a_{D}(\bomega - \bPc\bomega,\bDq)\nonumber\\
&\quad\,\, +b_{D}(\phi - \Pd \phi, \bDq )-b_{D}(\bDq, \bomega - \bPc\bomega ) -c_{D}(\phi - \Pd \phi,\Dp)-\frac{1}{G}((\bu - \bPb \bu)\cdot\btau_\Gamma,\bDu\cdot\btau_\Gamma)_\Gamma\nonumber\\
&\quad\,\,-g(\phi-\Pd\phi,\bDu\cdot\bn_{\Gamma})_{\Gamma} +g((\bu - \bPb \bu)\cdot\bn_{\Gamma},\Dp)_{\Gamma}- O_{NS}(\bDu,\bu,\bDu) -O_{NS}(\bu,\bDu ,\bDu)-O_{NS}(\bu,\bu-\bPb\bu ,\bDu)\nonumber\\
&\quad\,\,- O_{NS}(\bu-\bPb\bu,\bu,\bDu)\,\,= \sum_{i=1}^{16}S_i.\label{Add1}
\end{align}}
\noindent Now, using the definition of the norms and seminorms, the left hand side of (\ref{Add1}) can be bounded below by
\small{\begin{align*}
\mu|(\bDs,\bDu)|^2_{\mathcal{A}} +g |(\bDq,\Dp)|^2_{\mathcal{B}}+ \frac{1}{G}\|\bDu\cdot\btau_{\Gamma}\|^2_{L^2(\Gamma)}.
\end{align*}}
\noindent To bound the right hand side of (\ref{Add1}), we apply the standard Cauchy-Schwarz and Young's inequalities, Lemma \ref{boundedness}, the projection properties (\ref{projection1defn}), (\ref{projection1}), (\ref{projection3}), (\ref{projection2}), (\ref{N5}), (\ref{N6}), and the trace inequalities Lemma \ref{TL1} and \ref{TL2}, and obtain 
\small{\begin{align*}
&|S_1|
= \bigg| \sum_{e\in\Gamma_{fh}}\int_{e}(\{p-r_hp\}
-\textbf{D}_{12}\cdot[p-r_hp]
)[\bDu]\,{ d\bx} 
+ \int_{\Gamma_f} (p-r_hp)\bDu\cdot\bn \, ds\bigg|\leq \frac{\mu}{12}|(\bDs,\bDu)|^2_{\mathcal{A}}+Ch^{2k_1}|p|^2_{k_1,\Omega_f},\\
&|(S_2+S_5)+(S_6+S_{9})+S_3+S_{7}|\leq \frac{\mu}{12}|(\bDs,\bDu)|^2_{\mathcal{A}} +  \frac{g}{8}|(\bDq,\Dp)|^2_{\mathcal{B}}+C(h^{2k_1}(|\underline{\bs}|^2_{k_1,\Omega_f}+|\bu|^2_{k_1+1,\Omega_f})\\
&\quad\quad\quad\quad\quad\quad\quad\quad\quad\quad\quad\quad\quad\quad\quad
+h^{2k_2}(|\bomega|^2_{k_2,\Omega_p}+|\phi|^2_{k_2+1,\Omega_p}))+Ch^{2k_1}|\bu|^2_{k_1+1,\Omega_f}+Ch^{2k_2}|\phi|^2_{k_2+1,\Omega_p}.
\end{align*}}
\noindent A use of the Cauchy-Schwarz inequality, Lemma \ref{equivalence}, (\ref{N1}), (\ref{N3}), and Young's inequality yields
\small{\begin{align*}
&|S_4|
= \bigg|-\sum_{E\in \mathcal{E}_{fh}} \int_{E} \textbf{D}(\bDu):(\underline{\bs}- \bPa \underline{\bs}) d\bx + \sum_{e\in\Gamma_{fh}}\int_{e}(\{\underline{\bs}- \bPa \underline{\bs}\}-[\underline{\bs}- \bPa \underline{\bs}]\otimes\textbf{C}_{12}):\underline{[{\bDu}]}\,ds\\
&\quad + \int_{\Gamma_f}(\underline{\bs}- \bPa \underline{\bs}):(\bDu\otimes\bn)\, ds\bigg|
\leq \|\bDu\|_{ \bX_h} \|\underline{\bs}- \bPa \underline{\bs}\|_{L^2(\Omega_f)}+Ch^{k_1}|\bs|_{k_1,\Omega_f}c_{NS}(\bDu,\bDu)^{1/2}\\
&\quad\quad\leq Ch^{k_1}|\underline{\bs}|_{k_1,\Omega_f}|(\bDs,\bDu)|_{\mathcal{A}} +Ch^{k_1}|\underline{\bs}|_{k_1,\Omega_f}c_{NS}(\bDu,\bDu)^{1/2}
\leq \frac{\mu}{12}|(\bDs,\bDu)|^2_{\mathcal{A}}+Ch^{2k_1}|\underline{\bs}|^2_{k_1,\Omega_f},\\ 
&\,\,|S_{8}|\leq \frac{g}{8}|(\bDq,\Dp)|^2_{\mathcal{B}}+Ch^{2k_2}|\bomega|^2_{k_2,\Omega_p}.
\end{align*}}
\noindent Further, we apply Lemma \ref{TL1}, (\ref{T9}), and (\ref{projection2}) along with the standard inequalities to bound the interface terms as
\small{\begin{align*}
&|S_{10}|=\big|\frac{1}{G}((\bu - \bPb \bu)\cdot\btau_{\Gamma},\bDu\cdot\btau_{\Gamma})_{{\Gamma}}\big|\leq \frac{1}{2G}\|\bDu\cdot\btau_{\Gamma}\|^2_{L^2(\Gamma)}+Ch^{2k_1+1}|\bu|^2_{k_1,\Omega_f},\\
&|S_{11}+S_{12}|\leq \frac{\mu}{12}|(\bDs,\bDu)|^2_{\mathcal{A}}+ \frac{g}{8}|(\bDq,\Dp)|^2_{\mathcal{B}}+Ch^{2k_1+1}|\bu|^2_{k_1+1,\Omega_f}+Ch^{2k_2+1}|\phi|^2_{k_2+1,\Omega_p}.
\end{align*}}
\noindent Next, we bound the trilinear terms using (\ref{ctsineq}), (\ref{ST8}), (\ref{T8}), and (\ref{EO1}) and obtain

\small{\begin{align*}
|S_{13}|&=\Big|-\sum_{E\in \mathcal{E}_{fh}} \int_{E} (\nabla\cdot\bDu)\bu\cdot\bDu{ d\bx} -\sum_{E\in \mathcal{E}_{fh}} \int_{E} \bu\cdot(\bDu\cdot\nabla)\bDu{ d\bx}+ \sum_{e\subset\Gamma_f^+\cup\Gamma_{fh}} \int_{e}\bigg(\bDu\cdot\bn\{\bu\}\\
&-\frac{1}{2}|\bDu\cdot\bn|(\bu^{ext}-\bu)\bigg):\underline{[\bDu]}\, ds\Big|\leq d_1d_2\|\textbf{D}(\bu)\|_{L^2(\Omega_f)}\|\bDu\|^2_{\bX_h}\leq \frac{\mu}{4}|(\bDs,\bDu)|^2_{\mathcal{A}},\\
|S_{14}|&\leq \frac{\mu}{4}|(\bDs,\bDu)|^2_{\mathcal{A}}.
\end{align*}}
\noindent Also, the rest of the trilinear terms can be bounded using (\ref{trilinear1}) along with Young's inequality as
\begin{align*}
|S_{15}|+|S_{16}|
\leq \frac{\mu}{12}|(\bDs,\bDu)|^2_{\mathcal{A}} + Ch^{2k_1}|\bu|^2_{k_1+1,\Omega_f}.
\end{align*}
\noindent Incorporating all the values in (\ref{Add1}), we obtain
\begin{align*}
\mu|(\bDs,\bDu)|^2_{\mathcal{A}} +g |(\bDq,\Dp)|^2_{\mathcal{B}}+ \frac{1}{G}\|\bDu\cdot\btau_{\Gamma}\|^2_{L^2(\Gamma)} \leq C (h^{2k_1}(|\bu|^2_{k_1+1,\Omega_f} +|p|^2_{k_1,\Omega_f})+h^{2k_2}|\phi|^2_{k_2+1,\Omega_p}).
\end{align*}
\noindent Since $\underline{\bs} - \bOa\underline{\bs}_{L^2(\Omega_f)} =\underline{\bs} - \bPa\underline{\bs} + \bDs$, $\bu - \bOb \bu= \bu-\bPb\bu+ \bDu$, $\bomega - \bOc \bomega=\bomega - \bPc \bomega + \bDq$, and $\phi - \Od \phi=\phi-\Pd \phi+\Dp$, using approximation properties (\ref{projection3}), (\ref{projection2}), (\ref{N1}), and (\ref{N3}), we complete the following energy norm estimate:
\small{\begin{align}
|(\underline{\bs} - \bOa\underline{\bs},\bu - \bOb \bu)|^2_{\mathcal{A}} +|(\bomega - \bOc \bomega,\phi - \Od \phi)|^2_{\mathcal{B}}
\leq C& (h^{2k_1}(|\bu|^2_{k_1+1,\Omega_f}
+|p|^2_{k_1,\Omega_f}) +h^{2k_2}|\phi|^2_{k_2+1,\Omega_p}).\label{Oseenmid1}
\end{align}}%
\noindent For the $L^2$-norm estimate, we apply the Aubin-Nitsche duality argument.  
Given $\textbf{g}_1\in \bL^2(\Omega_f)$, $g_2\in L^2(\Omega_p)$, the solution $(\bz,s,\varrho)$ satisfies the following dual problem:
\small{\begin{align*}
   \nabla\cdot(2\mu \textbf{D}(\bz)-s{\bf I})- \nabla\cdot(\bz\otimes\bu)+(\nabla\bu)^T  \bz &= \textbf{g}_1, ~~\mbox{in}~ \Omega_f, \\
\nabla \cdot \bz &= 0,~~\mbox{in}~ \Omega_f,\\
-\nabla\cdot (\bK\nabla\varrho) & = g_2~~ \mbox{in}~~ \Omega_p.
	\end{align*}}%
    \noindent The boundary and interface conditions are the same as in the original problem. 
Observe that the two trilinear terms in the dual problem are bounded for a sufficiently large $\mu$, and the regularity assumption in \textbf{(A1)}. Now following the argument in \cite[pp. 3332]{CMX09}, we may assume that the solution of the dual problem satisfies the following regularity result:
\begin{align*}
\bigg(\|\bz\|^2_{2,\Omega_f}+\|\varrho\|^2_{2,\Omega_p}\bigg)^{\frac{1}{2}}+\|s\|_{1,\Omega_f}\leq  \|\textbf{g}_1\|_{L^2(\Omega_f)}+\|g_2\|_{L^2(\Omega_p)}.
\end{align*}
Now, we set $\underline{\bl}=-\textbf{D}(\bz)$ and $\bbeta=-\bK\nabla\varrho$ in the above equations. Then, multiplying the resulting equations by the test functions $(\underline{\btau},\bv, q,\bw,\psi)\in \underline{\bY}\times\bX\times M_f\times\textbf{W}\times M_p$ and integrating over the domains $\Omega_f$ and $\Omega_p$, we obtain the following mixed variational formulation: find $(\underline{\bl},\bz,s,{\bbeta},\varrho)\in \underline{\bY}\times\bV\times M_f \times\textbf{W}\times M_{p}$ such that  $\forall(\underline{\btau},\bv, q,\bw,\psi)\in \underline{\bY}\times\bX\times M_f\times\textbf{W}\times M_p$,
{\small\begin{align*}
(\underline{\bl},\underline{\btau}) + (\textbf{D}(\bz),\underline{\btau})&=0, \nonumber\\
-2\mu (\underline{\bl},\textbf{D}(\bv))- (\bu \cdot \nabla \bv, \bz)   - (\bv \cdot \nabla \bu, \bz)   
- (s , \nabla \cdot \bv)
&+\frac{1}{G}(\bz\cdot \btau_{\Gamma}, \bv\cdot \btau_{\Gamma})_{\Gamma}
+  g(\phi,\bv\cdot \bn_\Gamma)_{\Gamma}
=({\textbf{g}}_1,\bv)_{\Omega_f},\nonumber\\
\bK^{-1}(\bbeta,{\bw})&=-(\nabla\varrho,{\bw}){,}\nonumber\\
-  g( \bbeta , \nabla\psi) - g(\bz\cdot\bn_\Gamma , \psi)_{\Gamma} &= g(g_2,\psi)_{\Omega_p}.
\end{align*}}\normalsize%
\noindent 
Due to consistency, one can rewrite the above problem as follows:
{\small\begin{align*}
  a_{NS}(\underline{\bl},\underline{\btau}) -  b_{NS}(\bz,\underline{\btau}) &= 0,\\ 
   b_{NS}(\bv,\underline{\bl}) +  c_{NS}(\bv,\bz) + O_{NS}(\bu,\bv,\bz)+ O_{NS}(\bv,\bu,\bz)  
  + \frac{1}{G}(\bv\cdot\btau_{\Gamma},\bz\cdot\btau_{\Gamma})_{{\Gamma}}
&+d_{NS}(\bv,s)+  g(\bv\cdot{\bn_{\Gamma}},\varrho)_{{\Gamma}}  = (\textbf{g}_1,\bv)_{{\Omega}},\\
a_{D}({\bbeta},\bw)-b_{D}(\varrho,\bw)&=0,\\
b_{D}(\psi,{\bbeta})+ c_{D}(\psi,\varrho)-g(\psi,\bz\cdot\bn_{\Gamma})_{{\Gamma}}&= g(g_2,\psi_h)_{{\Omega}}.
\end{align*}}\normalsize%
\noindent Allowing $\textbf{g}_1=\bu-\bO_{fh}\bu$, $g_2=O_{ph}\phi-\phi$, we put
$\underline{\btau}=\underline{\bs} - \bOa\underline{\bs},\,\,
  \bv=\bu - \bOb \bu,\,\, {\bw}=  {\bomega} - \bOc{\bomega} ,\,\, \psi=\Od\phi -\phi$ in the above equations. We then subtract the resulting equations from (\ref{Oseen1})-(\ref{Oseen4}) after setting 
  ${\underline{\btau}_h}=\bPa\underline{\bl}$, $\bv_h=\bPb\bz$, ${\bw}_h=\bPc\bbeta$, and $\psi_h=-\Pd\varrho$ to obtain
 \small{ \begin{align}
  &\|\bu-\bOb\bu\|^2_{L^2(\Omega_f)}+\|\Od\phi-\phi\|^2_{L^p(\Omega_p)} =  a_{NS}(\underline{\bs} - \bOa\underline{\bs},\underline{\bl}-\bPa\underline{\bl})-  b_{NS}(\underline{\bs} - \bOa\underline{\bs},\bz - \bPb\bz) \nonumber \\
 &\,\,\,+  b_{NS}(\underline{\bl}-\bPa\underline{\bl},\bu - \bOb \bu)+  c_{NS}(\bu - \bOb \bu,\bz - \bPb\bz) +d_{NS}(\bz - \bPb\bz,p-r_hp) -d_{NS}(\bu - \bOb \bu,s)\nonumber  \\
& \,\,\,\,\,+ O_{NS}(\bu - \bOb \bu,\bu,\bz- \bPb\bz)  + O_{NS}(\bu, \bu - \bOb \bu,\bz- \bPb\bz) + ({\frac{1}{G}}((\bu - \bOb \bu)\cdot\btau_{\Gamma}),(\bz-\bPb\bz)\cdot\btau_{\Gamma})_{{\Gamma}} \nonumber\\
  &\,\,\,\,\,+  g((\bu - \bOb \bu)\cdot{\bn_{\Gamma}},\varrho-\Pd\varrho)_{{\Gamma}} + g(\phi-\Od\phi,(\bz- \bPb\bz)\cdot\bn_{\Gamma})_{{\Gamma}}
+ a_{D}(\bomega - \bOc\bomega,\bPc{\bbeta}-{\bbeta})
  \nonumber\\
  &\,\,\,\,\,+ b_{D}(\bomega - \bOc\bomega,\varrho-\Pd\varrho)- gb_{D}(\phi-\Od\phi,\Pd{\bbeta}-{\bbeta})- g c_{D}(\phi-\Od\phi,\Pd\varrho-\varrho).\label{Add2}
  \end{align}}\normalsize%
\noindent Here, we have used the fact $[\bz]=0$ in the fifth term of the right side of (\ref{Add2}). Similar to the terms on the right-hand side of (\ref{Add1}), we can bound all the terms except the sixth one. Using the definition of $\bZ_h$ and using standard inequalities with the projection properties (\ref{projection1defn}) and (\ref{projection2}), one finds 
 \small{ \begin{align*}
&|  d_{NS}(\bu - \bOb \bu,s) |= | d_{NS}(\bu - \bOb \bu,s-r_hs)|= | d_{NS}(\bu - \bPb\bu,s-r_hs)+ (\bPb\bu-\bOb \bu,s-r_hs)| \nonumber\\
&\leq \frac{1}{16} \|\bu-\bOb\bu\|^2_{L^2(\Omega_f)}+ \frac{1}{10}\|\Od\phi-\phi\|^2_{L^p(\Omega_p)} + Ch^{2k_1+2}(|\bu|^2_{k_1+1,\Omega_f}+|p|^2_{k_1,\Omega_f})).
\end{align*}}\normalsize%
\noindent Incorporating the resulting inequalities in (\ref{Add2}), we complete the estimate (\ref{Oseenprojection1}).
Repeating the above set of arguments, we arrive at the estimate (\ref{Oseenprojection2}) involving $(\bu-\bOb\bu)_t$, $((\underline{\bs} - \bOa\underline{\bs})_t$, $(\phi-\Od\phi)_t$, and $(({\bomega} - \bOc {\bomega})_t$. This completes the proof.
\end{proof}
\vspace*{-5mm}
\section{Error Analysis}\label{s5}
\vspace*{-2mm}
This section deals with the error bounds of the fully discrete scheme. The stress tensor, NS velocity, Darcy velocity, and piezometric head errors are denoted by ${\underline{\e}^i_{\bs}}=\underline{\bs}-{ \underline{\bSig}_h^i}$, $\e_{\bu}^i=\bu-\bU_h^i$, $\e_{\bomega}^i=\bomega^i-{\bQ}_h^i$, and $e_{\phi}^i=\phi^i-{\Phi}_h^i$ for $i = 1,2,\cdots,N$, respectively.
Due to consistency of the scheme (\ref{FSLNSD1})-(\ref{FSLNSD5}), we observe that the weak solution $(\underline{\bs},\bu,p,\bomega,\phi)$ satisfies (\ref{FSLNSD2}) and (\ref{FSLNSD5}) with $\partial_t\bU_h^i$ and $\partial_t\Phi_h^i$ replaced by $\bu_t^i$ and $\phi_t^i$ at $t=t_i$. Therefore, subtracting (\ref{DFSLNSD1})-(\ref{DFSLNSD4}) from the resulting equations, we obtain the error equation such that
\small{\begin{align}
 & \quad\quad  a_{NS}({ \underline{\e}^i_{\bs}},{\underline{\btau}_h}) +  b_{NS}(\e_{\bu}^i,{ \underline{\btau}_h}) = 0,\label{EFSLNSD1}\\ 
  ( \partial_t\e_{\bu}^i,\bv_h) - & b_{NS}(\bv_h,{ \underline{\e}^i_{\bs}}) +   c_{NS}(\e_{\bu}^i,\bv_h)   
 + O_{NS}(\bu^i,\bu^i,\bv_h)
  - O_{NS}(\p\bU_h^{i-1},\bU_h^i,\bv_h)\nonumber\\
  + \frac{1}{G}(\e_{\bu}^i\cdot&\btau_{\Gamma},\bv\cdot\btau_{\Gamma})_{{\Gamma}}
+ g(e_{\phi}^i,\bv_h\cdot{\bn_{\Gamma}})_{{\Gamma}}+d_{NS}(\bv_h,p^i)  = (\partial_t\bu^i-\bu_t^i,\bv_h)_{{\Omega}_f},\label{EFSLNSD2}\\
&\quad \quad a_{D}(\e_{\bomega}^i,{\bw}_h)-b_{D}(e_{\phi}^i,{\bw}_h)=0,\label{EFSLNSD3}\\
gS_0(\partial_te_{\phi}^i,\psi_h) &+  b_{D}(\psi_h,\e_{\bomega}^i)+ c_{D}(e_{\phi}^i,\psi_h)-g(\e_{\bu}^i\cdot\bn_{\Gamma},\psi_h)_{{\Gamma}}= gS_0(\partial_t\phi^i-\phi_t^i,\psi_h)_{{\Omega}_p},\label{EFSLNSD4}
\end{align}}
for { all $ (\underline{\btau}_h,\bv_h,q_h,\bw_h)\in(\underline{\bY}_h,\bZ_h,M_{fh},\bW_h)$.}\par
\noindent We observe that the error due to the first-order approximation of the derivative for any function $\nu(t)$ is
\begin{align}
(\partial_t \nu^i,v)- ( \nu_t^i,v) = -\frac{1}{\Delta t} \int_{t_i-1}^{t_i} (s-t_{i-1})(\nu_{ss}(s),v)ds ,\label{Taylor}
\end{align}
where $\nu_t^i=\nu_t(x_i)$.
First, we split the errors as
\small{\begin{eqnarray*}
{\underline{\e}^i_{\bs}}&=\underline{\bs}-\textbf{P}_h^{\bs}\underline{\bs}^i+\textbf{P}_h^{\bs}\underline{\bs}^i- { \underline{\bSig}_h^i}=\bca+\bka,\quad
\e_{\bu}^i=\bu^i-\textbf{P}_h^{\bu}\bu^i-\textbf{P}_h^{\bu}\bu^i-\bU_h^i  =\bcb+\bkb,\\
\e_{\bomega}^i&=\bomega^i-\textbf{P}_h^{\bomega}\bomega^i+\textbf{P}_h^{\bomega}\bomega^i-{\bQ}_h^i = \bcc+\bkc,\quad
e_{\phi}^i=\phi^i-P_h^\phi\phi^i+P_h^\phi\phi^i -\Phi_h^i =\cd+\kd .
\end{eqnarray*}}
By means of the split in the errors, we rewrite the error equations (\ref{EFSLNSD1})-(\ref{EFSLNSD4}) as

\small{\begin{equation}
    \left. %
\begin{aligned}
 & \quad\quad  a_{NS}(\bka,{\underline{\btau}_h}) +  b_{NS}(\bkb,{ \underline{\btau}_h}) = -a_{NS}({ \underline{\boldsymbol{\zeta}}_{fh}^{\bs i}},{\underline{\btau}_h})- b_{NS}(\bcb,\underline{\btau}_h) ,\\ 
  ( \partial_t\bkb,\bv_h) - & b_{NS}(\bv_h,{\bka}) +   c_{NS}(\bkb,\bv_h)   
 + O_{NS}(\bu^i,\bu^i,\bv_h)
  - O_{NS}(\p\bU_h^{i-1},\bU_h^i,\bv_h)\\
  + \frac{1}{G}(\bkb\cdot&\btau_{\Gamma},\bv_h\cdot\btau_{\Gamma})_{{\Gamma}}
+ g(\Kd,\bv_h\cdot{\bn_{\Gamma}})_{{\Gamma}} = (\partial_t\bu^i-\bu_t^i,\bv_h)-( \partial_t\bcb,\bv_h)+b_{NS}(\bv_h,{\bca})\\
&\quad\quad -c_{NS}(\bcb,\bv_h) -\frac{1}{G}(\bcb\cdot\btau_{\Gamma},\bv_h\cdot\btau_{\Gamma})_{{\Gamma}}-g(\cd,\bv_h\cdot{\bn_{\Gamma}})_{{\Gamma}}-d_{NS}(\bv_h,p^i) ,\\
&\quad \quad a_{D}(\bkc,{\bw}_h)-b_{D}(\kd,{\bw}_h)=-a_{D}(\bcc,{\bw}_h)+b_{D}(\cd,{\bw}_h),\\
gS_0(\partial_t\kd,\psi_h) &+  b_{D}(\psi_h,\bkc)+ c_{D}(\kd,\psi_h)-g(\bkb\cdot\bn_{\Gamma},\psi_h)_{{\Gamma}}= gS_0(\partial_t\phi^i-\phi_t^i,\psi_h)\\
&\quad\quad-gS_0(\partial_t\cd,\psi_h)- b_{D}(\psi_h,\bcc)-c_{D}(\cd,\psi_h)+g(\bcb\cdot\bn_{\Gamma},\psi_h)_{{\Gamma}},\label{standard1}
\end{aligned}\right\}%
\end{equation}}
for { all $ (\underline{\btau}_h,\bv_h,q_h,\bw_h)\in(\underline{\bY}_h,\bZ_h,M_{fh},\bW_h)$.}\par
 \noindent Now, we are going to state one of the main results of our work. The next theorem gives the optimal error estimates for velocity, stress tensor, piezometric head, and Darcy velocity in the energy norm. Additionally, we take note of the fact that henceforth the positive constant $K_T$, which depends exponentially on time but is independent of $h$ and $E$, will be used instead of the constant $C$ in the subsequent analysis.
\begin{theorem}\label{Errorresult1}
Let the assumption $(\textbf{A}\textbf{1})$ and the hypothesis of Lemma \ref{apriori} be satisfied. Furthermore, assume that the initial values of the discrete velocity and piezometric head are $\bU_h^0=\bP_h\bu_0$ and $\Phi_h^0=L_h\phi_0$, respectively. Then, for $K_T > 0$, the errors satisfy
\small{ \begin{align*}
\|\e_{\bu}^n\|^2_{L^2(\Omega_f)} +  gS_0\|e_{\phi}^n\|^2_{L^2(\Omega_p)} + \mu  e^{-2\alpha t_n}\Delta t \sum_{i=1}^{n}& e^{2\alpha t_n}|({\underline{\e}^i_{\bs}},\e_{\bu}^i)|^2_{\mathcal{A}} +   gC_4 e^{-2\alpha t_n}\Delta t \sum_{i=1}^{n} e^{2\alpha t_n}|(\e_{\bQ}^i,e_{\phi}^i)|^2_{\mathcal{B}}\nonumber\\
+ e^{-2\alpha t_n} \frac{\Delta t}{G} \sum_{i=1}^{n} e^{2\alpha t_n}\|\e_{\bu}^i\cdot\btau\|^2_{L^2(\Gamma)} &\leq K_T(h^{2k_1}+h^{2k_2}+\Delta t^2).
\end{align*}}
\end{theorem}
\begin{proof}
\noindent We now choose 
${\underline{\btau_h} }=2\mu{\underline{\boldsymbol{\chi}}^{\bs i}}$, $\bv_h={\boldsymbol{\chi}^{\bu i}}$, ${\bw}_h=g{ \boldsymbol{\chi}^{\bomega i}}$, and $\psi_h={ \chi^{\phi i}}$, respectively, in (\ref{standard1}) and add the resulting equations. Next, we observe that
\small{\begin{align*}
(\partial_t \bkb,\bkb)\geq \frac{1}{2} \partial_t\| \bkb\|^2_{L^2(\Omega_f)}, \quad  (\partial_t \kd,\kd)\geq \frac{1}{2} \partial_t\| \kd\|^2_{L^2(\Omega_p)}.
\end{align*}}
\noindent Next, an application of (\ref{Taylor}) with standard inequalities leads to
{\small\begin{align*}
|(\partial_t\bu^i-\bu_t^i,\bkb)|&\leq \frac{\mu}{18} |(\bka,\bkb)|^2_{\mathcal{A}} + C\Delta t  \int_{t_{i-1}}^{t_i} \|\bu_{tt}\|^2_{L^2(\Omega_f)}  ds,\\
|(\partial_t\phi^i-\phi_t^i,\kd)|&\leq\frac{g}{10} |(\bkc,\kd)|^2_{\mathcal{B}}  + C\Delta t  \int_{t_{i-1}}^{t_i} \|\phi_{tt}\|^2_{L^2(\Omega_p)}  ds.
\end{align*}}
\noindent Since $ \partial_t \bcb   =  \frac{1}{\Delta t} \int_{t_{i-1}}^{t_i} \frac{\partial}{\partial t} (\bcb(s))ds$, applying (\ref{projection2}) and the standard inequalities, we obtain
\small{\begin{align*}
|A_3|&=\Big| \frac{1}{\Delta t} \int_{t_{i-1}}^{t_i} (\frac{\partial}{\partial t} (\bcb(s)),\bkb)ds\Big|
\leq \frac{\mu}{18} |(\bka,\bkb)|^2_{\mathcal{A}} +   C\frac{h^{2k_1}}{\Delta t} |\bu_t|^2_{L^2(t_{i-1},t_i;H^{k_1}(\Omega_f))},\\
|A_4|&=\Big| \frac{1}{\Delta t} \int_{t_{i-1}}^{t_i} (\frac{\partial}{\partial t} (\cd(s)),\kd)ds\Big|\leq \frac{g}{10} |(\bkc,\kd)|^2_{\mathcal{B}}  + C\frac{h^{2k_2}}{\Delta t} |\phi_t|^2_{L^2(t_{i-1},t_i;H^{k_2}(\Omega_p))}.
\end{align*}}
\noindent The bilinear terms and the interface terms can be estimated similarly, as shown in Lemma \ref{Opproperty}.\par
\noindent Incorporating all the above bounds, multiplying by $\Delta t e^{2\alpha t_i}$, and finally summing the resulting inequality from $i=1\,\,\text{to}\,\,n ~(\leq N)$, we arrive at
{\scriptsize\begin{align}
 &e^{2\alpha t_n} \|{\boldsymbol{\chi}}_{fh}^{\bu n}\|^2_{L^2(\Omega_f)}+  gS_0 e^{2\alpha t_n} \| \chi_{ph}^{\phi n}\|^2_{L^2(\Omega_p)} +  \mu  \Delta t \sum_{i=1}^{n}  e^{2\alpha t_i} |({\bka},\bkb)|^2_{\mathcal{A}} + g\Delta t \sum_{i=1}^{n}  e^{2\alpha t_i} |({\bkc},\kd)|^2_{\mathcal{B}}  +  \frac{\Delta t}{{G}}  \sum_{i=1}^{n} e^{2\alpha t_i}\|\bkb\cdot\btau_{\gamma}\|^2_{L^2(\Gamma)}\nonumber\\
&\,\,
\leq C {\Delta t}\big( \sum_{i=1}^{n}  e^{2\alpha t_i} \|\bkb\|^2_{L^2(\Omega_f)} +  gS_0 \sum_{i=1}^{n}  e^{2\alpha t_i}  \|\kd\|^2_{L^2(\Omega_p)} \big ) +\Delta t \sum_{i=1}^{n}  e^{2\alpha t_i}\bigg(\sum_{m=1}^{5} T_m\bigg) +  C(h^{2k_1}+h^{2k_2})e^{2\alpha t_n} + 
C {\Delta t}^2 e^{2\alpha t_{n+1}}\,.\label{EFNSED27}
\end{align}}
\noindent Note that the $T_1,\cdot\cdot\cdot,T_5$ appearing above are the trilinear terms, once we rewrite them as follows:
\small{\begin{align*}
& O_{NS}(\p\bU_h^{i-1},\bU_h^i,\bkb)-O_{NS}(\bu^i,\bu^i,\bkb)
=- O_{NS}(\p\bU_h^{i-1},\bkb,\bkb)- O_{NS}(\p\bU_h^{i-1},\bcb,\bkb)\\
 &
 + O_{NS}(\p\bU_h^{i-1}-\bar{\bU}_h^{i-1},\bu^i,\bkb)+ O_{NS}(\bar{\bU}_h^{i-1}-\bu^{i-1},\bu^i,\bkb)+  O_{NS}(\bu^{i-1}-\bu^i,\bu^i,\bkb) {= \sum_{m=1}^{5} T_m.} 
\end{align*}}
\noindent In obtaining (\ref{EFNSED27}), we have made use of the following facts: (i) Since $\bchi_{fh}^{\bu 0}=\bU_h^0-\textbf{P}_h^{\bu}\bu_0=0$, we have
\small{\begin{align}
\sum_{i=1}^{n} \Delta t e^{2\alpha t_i} \partial_t\|{ \boldsymbol{\chi}_{fh}^{\bu i}}\|^2_{L^2(\Omega_f)} &= e^{2\alpha t_n} \|{\boldsymbol{\chi}_{fh}^{\bu n}}\|^2_{L^2(\Omega_f)}
 - \sum_{i=1}^{n-1} (e^{2\alpha \Delta t}-1)  e^{2\alpha t_i} \|{ \boldsymbol{\chi}_{fh}^{\bu i}}\|^2_{L^2(\Omega_f)},\label{sum1}
\end{align}}
(ii) Since $\chi_{ph}^{\phi 0}=\Phi_h^0-P_h^{\phi}\phi_0=0$, a similar break for $\partial_t\|\kd\|^2_{L^2(\Omega_f)}$ as well. (iii) Also, we have applied the fact that $e^{2\alpha \Delta t}-1\leq C(\alpha)\Delta t$.\par
\noindent We now begin to bound the trilinear terms $T_m$'s.
 Using (\ref{non}), (\ref{nonliearity2}), (\ref{e1}) and (\ref{aprioricondition1}), the first trilinear term can be bounded as follows:
{\scriptsize\[\Delta t \sum_{i=1}^{n}  e^{2\alpha t_i}|T_1|\leq C_0\Delta t \sum_{i=1}^{n}  e^{2\alpha t_i}\|\p\bU_h^{i-1}\|_{\bX_h}\|\bkb\|_{\bX_h}^2 \leq \frac{\mu}{2} \Delta t \sum_{i=1}^{n}  e^{2\alpha t_i} |(\bka,\bkb)|^2_{\mathcal{A}} . \]}
\noindent For the next two terms, applying (\ref{Postprocess1}), (\ref{Postprocess2}), (\ref{ST8}), (\ref{trilinear1}), {(\ref{projection2})}, and Young's inequality, we obtain 
{\scriptsize\begin{align*}
\Delta t \sum_{i=1}^{n} & e^{2\alpha t_i}|T_2|
\leq \frac{\mu}{18}\Delta t \sum_{i=1}^{n}  e^{2\alpha t_i} |( {\bka},{\bkb})|^2_{\mathcal{A}}+ C\Delta t \sum_{i=1}^{n}  e^{2\alpha t_i}\|\bu^i\|^2_{2,\Omega_f}\|\boldsymbol{\chi}_{fh}^{\bu {(i-1)}}\|^2_{L^2(\Omega_f)}+C\bigg(\|\boldsymbol{\bzeta}_{fh}^{\bu (i-1)}\|_{\bX_h}^2\\
&\quad\quad\quad\quad\quad +\|\bu^{i-1}-\bar{\bU}_h^{i-1}\|_{\bX_h}+\|\bu^{i-1}\|^2_{2,\Omega_f}\bigg)\|\bcb\|_{\bX_h}^2,\\
\Delta t \sum_{i=1}^{n} & e^{2\alpha t_i}|T_4| \leq \frac{\mu}{18}\Delta t \sum_{i=1}^{n}  e^{2\alpha t_i} |(\bka,\bkb)|^2_{\mathcal{A}}+ C\Delta t \sum_{i=1}^{n}  e^{2\alpha t_i}\|\bar{\bU}_h^{i-1}-\bu^{i-1}\|^2_{\bX_h}.
\end{align*}}
\noindent Applying (\ref{brokenpoincare}) and (\ref{trilinear2}) along with standard inequalities, we find
\begin{equation*}\scalebox{0.88}{$
 \Delta t \sum_{i=1}^{n}  e^{2\alpha t_i}|T_3| 
  \leq \frac{\mu}{18}\Delta t \sum_{i=1}^{n}  e^{2\alpha t_i} |({\bka},\bkb)|^2_{\mathcal{A}}+ C\Delta t \sum_{i=1}^{n}  e^{2\alpha t_i}\|\bu^i\|^2_{2,\Omega_f}\Big(\|{\boldsymbol{\chi}}_{fh}^{\bu (i-1)}\|^2_{L^2(\Omega_f)}
  +\|\boldsymbol{\bzeta}_{fh}^{\bu (i-1)}\|^2_{\bX_h}+\|\bar{\bU}_h^{i-1}-\bu^{i-1}\|^2_{\bX_h}\Big).$}
  \end{equation*}
\noindent To bound $T_5$, we use $\bu^i-\bu^{i-1}=\int_{t_{i-1}}^{t_i}\bu_t(t)\,dt$ and standard inequalities to arrive at
\small{ \begin{align*}
 \Delta t \sum_{i=1}^{n}  e^{2\alpha t_i}|T_5|
 \leq\frac{\mu}{18}\Delta t \sum_{i=1}^{n}  e^{2\alpha t_i} |({\bka},\bkb)|^2_{\mathcal{A}} + C \Delta t^2 \sum_{i=1}^{n}  e^{2\alpha t_i}\|\bu^i\|^2_{2,\Omega_f}\int_{t_{i-1}}^{t_i}\|\bu_t\|_{L^2(\Omega_f)}^2.
 \end{align*}}
\noindent Now, applying assumption $(\textbf{A}\textbf{1})$ for the terms appearing in $T_2$, $T_5$, etc. can be estimated as

{\small\begin{align*}
&C {\Delta t}^2 \sum_{i=1}^{n}  e^{2\alpha t_i}\|\bu^i\|^2_{2,\Omega_f} \int_{t_{i-1}}^{t_i} \|\bu_t\|^2_{L^2(\Omega_f)} dt \leq C {\Delta t}^2 \sum_{i=1}^{n}   \int_{t_{i-1}}^{t_i} e^{2\alpha (t_i-s)} e^{2\alpha s}\|\bu_{t}\|^2_{L^2(\Omega_f)} dt 
\leq C {\Delta t}^2 e^{2\alpha t_{n+1} },\\
&C\Delta t \sum_{i=1}^{n}  e^{2\alpha t_i} \|\bcbI\|^2_{\bX_1}\|\bcb\|^2_{\bX_1} \leq Ch^{4k_1+4k_2} \Delta t \sum_{i=1}^{n}  e^{2\alpha t_i} \leq Ch^{4k_1+4k_2} \Delta t e^{2\alpha\Delta t}\dfrac{(e^{2\alpha n\Delta t}-1)}{(e^{2\alpha\Delta t}-1)}
 \leq Ch^{4k_1+4k_2} e^{2\alpha t_{n+1}}.
\end{align*}}%
\noindent Finally, collecting the nonlinear bounds, applying the discrete Gr\"{o}nwall's inequality, multiplying the resulting equation by $e^{-2\alpha t_n}$, using the triangle inequality and approximation properties (\ref{projection3}), (\ref{projection2}), (\ref{N5}), and (\ref{N6}) along with assumption $(\textbf{A1})$, the desired estimate is achieved. This concludes the proof. 
 \end{proof}
\subsection{Error Estimates for the
Velocity and Piezometric Head in $L^{\infty}(\textbf{L}^2)$-norm} 
\noindent Our next goal is to derive optimal estimates in the $L^\infty(\bL^2)$ and $L^\infty(L^2)$-norms, respectively. A standard approach to derive these estimates for NSEs is to
perform integration by parts to handle the nonlinear term. But after performing the mechanism, the interface terms will appear, which will create difficulty in recovering the desired optimal results. Therefore, instead of using standard operators, we will use modified Oseen operators henceforth. For our subsequent analysis, we split the errors as
{\small\begin{align*}
 {\underline{\e}^i_{\bs}}&=\underline{\bs}^i-\bOa\underline{\bs}^i+\bOa\underline{\bs}^i-{ \underline{\bSig}_h^i}=\bTsi+\bups, \quad
  \e_{\bu}^i=\bu^i-\bOb\bu^i+\bOb\bu^i-\bU_h^i=\bTui+\bupu,\\
   \e_{\bomega}^i&=\bomega^i-\bOc\bomega^i+\bOc\bomega^i-\bomega_h^i=\bTqi+\bupq,\quad
    e_{\phi}^i=\phi-\Od\phi^i+\Od\phi^i-\Phi_h^i=\Tpi+\upp.  
 \end{align*}}
\noindent By means of the split in the errors, we rewrite the error equations (\ref{EFSLNSD1})-(\ref{EFSLNSD4}) as
 {\scriptsize\begin{equation}
    \left. %
\begin{aligned}
 & \quad\quad  a_{NS}(\bups,\underline{\btau}_h) +  b_{NS}(\bupu,{ \underline{\btau}_h}) = -a_{NS}(\bTsi,{\underline{\btau}_h})- b_{NS}(\bTui,\underline{\btau}_h) ,\\ 
  ( \partial_t\bupu,\bv_h) &-  b_{NS}(\bv_h,\bups) +   c_{NS}(\bupu,\bv_h)   
 + O_{NS}(\bu^i,\bu^i,\bv_h)
  - O_{NS}(\p\bU_h^{i-1},\bU_h^i,\bv_h)\\
  + \frac{1}{G}(\bupu\cdot\btau_{\Gamma},&\bv_h\cdot\btau_{\Gamma})_{{\Gamma}}
+ g(\upp,\bv_h\cdot{\bn_{\Gamma}})_{{\Gamma}} = (\partial_t\bu^i-\bu_t^i,\bv_h)-( \partial_t\bcb,\bv_h)+b_{NS}(\bv_h,\bTsi)\\
 &\quad-c_{NS}(\bTui,\bv_h) -\frac{1}{G}(\bTui\cdot\btau_{\Gamma},\bv_h\cdot\btau_{\Gamma})_{{\Gamma}}-g(\Tpi,\bv_h\cdot{\bn_{\Gamma}})_{{\Gamma}}-d_{NS}(\bv_h,p^i) ,\\
& \quad\quad a_{D}(\bupq,{\bw}_h)-b_{D}(\upp,{\bw}_h)=-a_{D}(\bTqi,{\bw}_h)+b_{D}(\Tpi,{\bw}_h),\\
gS_0(\partial_t\upp,\psi_h) &+  b_{D}(\psi_h,\bupq)+ c_{D}(\upp,\psi_h)-g(\bupu\cdot\bn_{\Gamma},\psi_h)_{{\Gamma}}= gS_0(\partial_t\phi^i-\phi_t^i,\psi_h)\\
&\quad\quad-gS_0(\partial_t\Tpi,\psi_h)- b_{D}(\psi_h,\bTqi)-c_{D}(\Tpi,\psi_h)+g(\bTui\cdot\bn_{\Gamma},\psi_h)_{{\Gamma}},\label{erroseen1}
\end{aligned}
\right\}
\end{equation}}
for all $ (\underline{\btau}_h,\bv_h,q_h,\bw_h)\in(\underline{\bY}_h,\bX_h,M_{fh},\bW_h)$.
\noindent Now we present here the optimal estimates for $\bups$ and $\Tpi$ in the $L^\infty(\bL^2)$ and $L^\infty(L^2)$-norms, respectively.
  \begin{lemma}\label{Errorresult2}
 Let the assumption $(\textbf{A}\textbf{1})$ and the hypothesis of Lemma \ref{apriori} be satisfied. Furthermore, assume that the initial values of the discrete velocity and piezometric head are $\bOb \bu^0=\bP^u_h\bu_0$ and $\Od\phi^0=P^{\phi}_h\phi_0$, respectively. Then, for $K_T > 0$, the errors $\bupu=\bOb\bu^i-\bU_h^i$ and $\upp=\Od\phi^i-\Phi_h^i$ satisfy
\small{\begin{align*}
\|\bUpun\|^2_{L^2(\Omega_f)} +  gS_0\|\Uppn\|^2_{L^2(\Omega_p)} + \mu  e^{-2\alpha t_n}&\Delta t \sum_{i=1}^{n} e^{2\alpha t_n}(\bups,\bupu)|^2_{\mathcal{A}} +   g e^{-2\alpha t_n}\Delta t \sum_{i=1}^{n} e^{2\alpha t_n}|(\bupq,\upp )|^2_{\mathcal{B}}\nonumber\\
+ e^{-2\alpha t_n} \frac{\Delta t}{{G}} \sum_{i=1}^{n} e^{2\alpha t_n}\|\bupu\cdot\btau\|^2_{L^2(\Gamma)} &\leq K_T(h^{2k_1+2}+h^{2k_2+2}+\Delta t^2).
\end{align*}}
\end{lemma}
\begin{proof}
We now use the definition of modified Oseen operators (\ref{Oseen1})-(\ref{Oseen4}) in (\ref{erroseen1}), add the resulting equations, and plug ${ \underline{\btau}_h}=2\mu{\bups}$, $\bv_h=\bupu$, ${\bw}_h=g\bupq$, and $\psi_h=\upp$. 
The bilinear terms can be bounded similarly as in Theorem \ref{Errorresult1}. But here we need $O(h^{k+1})$ for optimal results. Therefore, the next two terms can be bounded by applying (\ref{projection2}) and standard inequalities as follows:
{\scriptsize\begin{align*}
|(\partial_t\bTui,\bv_h)|&=\Big| \frac{1}{\Delta t} \int_{t_{i-1}}^{t_i} (\frac{\partial}{\partial t} (\bTui(s)),\bkb)ds\Big|
\leq \frac{\mu}{20} |(\bups,\bupu)|^2_{\mathcal{A}} +   C\frac{h^{2k_1+2}}{\Delta t} |\bu_t|^2_{k_1+1,\Omega_f}.\\
|(\partial_t\Tpi,\psi_h)|&=\Big| \frac{1}{\Delta t} \int_{t_{i-1}}^{t_i} (\frac{\partial}{\partial t} (\Tpi(s)),\kd)ds\Big|\leq \frac{g}{4} |(\bupq,\upp)|^2_{\mathcal{B}}  + C\frac{h^{2k_2+2}}{\Delta t} |\phi_t|^2_{k_1+1,\Omega_p}.
\end{align*}}\normalsize
\noindent 
Incorporating all the above bounds, multiplying by $\Delta t e^{2\alpha t_i}$, and finally summing the resulting inequality from $i=1\,\,\text{to}\,\,n ~(\leq N)$, we arrive at
{\scriptsize \begin{align}
 &e^{2\alpha t_n} \|\bUpun\|^2_{L^2(\Omega_f)}+ e^{2\alpha t_n}gS_0\|\Uppn\|^2_{L^2(\Omega_p)}+ \mu \Delta t \sum_{i=1}^{n} e^{2\alpha t_i}|(\bups,\bupu)|^2_{\mathcal{A}}+ \Delta t \sum_{i=1}^{n} e^{2\alpha t_i} |(\bupq,\upp)|^2_{\mathcal{B}}
  \nonumber\\
 &\quad+ \frac{\Delta t}{{G}}\sum_{i=1}^{n} e^{2\alpha t_i}\|\bupu\cdot\btau_{\Gamma}\|^2_{L^2(\Gamma)}\leq C {\Delta t}\Big( \sum_{i=1}^{n}  e^{2\alpha t_i} \|\bupu\|^2_{L^2(\Omega_f)} +  gS_0 \sum_{i=1}^{n}  e^{2\alpha t_i}  \|\upp\|^2_{L^2(\Omega_p)} \Big )\nonumber\\
 &\quad\quad\quad\quad\quad\quad\quad\quad\quad+ \Delta t \sum_{i=1}^{n}  e^{2\alpha t_i}\bigg(\sum_{i=1}^{10} B_i\bigg) +  C(h^{2k_1+2}+h^{2k_2+2})e^{2\alpha t_n} + 
C {\Delta t}^2 e^{2\alpha t_{n+1}}\,.\label{auxillary2}
 \end{align}} \normalsize
\noindent
{\xm Here, the terms $B_1,\cdot\cdot\cdot,B_{10}$ appearing above are the trilinear terms, which can be rewritten as follows:}
{\scriptsize\begin{align}
& O_{NS}(\p\bU_h^{i-1},\bU_h^i,\bupu)-O_{NS}(\bu^i,\bu^i,\bupu)+  O_{NS}(\bu^i,\bTui ,\bupu) + O_{NS}(\bTui,\bu^i,\bupu) = O_{NS}(\bu^{i}-\bu^{i-1},\bTu^i,\bupu)\nonumber\\
&\quad +O_{NS}(\bu^{i-1}-\bar{\bU}_h^{i-1},\bTui,\bupu)
+ O_{NS}({\bu}^{i-1}-{\bu}^{i},\bu^i,\bupu)+O_{NS}(\p\bU_h^{i-1}-\bar{\bU}_h^{i-1},\bu^i,\bupu)\nonumber\\
&\quad+ O_{NS}(\bar{\bU}_h^{i}-\bu^{i},\bu^i,\bupu)+ O_{NS}(\bTui,\bu^i,\bupu)-O_{NS}(\p\bU_h^{i-1},\bupu,\bupu)-O_{NS}(\p\bU_h^{i-1}-\bar{\bU}_h^{i-1},\bTui,\bupu)
 \nonumber\\
&\quad+ O_{NS}(\bar{\bU}_h^{i-1}-\bu^{i-1},\bu^i,\bupu)-O_{NS}(\bar{\bU}_h^{i}-\bu^{i},\bu^i,\bupu)=\sum_{i=1}^{10}B_i.\label{rearrange1}
\end{align}}\normalsize
To obtain (\ref{auxillary2}), we have used similar steps as (\ref{sum1}) of Theorem \ref{Errorresult1} with the fact that $\bUpuz =\bOb \bu^0-\bP_h^{\bu}\bu^0=\textbf{0}$ and $\Uppz=\Od\phi^0-P_h^{\phi}\phi^0=0$.
 
\noindent Next, we will bound the trilinear terms. An application of (\ref{Postprocess1}), (\ref{Postprocess2}),(\ref{T8}), (\ref{e1}), (\ref{trilinear1}), (\ref{continuity}), (\ref{nonliearity2}), and standard inequalities leads to
{\scriptsize\begin{align*}
&\Delta t \sum_{i=1}^{n} e^{2\alpha t_i}|B_1|\leq \frac{\mu}{20}\Delta t \sum_{i=1}^{n} e^{2\alpha t_i}|(\bups,\bupu)|^2_{\mathcal{A}} + C\Delta t^2 \sum_{i=1}^{n} e^{2\alpha t_i}\| \bTui \|^2_{\bX_h}\bigg(\int_{t_{i-1}}^{t_i}\|\bu_t\|^2_{1,\Omega_f}\,ds\bigg),\\
&\Delta t \sum_{i=1}^{n} e^{2\alpha t_i}|B_2|\leq \frac{\mu}{20}\Delta t \sum_{i=1}^{n} e^{2\alpha t_i}|(\bups,\bupu)|^2_{\mathcal{A}} + C\Delta t \sum_{i=1}^{n} e^{2\alpha t_i} \|\bu^{i-1}-\bar{\bU}_h^{i-1}\|^2_{\bX_h}\| \bTui \|^2_{\bX_h},\\
&\Delta t \sum_{i=1}^{n} e^{2\alpha t_i}|B_3|\leq \frac{\mu}{20}\Delta t \sum_{i=1}^{n} e^{2\alpha t_i}|(\bups,\bupu)|^2_{\mathcal{A}} + C\Delta t^2 \sum_{i=1}^{n} e^{2\alpha t_i} \|\bu^i\|_{2,\Omega_f}\int_{t_{i-1}}^{t_i} \|\bu_{t}\|_{L^2(\Omega_f)}^2\,dt,\\
&\Delta t \sum_{i=1}^{n} e^{2\alpha t_i}|B_7|\leq \frac{\mu}{20}\Delta t \sum_{i=1}^{n} e^{2\alpha t_i}|(\bups,\bupu)|^2_{\mathcal{A}},\\
&\Delta t \sum_{i=1}^{n} e^{2\alpha t_i}|B_8|\leq \frac{\mu}{20}\Delta t \sum_{i=1}^{n} e^{2\alpha t_i}|(\bups,\bupu)|^2_{\mathcal{A}} + C\Delta t \sum_{i=1}^{n} e^{2\alpha t_i}\|\bu^i\|^2_{2,\Omega_f}\|\bUpui\|^2_{L^2(\Omega_f)}\\
&\quad\quad\quad\quad\quad\quad\quad\quad+C\Delta t \sum_{i=1}^{n} e^{2\alpha t_i}\bigg(\|\bTuI\|^2_{\bX_h}+ \|\bu^{i-1}-\bar{\bU}_h^{i-1}\|^2_{\bX_h}\bigg)\|\bTui\|^2_{\bX_h},\\
&\Delta t \sum_{i=1}^{n} e^{2\alpha t_i}|B_9-B_{10}| \leq \frac{\mu}{20}\Delta t \sum_{i=1}^{n} e^{2\alpha t_i}|(\bups,\bupu)|^2_{\mathcal{A}}  + C\Delta t^2 \sum_{i=1}^{n} e^{2\alpha t_i} \|\bu^i\|_{2,\Omega_f}\int_{t_{i-1}}^{t_i} \|(\bar{\bU}_h^{i}-\bu^{i})_{t}\|_{\bX_h}^2\,dt.
\end{align*}}
\noindent Applying (\ref{Postprocess1}), (\ref{trilinear1}), and the triangle inequality, 
the next term can be bounded as 
{\scriptsize\begin{align*}
   \Delta t \sum_{i=1}^{n} e^{2\alpha t_i} |B_4|  &\leq\frac{\mu}{20}\Delta t \sum_{i=1}^{n} e^{2\alpha t_i}|(\bups,\bupu)|^2_{\mathcal{A}}+C\Delta t \sum_{i=1}^{n} e^{2\alpha t_i}\|\bu^i\|_{2,\Omega_f}^2\bigg(\Delta t\big(\|(\bU_h^i-\bu^i)_t\|^2_{\bX_h}\\
   &\quad+\|(\bu^i-\bar{\bU}_h^{i})_t\|^2_{\bX_h}\big)+\|\bUpui\|^2_{L^2(\Omega_f)}+\|\bTuI\|^2_{L^2(\Omega_f)}\bigg).
    \end{align*}}
Moreover, a use of (\ref{trilinear2}) and Young's inequality yields {\scriptsize\begin{align*}
 \Delta t \sum_{i=1}^{n} &e^{2\alpha t_i} |B_5+B_6| 
  \leq \frac{\mu}{20}\Delta t \sum_{i=1}^{n} e^{2\alpha t_i}|(\bups,\bupu)|^2_{\mathcal{A}}+C\Delta t \sum_{i=1}^{n} e^{2\alpha t_i}\|\bu^i\|_{2,\Omega_f}^2\bigg(\|\bar{\bU}_h^{i}-\bu^i\|^2_{L^2(\Omega_f)}+\|\bTui\|^2_{L^2(\Omega_f)}\bigg).
    \end{align*}}
   \noindent Also, we notice that $ \|(\bar{\bU}_h^{i}-\bu^{i})_{t}\|_{\bX_h}\leq C $ from (\ref{Postprocess2}). 
\noindent 
 Finally, collecting the nonlinear bounds, applying the discrete Gr\"{o}nwall's inequality, multiplying by $e^{-2\alpha t_n}$, applying the triangle inequality, approximation properties (\ref{Oseenprojection1}), (\ref{Oseenprojection2}), and assumption $(\textbf{A}\textbf{1})$, we obtain the desired result. 
\end{proof}
\noindent The next theorem states optimal $L^\infty(\bL^2)$ and $L^\infty(L^2)$-norm results for velocity and piezometric head errors, respectively.
\begin{theorem} \label{velo1}
Let the hypothesis of Lemma \ref{Errorresult2} be satisfied. Then there exists a constant $K_T>0$ such that 
\small{ \begin{align*}
\|\bu^i-\bU_h^i\|^2_{L^2(\Omega_f)} +  gS_0\|\phi^i-\Phi_h^i\|^2_{L^2(\Omega_p)}\leq K_T(h^{2k_1+2}+h^{2k_2+2}+\Delta t^2). 
\end{align*}}
\end{theorem}
\begin{proof}
A combination of Lemmas \ref{Opproperty} and \ref{Errorresult2} completes the proof.
\end{proof}
\begin{Remark}
As demonstrated in Lemma \ref{Errorresult2}, all interface terms are wiped away by the definition of the modified Oseen operators. Consequently, they are not visible in the error equation. Therefore, the difficulties that arise regarding them, as discussed in the 'Introduction', have been taken care of. Apart from that, we have noticed in the above lemma that integration by parts has been avoided by means of the modified Oseen operators. Thus, instead of suboptimal estimates in $h$, we have taken contributions in terms of $\Delta t$ from the trilinear terms in (\ref{rearrange1}) and managed to obtain optimal results. 
\end{Remark}%
\subsection{Error Estimates for Pressure}
\noindent 
In order to obtain an optimal fully discrete pressure error bound in $L^2(L^2)$-norm, we first need to derive the error bounds for $\partial_t\e_f^i$ and $\partial_t e_p^i$ in $\bL^2(\bL^2)$ and $L^2(L^2)$-norms. These results are established in the following lemma. The main difficulty here is that, unlike earlier, the interface terms of the form $(\cdot,\cdot)_{\Gamma}$ do not get cancelled out. To overcome this technical challenge, we have applied (\ref{T12}) and summation by parts and achieved the desired result.
\begin{lemma}\label{energy}
 Suppose the assumptions of Theorem \ref{velo1} are satisfied. 
 Then the following estimate holds true:
\small{ \begin{align*}
&e^{-2\alpha t_n}\Delta t \sum_{i=1}^{n} e^{2\alpha t_i}\|\partial_t\e_{\bu}^i\|^2_{L^2(\Omega_f)} +  gS_0e^{-2\alpha t_n}\Delta t \sum_{i=1}^{n} e^{2\alpha t_i}\|\partial_t e_{\phi}^i\|^2_{L^2(\Omega_p)} + \mu  |(\underline{\e^n_{\bs}},\e_{\bu}^n)|^2_{\mathcal{A}} +   g |(\e_{{\bQ}}^n,e_{\phi}^n)|^2_{\mathcal{B}}\nonumber\\
&+ e^{-2\alpha t_n} \frac{\Delta t}{G} \sum_{i=1}^{n} e^{2\alpha t_n}\|\partial_t\e_{\bu}^i\cdot\btau\|^2_{L^2(\Gamma)} + \frac{1}{G}\|\e_{\bu}^n\cdot\btau\|^2_{L^2(\Gamma)} \leq K_T(h^{2k_1}+h^{2k_2}+\Delta t^2).
\end{align*}}
\end{lemma}
\begin{proof}
\noindent 
Consider the error equation (\ref{erroseen1}), subtract the first and third equations of (\ref{erroseen1}) at time $t=t_{i-1}$ from the one at time $t=t_{i}$. Insert the definition of modified Oseen operators (\ref{Oseen1})-(\ref{Oseen4}), add all the resulting equations, and choose ${ \underline{\btau}_h}=2\mu{\bups}$, $\bv_h=\partial_t\bupu$, ${\bw}_h=g\bupq$, $\psi_h=\partial_t\upp$. 
\noindent The bilinear terms on the right hand side of the resulting equation can be bounded following similar steps as in Lemma \ref{Errorresult2}.\par
\noindent Inserting all the bounds, multiplying by $\Delta t e^{2\alpha t_i}$ and summing from $i=1$ to $i=n(\leq N)$, we arrive at
{\scriptsize\begin{align}
&\Delta t\sum_{i=1}^{n} e^{2\alpha t_i}\Bigg(\|\partial_t\bupu\|_{L^2(\Omega_f)}+gS_0\|\partial_t\upp\|_{L^2(\Omega_p)}
 + 2\mu  a_{NS}(\partial_t{\bups},{\bups}) 
 + 2\mu  c_{NS}(\bupu,\partial_t\bupu) + ga_{D}(\partial_t\bupq,\bupq)\nonumber\\
&\,\,+ g c_{D}(\upp,\partial_t\upp) 
 + {\frac{1}{G}}(\bupu\cdot\btau_{\Gamma},\partial_t\bupu\cdot\btau_{\Gamma})_{{\Gamma}}  +  g(\upp,\partial_t\bupu\cdot\bn_{\Gamma})_{\Gamma}  -g(\bupu\cdot\bn_{\Gamma},\partial_t\upp)_{\Gamma}  \Bigg) 
\leq  C {\Delta t}^2 e^{2\alpha t_{n+1}}  \nonumber\\
 &\,\,\,\,+ C(h^{2k_1+2}+h^{2k_2+2})e^{2\alpha t_n} +\Delta t\sum_{i=1}^{n} e^{2\alpha t_i}\bigg(\frac{12}{28}  \|\partial_t\bupu\|_{L^2(\Omega_f)} + \frac{2gS_0}{4}  \|\partial_t\upp\|_{L^2(\Omega_p)} +C  \sum_{i=1}^{9}M_i\bigg) 
\,.\label{IEFNSED2.O}
\end{align}}
\noindent Note that $M_1,\cdot\cdot\cdot,M_9$ are the trilinear terms. We now rewrite the trilinear terms as follows:
{\scriptsize\begin{align}
 &O_{NS}(\p\bU_h^{i-1},\bU_h^i,\partial_t\bupu)-O_{NS}(\bu^i,\bu^i,\partial_t\bupu)+  O_{NS}(\bu^i,\bTui ,\partial_t\bupu) + O_{NS}(\bTui,\bu^i,\partial_t\bupu) \nonumber\\
&=   O_{NS}(\p\bU_h^{i-1}-\p\bar{\bU}_h^{i-1},\bTui,\partial_t\bupu) + O_{NS}(\bar{\bU}_h^{i-1}-{\bu}^{i-1},\e_u^i,\partial_t\bupu)-O_{NS}(\bu^{i-1}-\bu^i,\bu^i,\partial\bupu)\nonumber\\
&\quad + O_{NS}(\p\bU_h^{i-1}-\p\bar{\bU}_h^{i-1},\bupu,\partial_t\bupu)- O_{NS}(\p\bU_h^{i-1}-\p\bar{\bU}_h^{i-1},\bu^i,\partial_t\bupu)
+O_{NS}(\bTui,\bu^i,\partial_t\bupu)  \nonumber\\
 & \quad-  O_{NS}(\bar{\bU}_h^{i-1}-{\bu}^{i-1},\bu^i,\partial_t\bupu)
  + O_{NS}({\bu}^{i-1}-{\bu}^{i},\e_u^i,\partial_t\bupu)+ O_{NS}({\bu}^{i},\bupu,\partial_t\bupu).\label{rearrange2}
\end{align}}\normalsize%
\noindent
Similar to Theorem \ref{Errorresult2}, all the trilinear terms on the right hand side of (\ref{rearrange2}) can be bounded except the fourth term. 
Applying standard inequalities with (\ref{Postprocess1}) and (\ref{T6}), the term can be bounded as
{\scriptsize\begin{align*}
&\Delta t\sum_{i=1}^{n} e^{2\alpha t_i}| O_{NS}(\p\bU_h^{i-1}-\p\bar{\bU}_h^{i-1},\bupu,\partial_t\bupu)| 
\leq \frac{\mu}{28} \Delta t\sum_{i=1}^{n} e^{2\alpha t_i}\|\partial_t\bupu\|^2_{L^2(\Omega_f)}\\
&\quad\quad\quad+ \frac{1}{min h_E}\Delta t\sum_{i=1}^{n} e^{2\alpha t_i}\bigg( \|\bUpui\|^2_{\bX_h}+\|\bu^{i-1}-\bar{\bU}^{i-1}_h\|_{\bX_h}^2+\|\bTuI\|^2_{\bX_h}\bigg)\|\bupu\|^2_{\bX_h}.
\end{align*}}\normalsize%
\noindent For the terms on left hand side of (\ref{IEFNSED2.O}), we shift the discrete time derivative to the other arguments in the bilinear forms by using summation by parts. Since $a_{NS}(.,.)$ and $c_{NS}(.,.)$ are symmetric, we obtain
{\scriptsize\begin{align*}
&\Delta t \sum_{i=1}^{n}  e^{2\alpha t_i} \Big|    a_{NS}(\partial_t{\bups},{\bups}) 
 +   c_{NS}(\bups,\partial_t\bupu) \Big| = \frac{1}{2}   e^{2\alpha t_n} |  ({{\bUpsn}},\bUpun)|^2_{\mathcal{A}}  - \frac{1}{2} | ({{\bUpsz}},\bUpuz)|^2_{\mathcal{A}} \nonumber\\
& \quad\quad\quad  - \frac{1}{2}  \sum_{i=0}^{n-1} e^{2\alpha t_i} (e^{2\alpha \Delta t}-1)     |({{\bups}},\bupu)|^2_{\mathcal{A}} + \frac{{\Delta t}^2}{2}  \sum_{i=1}^{n} e^{2\alpha t_i} | (\partial_t{{\bups}},\partial_t\bUpu^i)|^2_{\mathcal{A}}.
\end{align*}}\normalsize%
\noindent Similar calculations will go for $a_D(\cdot,\cdot)$. Regarding the interface terms, we use the same strategy. For instance, the seventh term on the left hand side of (\ref{IEFNSED2.O}) can be estimated as
{\scriptsize\begin{align*}
&\Delta t \sum_{i=1}^{n}  \frac{1}{G} e^{2\alpha t_i} \Big|    (\bupu\cdot\btau_{\Gamma},\partial_t\bupu\cdot\btau_{\Gamma})_{{\Gamma}}  \Big| = \frac{1}{2G}  ( e^{2\alpha t_n} (\bUpun\cdot\btau_{\Gamma},\bUpun\cdot\btau_{\Gamma})_{{\Gamma}}   - (\bUpuz\cdot\btau_{\Gamma},\bUpuz\cdot\btau_{\Gamma})_{{\Gamma}} )  \nonumber\\
&\quad- \frac{1}{2G}  \sum_{i=0}^{n-1} e^{2\alpha t_i} (e^{2\alpha \Delta t}-1)   (\bupu\cdot\btau_{\Gamma},\bupu\cdot\btau_{\Gamma})_{{\Gamma}}   + \frac{{\Delta t}^2}{2G}  \sum_{i=1}^{n} e^{2\alpha t_i} (\partial_t\bupu\cdot\btau_{\Gamma},\partial_t\bupu\cdot\btau_{\Gamma})_{{\Gamma}}.
\end{align*}}
\noindent To bound the eighth term on the left hand side of (\ref{IEFNSED2.O}), we employ (\ref{T12}), (\ref{auxillary2}), and Young's inequality as follows:
{\scriptsize\begin{align*}
\Delta t \sum_{i=1}^{n}  e^{2\alpha t_i} \big|   g(\upp,\partial_t\bupu\cdot\bn_{\Gamma})_{\Gamma}\big| &\leq \frac{\Delta t}{28}  \sum_{i=1}^{n}  e^{2\alpha t_i} \|\partial_t\bupu\|^2_{L^2(\Omega_f)}  + C(h^{2k_1+2}+h^{2k_2+2}+{\Delta t}^2). 
\end{align*}}
\noindent However, to bound the ninth term, summation by parts has been utilized, since (\ref{T12}) is no longer applicable.
{\scriptsize\begin{align*}
&\Delta t \sum_{i=1}^{n}  e^{2\alpha t_i} \big|  g(\bupu\cdot\bn_{\Gamma},\partial_t\upp)_{{\Gamma}} \big| \leq g e^{2\alpha t_n}\big|(\bUpun\cdot\bn_{\gamma}, \Uppn)_{{\Gamma}}\big| + ge^{2\alpha t_0}\big|(\bUpuz\cdot\bn_{\gamma}, \Uppz)_{{\Gamma}} \big|\nonumber\\
&\quad \quad + g\sum_{i=1}^{n} e^{2\alpha t_{i-1}}(e^{2\alpha \Delta t}-1) \big|(\bUpui\cdot\bn_{\gamma}, \Uppi)_{{\Gamma}}\big| + g\Delta t \sum_{i=1}^{n}e^{2\alpha t_i} \big|(\partial_t\bupu\cdot\bn_{\gamma}, \Uppi)_{{\Gamma}} \big|\nonumber\\
&\qquad\quad \leq    \frac{\Delta t}{28}\sum_{i=1}^{n}  e^{2\alpha t_i} \|\partial_t \bUpun\|^2_{L^2(\Omega_f)} + \frac{gC_4}{2} e^{2\alpha t_n}\|\Upp \|^2_{M_{ph}} + C e^{2\alpha t_n} \frac{1}{2} \|\bUpun  \|^2_{L^2(\Omega_f)}\nonumber\\
 &\qquad\quad\quad\,\,\,\, + C \Delta t \sum_{i=1}^{n} ( e^{2\alpha t_{i-1}}+e^{2\alpha t_i}) \|\Uppi  \|^2_{M_{ph}} + C\Delta t \sum_{i=1}^{n}  e^{2\alpha t_i} \|\bUpui  \|^2_{\bX_h}.
 \end{align*}}\normalsize%
\noindent
Now, we insert all the above bounds in (\ref{IEFNSED2.O}), apply the discrete Gr\"{o}nwall inequality, (\ref{auxillary2}), Theorem \ref{Errorresult2}, and assumption $(\textbf{A1})$ to arrive at
{\scriptsize\begin{align*}
&\Delta t \sum_{i=1}^{n} e^{2\alpha t_i} \|\partial_t\bupu\|^2_{L^2(\Omega_f)}+gS_0\Delta t \sum_{i=1}^{n} e^{2\alpha t_i}\|\partial_t\upp\|^2_{L^2(\Omega_p)}+ \mu e^{2\alpha t_n} |({{\bUpsn}},\bUpun)|^2_{\mathcal{A}}+ g e^{2\alpha t_n}|(\bUpqn,\Uppn)|^2_{\mathcal{B}}\\
 &\quad\quad\quad
  + {\frac{\Delta t}{{G}}} \sum_{i=1}^{n} e^{2\alpha t_i}\|\partial_t\bupu\cdot\btau\|^2_{L^2(\Gamma)}+  e^{2\alpha t_n}{\frac{1}{{G}}}\|\bUpun\cdot\btau\|^2_{L^2(\Gamma)}\leq Ce^D(h^{2k_1}+h^{2k_2}+\Delta t^2),
 \end{align*}}\normalsize%
\noindent where $D=\frac{h^{2k_1}+(\Delta t)^2}{min h_E}$. In this case, we need that the mesh satisfies, for some constant $L>0$, independent of $h$ and $\Delta t$, $D\leq L$. Finally, a use of triangular inequality and Lemma \ref{Opproperty}, assumption (\textbf{A1}), and multiplication by $e^{-2\alpha t_n}$ gives our desired result.
\end{proof}
\begin{Remark}
We would like to mention that the restriction $D\leq L$ is milder than the quasi-uniformity assumption in space (see \cite[pp. 43, (109)]{CGPR17}).
\end{Remark}
\begin{Remark}\label{pressure}
In Lemma \ref{Errorresult2}, we have observed that the interface terms on the right hand side of the error equation are absorbed by modified Oseen operators, and one of the interface terms on the left hand side of that equation gives a norm, and the other two get cancelled out. However, in this lemma, the scenario is quite different, as the interface terms on the left hand side of the error equation do not cancel out. To bound the first one that is $(\upp,\partial_t\bupu\cdot\bn_{\Gamma})_{\Gamma}$, we employ (\ref{T12}). But for the other one, i.e., $(\bupu\cdot\bn_{\Gamma},\partial_t\upp)_{{\Gamma}}$, we can not follow a similar technique, as we need to bound $\|\partial_t \xi^i\|$ in the $L^2$-norm, not in the energy norm. Therefore, to overcome the difficulty, we have used summation by parts and achieved the desired result. 
\end{Remark}
\noindent The next theorem establishes the fully discrete error estimate for pressure, which can be derived utilizing the discrete
inf-sup condition in Lemmas \ref{infsup}, \ref{energy}, and Theorem \ref{Errorresult1}. 
\begin{theorem}
    Under the assumptions of Lemma \ref{energy}, the following hold true:
   {\small \begin{align*}
       e^{-2\alpha t_n}\Delta t \sum_{i=1}^{n} e^{2\alpha t_i}\|p^i-P_h^i\|_{L^2(\Omega_f)}\leq K_T(h^{2k_1}+h^{2k_2}+\Delta t^2).
    \end{align*}}
\end{theorem}
\noindent Next, we prove a lemma that plays a crucial role in obtaining pressure error estimates. Since we are now only in need of recovering $O(h^k)$ instead of $O(h^{k+1})$, we utilize standard projections here. Also, we want to avoid the additional technicalities that arise due to the appearance of two extra nonlinear terms from the definition of the modified Oseen operators.\par
\noindent The next lemma provides $L^{\infty}(\bL^2)$ and $L^{\infty}(L^2)$-norm estimates for $\partial_t \e_f^i$ and $\partial_t e_p^i$, respectively.
\begin{lemma}\label{auxilarylemma}
The errors $\partial_t\e_\sigma^i$, $\partial_t\e_u^i$, $\partial_t\e_q^i$, and $\partial_te_\phi^i$ satisfy the following bound:
{\small\begin{align*}
&\|\partial_t\e_\bu^n\|^2_{L^2(\Omega_f)} + gS_0e^{2\alpha t_n}\|\partial_te_\phi^n\|^2_{L^2(\Omega_p)} + \mu \Delta t e^{-2\alpha t_n}\sum_{i=1}^{n} e^{2\alpha t_i} ({\partial_t\underline{\e_\sigma^i}},\partial_t\e_\bu^i)|^2_{\mathcal{A} }+ \Delta t \sum_{i=1}^{n} e^{2\alpha t_i}\frac{1}{{G}}\|\partial_t{ \boldsymbol{\chi}_{fh}^{\bu i}}\cdot\btau\|^2_{L^2(\Gamma)}\nonumber\\
&\quad\quad+   g \Delta t e^{-2\alpha t_n}\sum_{i=1}^{n} e^{2\alpha t_i} |(\partial_t\e_q^i,\partial_t e_\phi^i )|^2_{\mathcal{B}} \leq K_T(h^{2k_1}+h^{2k_2}+\Delta t).
\end{align*}}
\end{lemma}
\begin{proof}
 Subtract (\ref{standard1}) at time $t=t_{i-1}$ from the one at time $t=t_{i}$, substitute  $\underline{\btau}=2\mu\bka$, $\bv_h=\partial_t\bkb$, ${\bw}_h=g\partial_t\bkc$, and $\psi_h=\partial_t\kd$ in the resulting equations, and add all those equations. 

\noindent Now, let us rewrite the trilinear terms as follows:
{\scriptsize\begin{align*}
&-O_{NS}(\bu^i,\bu^i,\partial_t{ \boldsymbol{\chi}_{fh}^{\bu i}}) + O_{NS}(\bu^{i-1},\bu^{i-1},\partial_t\bkb)+ O_{NS}(\mathfrak{P}\bU_h^{i-1},\bU_h^i,\partial_t\bkb)-O_{NS}(\mathfrak{P}\bU_h^{i-2},\bU_h^{i-1},\partial_t{\bkb})\\
&\,\,\, = \Delta t O_{NS}(\partial_t (\mathfrak{P}\bU_h^{i-1}-\mathfrak{P}\bar{\bU}_h^{i-1}),\bU_h^i,\partial_t\bkb) - \Delta t O_{NS}(\partial_t (\mathfrak{P}\bU_h^{i-2}-\mathfrak{P}\bar{\bU}_h^{i-2}),\partial_t\bU_h^i,\partial_t{\bkb}) \\
&\quad -\Delta t O_{NS}(\partial_t (\bar{\bU}_h^{i-1}-\bu^{i-1}),\e_\bu^i,\partial_t\bkb) -\Delta t O_{NS}( (\bar{\bU}_h^{i-2}-\bu^{i-2}),\partial_t\e_\bu^i,\partial_t\bkb) \\
&\quad +\Delta t O_{NS}( \partial_t\bu^{i-1},\bU_h^{i-1}-\bu^{i-1},\partial_t\bkb)  +\Delta t O_{NS}( \bu^{i-1}, \partial_t \e_u^{i},\partial_t\bkb)+\Delta t O_{NS}( (\bar{\bU}_h^{i-1}-\bu^{i-1}), \partial_t\e_u^i,\partial_t\bkb)  \\
&\quad 
  - \Delta t O_{NS}( \partial_t(\bu^{i-1}-\bu^{i}),\bu^{i-1},\partial_t\bkb)
- \Delta t O_{NS}( \partial_t (\bar{\bU}_h^{i-1}-\bu^{i-1}),\bu^{i-1},\partial_t\bkb). 
\end{align*}}\normalsize
\noindent To bound the first trilinear term, we perform integration by parts. Then, using the standard inequalities, (\ref{T9}), 
(\ref{nonlinearity8}), (\ref{Postprocess1}), (\ref{projection2}), and Lemma \ref{TL2}, we bound the terms as
{\scriptsize\begin{align*}
 &\big|\Delta t O_{NS}(\partial_t(\p\bU_h^{i-1}-\p\bar{\bU}_h^{i-1}),\bU_h^i,\partial_t\bkb)\big|
\leq \frac{\mu\Delta t}{32} |({\bka},{ \boldsymbol{\chi}_{fh}^{\bu i}})|^2_{\mathcal{A}} + \frac{\mu\Delta t}{2} t e^{-2\alpha \Delta t}\|\partial_t \bkb\|_{\bX_h}\\
&\quad\quad\qquad + \frac{\mu\Delta t}{2} e^{2\alpha \Delta t}\|\bU_h^i\|^4_{\bX_h}\|\partial_t{ \boldsymbol{\chi}_{fh}^{\bu i}}\|^2_{\bX_h}+C \Delta t \|\bU_h^i\|^2_{\bX_h}|(\partial_t\bka,\partial_t\bkb)|^2_{\mathcal{A}}+ Ch^{2k_1}\|\bU_h^i\|^2_{\bX_h}\int_{t_{i-2}}^{t_{i-1}} |\bu_t|^2_{k_1+1,\Omega_f} .
\end{align*}}\normalsize
\noindent The other trilinear terms can be bounded in a similar manner. For the bilinear terms, we have utilized Lemma \ref{allprojections} along with the standard inequalities. Now, incorporate the above bounds, multiply the resulting inequality by $e^{2\alpha t i}$, and sum over $2 \leq i \leq n\leq N$. Also, we use the fact that $\|\partial_t \bkb\|_{L^2(\Omega_f)} \leq \|\partial_t\e_{\bu}^i\|_{L^2(\Omega_f)} + \|\partial_t(\bu^i-\bP_h^{\bu} \bu^i)\|_{L^2(\Omega_f)}$ and $\|\bkb\|_{\bX_h}\leq \|\e_{\bu}^i\|_{\bX_h}+ \|\bu^i-\bP_h^{\bu}\bu^i\|_{\bX_h}$. Similar inequalities hold for $\|\kd\|^2_{L^2(\Omega_p)}$ and $\|\kd\|^2_{M_{ph}}$. Along with these bounds, using (\ref{stability}) and the discrete Gr\"{o}nwall inequality, we arrive at
{\scriptsize\begin{align*}
&e^{2\alpha t_n}\|\partial_t\bchi^n\|^2_{L^2(\Omega_f)} + gS_0e^{2\alpha t_n}\|\partial_t\xi^n\|^2_{L^2(\Omega_p)} + \mu \Delta t\bigg(\sum_{i=2}^{n} e^{2\alpha t_i} |({\partial_t\bKa^i},\partial_t{ \boldsymbol{\chi}_{fh}^{\bu i}})|^2_{\mathcal{A} }+   g  \sum_{i=2}^{n} e^{2\alpha t_i} |(\partial_t{\boldsymbol{\chi}_{ph}^{\bomega i}},\partial_t{ \chi_{ph}^{\phi i}} )|^2_{\mathcal{B}}\bigg)\nonumber\\
&\quad+ \Delta t \sum_{i=2}^{n} e^{2\alpha t_i}\frac{1}{{G}}\|\partial_t{ \boldsymbol{\chi}_{fh}^{\bu i}}\cdot\btau
\|^2_{L^2(\Gamma)}  \leq K_T(\|\partial_t \bkB\|^2_{L^2(\Omega_f)}+ \|\partial_t \kD\|^2_{L^2(\Omega_p)}+(h^{2k_1}+h^{2k_2}+\Delta t)).
\end{align*}}\normalsize
To bound $\partial_t \bkB$ and $\partial_t\kD$, use the definition of modified Oseen operators (\ref{Oseen1})-(\ref{Oseen4}) in (\ref{erroseen1}) with $i=1$, $\underline{\btau}=2\mu\bkA$, $\bv_h=\partial_t\bkB$, ${\bw}_h=g\partial_t\bkC$, $\psi_h=\partial_t\kD$. Then the terms in the resulting equation can be bounded by applying similar techniques used in Lemma \ref{energy} to achieve the desired estimate. 
\end{proof}
\noindent Our next aim is to establish the main theorem of the fully discrete $L^\infty(L^2)$-norm estimate for the pressure error. 
\begin{theorem}\label{pressurefinal1}
Under the assumptions of Theorem \ref{velo1}, there exists a positive constant $K_T$ such that the following estimate holds:
\small{\begin{align*}
\|p^i-P_h^i\|_{L^2(\Omega_f)} \leq K_T(h^{k_1}+h^{k_2}+{\Delta t}^{1/2} ).
\end{align*}}
\end{theorem} 
\begin{proof}
Due to consistency of the numerical scheme, we observe that the weak solution $(\bu,p,\phi)$ satisfies (\ref{DFSLNSD1})-(\ref{DFSLNSD4}) with $\partial_t\bU_h^i$ replaced by $\bu_t^i$ at $t=t_i$. Therefore, subtracting (\ref{DFSLNSD2}) from the resulting equation and choosing test functions $(\bv_h,0)$, we obtain the error equation as follows:
\small{\begin{align}
-&d_{NS} (\bv_h,r_hp^i-P_h^i)= (\partial_t\bu^i-\bu_t^i,\bv_h) - (\partial_t\e_u^i,\bv_h) -2\mu b_{NS}(\bv_h, \underline{\e_\sigma^i})
 + 2\mu c_{NS}(\e_\bu^i,\bv_h)+ d_{NS}(\bv_h,p^i-r_h p^i)\nonumber\\
&- {\frac{1}{G}}(\e_\bu^i\cdot\btau_{\Gamma},\bv_h\cdot\btau_{\Gamma})_{{\Gamma}} + g(e_\phi^i,\bv_h\cdot\bn)_{\Gamma} -O_{NS}(\bu^i,\bu^i,\bv_h)  +O_{NS}(\p\bU_h^{i-1},\bU_h^i,\bv_h).\label{F5}
\end{align}}
\noindent Here, the sixth term can be bounded using the standard inequalities with (\ref{T9}) and Lemma \ref{energy} as
\small{\begin{align*}
\Big| \frac{1}{G}(\e_{\bu}^i\cdot\btau_{\Gamma},\bv_h\cdot\btau_{\Gamma})_{{\Gamma}}\Big| &\leq \frac{1}{{G}}\|\e_{\bu}^i\cdot\btau_{\Gamma}\|_{L^2(\Gamma)}\|\bv_h\cdot\btau_{\Gamma}\|_{L^2(\Gamma)}\leq C \frac{1}{{G}}\|\e_{\bu}^i\|_{L^2(\Gamma)}\|\bv_h\|_{\bX_h}\leq  K_T(h^{k_1}+h^{k_2}+\Delta t)\|\bv_h\|_{\bX_h}.
\end{align*}} \normalsize
\noindent The other interface terms can be treated similarly.
Now, rewriting the trilinear terms and applying Lemma \ref{con}, we find that
{\scriptsize\begin{align*}
&O_{NS}(\p\bU_h^{i-1},\bU_h^i,\bv_h)-O_{NS}(\bu^i,\bu^i,\bv_h)  = O_{NS}(\bu^i,\e_u^i,\bv_h) + O_{NS}(\bar{\bU}_h^{i-1}-\bu^{i-1},\bU_h^i,\bv_h)\\
&\quad\quad\quad\quad\quad\quad\quad\quad\quad
+ O_{NS}(\bu^{i-1}-\bu^{i},\bU_h^i,\bv_h) + O_{NS}(\p\bU_h^{i-1}- \bar{\bU}_h^{i-1},\bU_h^i,\bv_h)\\
&\quad\quad \leq C (\|\bu^i\|_{1,\Omega_f}\|\e_f^i\|_{\bX_h} + \Delta t\|\bu_t\|_{L^{\infty}(0,T,H^1(\Omega_f))}\|\bU_h^{i}\|_{\bX_h} + \|\e_f^{i-1}\|_{\bX_h}\|\bU_h^{i}\|_{\bX_h}+\|\bar{\bU}_h^{i-1}-\bu^{i-1}\|_{\bX_h}\|\bU_h^{i}\|_{\bX_h})\|\bv_h\|_{\bX_h}.
\end{align*}}\normalsize
\noindent The remaining terms on the right hand side of (\ref{F5}) can be bounded like Lemmas \ref{energy} and \ref{auxilarylemma}. Thus, employing the inf-sup condition from Lemma \ref{infsup}, combining all the above estimates in (\ref{F5}), and applying Lemma \ref{energy}, we arrive at
\small{\begin{align*}
\|r_hp^i-P_h^i\|_{L^2(\Omega_f)} \leq  \|\partial_t \e_\bu^i\|_{L^2(\Omega_f)} + K_T(h^{k_1}+h^{k_2}+{\Delta t}^{1/2}). 
\end{align*} }
\noindent Since $\|p^i-P_h^i\|\leq\|p^i-r_hp^i\|_{L^2(\Omega_f)}+ \|r_hp^i-P_h^i\|_{L^2(\Omega_f)} $, using the approximation property (\ref{projection1}) and Lemma \ref{auxilarylemma}, we arrive at the desired result.
\end{proof}
\begin{Remark}
If we assume $\bu_{tt}\in {L^{\infty}(0,T,\bL^2(\Omega_f))}$, then we obtain an improved pressure error estimate:
\small{\begin{align*}
\|p^i-P_h^i\|_{L^2(\Omega_f)} \leq K_T(h^{k_1}+h^{k_2}+{\Delta t} ).
\end{align*}}
\end{Remark}
\begin{Remark}\label{3D}
We can easily extend the optimal results discussed in this study to the three dimensional case under certain conditions. The analysis would be modified according to trace inequalities \cite{KJ98, DDE11}, Sobolev embeddings \cite{LS03}, inverse inequalities \cite{DDE11}, etc. Mainly, the estimates of trilinear terms will be changed. To estimate the $L^6$-norm, the existing $L^\mathfrak{p}$ bounds for the three dimensional case \cite{LS03} can be used:
\small{\begin{align}
\|\bv\|_{L^\mathfrak{p}(\Omega_f)}\leq C_\mathfrak{p}\|\bv\|_{\bX_1}, \,\,\, \forall\bv\in\bX_1, \,\,\, \mathfrak{p}\in [2,6].\label{3D1}
\end{align}}
\noindent For the $L^3$-norm estimate, we need a modified bound for the elements in the discrete space $\bX_h$ (\cite[(25)]{KHR23}):
\small{\begin{align*}
\|\bv_h\|_{L^3(\Omega)}\leq C \|\bv_h\|^{1/2}_{L^2(\Omega)}\|\bv_h\|^{1/2}_{\bX_1},\quad \forall\bv_h\in\bX_h .
\end{align*}}
\noindent  
In the three dimensional case, we have to rely on the $L^3$ and $L^6$-norm bounds. So, using the aforementioned estimate, we can modify the lemma for the $L^3$-norm. For the $L^6$-norm, there is only an estimate (\ref{3D1}).
Consequently, we can estimate the nonlinear terms accordingly and derive our desired estimates in the three dimensional case.  
\end{Remark}
\vspace*{-4mm}
\section{Numerical Results}\label{s6}
This section provides the results of numerical simulations for the unsteady NSD model, offering empirical validation of the theoretical findings stated in Theorems \ref{Errorresult1}, \ref{velo1}, and \ref{pressurefinal1}. Spatial discretization is achieved using the LDG method, while temporal discretization employs the backward Euler scheme. To approximate the stress tensor $\underline{\bs}$, NS velocity $\bu$, pressure $p$, discontinuous $P_0-P_1-P_0$, and $P_1-P_2-P_1$ mixed finite element spaces are employed. For the approximation of the piezometric head $\phi$ and Darcy velocity ${\bomega}$, discontinuous $P_1-P_0$ and $P_2-P_1$ mixed finite element spaces are utilized. We choose the constants $\textbf{C}_{12}=0$ and $\textbf{D}_{12}=0$ \cite{CKSS02}. We set all parameters equal to $1$ for the first two examples (\cite{BDQ10,WLS18}).  
These examples are computed on the time interval $[0, 1]$ with the final time $T = 1$. Also, we have provided a model test in the last example.
We note that $|(\underline{\bs}-{ \underline{\bSig}_h^i},\bu-\bU_h^i)|_{\mathcal{A}}= |(\e^i_{\bs},\e_{\bu}^i)|_{\mathcal{A}}$ and $|(\bomega-{\bQ_h^i},\phi-\Phi_h^i)|_{\mathcal{A}}= |(\e^i_{\bomega},\e_{\phi}^i)|_{\mathcal{B}}$.

\noindent\textbf{Example 7.1} Suppose that the computational domain $\Omega$ is composed of $(0,1)\times(1,2)$ and $\Omega_p = (0,1)\times (0,1)$. The interface is $\Gamma=(0,1)\times \{1\} $.
 Choose the exact solution as
{\small\begin{align*}
     \bu &= [(x^2(y-1)^2+y)cos(t),(-\frac{2}{3}x(y-1)^3+(2-\pi sin(\pi x)))cos(t)],  \\
     p &=(2-\pi sin(\pi x))sin(0.5\pi y)cos(t), \quad
     \phi = (2-\pi sin(\pi x))(1-y-cos(\pi y))cos(t).
\end{align*}}
\noindent The initial condition and the source term are derived from the exact solution.  
	\begin{table}[ht!]
		\centering
		\caption{Numerical errors and convergence rates for Example $7.1$ with ${P}_0$--${P}_1$--${P}_0$ elements.}
		\label{tab:Table1}
		
\small{		\begin{tabular}{|c|c c| c c| c c|}
		
			\hline
			
			 $h$ & $|(\e^i_{\bs},\e_{\bu}^i)|_{\mathcal{A}} $ & Rate & $\|\bu^i-\bU^i_h\|_{L^2(\Omega)}$ & Rate & $\|p^i-P_h^i\|_{L^2(\Omega_f)}$ & Rate\\[1ex]
			
	\hline   1/2  & $2.0038e^{0}$  &          & $1.8170e^{-01}$ &        & $9.9704e^{-01}$ &     \\
			 1/4  & $1.0050e^{0}$  &  0.995  & $6.3951e^{-02}$ & 1.506 & $4.4052e^{-01}$ & -\\
			 1/8  & $ 4.9434e^{-01}$  &  1.02  & $1.2242e^{-02}$ & 1.969 & $3.3277e^{-01}$ & 0.575\\
			 1/16 & $2.8081e^{-01}$  &  0.815   & $3.1268e^{-03}$ & 2.008 & $2.2333^{-01}$ & 0.840 \\
			 1/32 & $1.3940e^{-01}$  &  1.010  & $7.5242e^{-04}$ & 2.038 & $1.1403e^{-01}$ & 0.969 \\
			 1/64 & $6.8166e^{-02}$  & 1.032   & $1.9094e^{-04}$ & 2.055 & $5.5664e^{-02}$ & 1.034 \\
			\hline

		\end{tabular}}


	\end{table}
\begin{table}[ht!]
		\centering
		\caption{Numerical errors and convergence rates for Example $7.1$ with ${P}_1$--${P}_0$ elements} 
		\label{tab:Table2}
		
	\small{	\begin{tabular}{|c|c c| c c| c c|}
		
			\hline
			
			 $h$ & $\|\phi-\Phi^i_h\|_{L^2(\Omega_p)}$ & Rate & $|(\e^i_{\bomega},\e_{\phi}^i)|_{\mathcal{B}}$ & Rate \\[1ex]
			
	\hline   1/2  & $7.6361e^{-02}$  &          & $1.4791e^{0}$ &        \\
			 1/4  & $2.4652e^{-02}$  &  1.631  & $ 7.3933e^{-01}$ & 1.000 \\
			 1/8  & $ 4.2891e^{-03}$  & 1.928   & $ 3.2385e^{-01}$ & 1.019 \\
			 1/16 & $1.1265e^{-03}$  &   1.966  & $1.5609e^{-01}$ & 0.987 \\
			 1/32 & $2.8828e^{-04}$  &  2.038  & $7.8699e^{-02}$ & 1.025 \\
			 1/64 & $7.0204e^{-05}$  &  2.052  & $3.8660e^{-02}$ & 1.054 \\
			\hline

		\end{tabular}}


	\end{table}
     \begin{table}[ht!]
		\centering
		\caption{Numerical errors and convergence rates for Example $7.1$ with ${P}_1$--${P}_2$--${P}_1$ elements.}
		\label{tab:Table3}
		
	\small{	\begin{tabular}{|c|c c| c c| c c|}
		
			\hline
			
			 $h$ & $|(\e^i_{\bs},\e_{\bu}^i)|_{\mathcal{A}}$ & Rate & $\|\bu^i-\bU^i_h\|_{L^2(\Omega)}$  & Rate & $\|p^i-P_h^i\|_{L^2(\Omega_f)}$ & Rate\\[1ex]
			
	\hline   1/2  & $1.2745e^{0}$  &          & $2.2755e^{-02}$ &        & $3.8440e^{-01}$ &     \\
			 1/4  & $3.2301e^{-01}$  &  1.947  & $3.4608e^{-03}$ & 2.693 & $1.1670e^{-01}$ & 1.719\\
			 1/8  & $ 8.3756e^{-02}$  &  1.980  & $4.9149e^{-04}$ & 2.817 & $3.0947e^{-02}$ & 1.914\\
			 1/16 & $1.8793e^{-02}$  &  2.015   & $6.7634e^{-05}$ & 2.915 & $7.6678^{-03}$ &  2.012 \\
			 1/32 & $4.6392e^{-03}$  & 2.018   & $1.1208e^{-05}$ & 2.986 & $1.8504e^{-03}$ & 2.050 \\
			\hline

		\end{tabular}}


	\end{table}
\begin{table}[ht!]
		\centering
		\caption{Numerical errors and convergence rates for Example $7.1$ with ${P}_2$--${P}_1$ elements} 
		\label{tab:Table4}
		
\small{		\begin{tabular}{|c|c c| c c| c c|}
		
			\hline
			
			 $h$ & $\|\phi-\Phi^i_h\|_{L^2(\Omega_p)}$ & Rate & $|(\e^i_{\bomega},\e_{\phi}^i)|_{\mathcal{B}}$ & Rate \\[1ex]
			
	\hline   1/2  & $3.5787e^{-02}$  &          & $2.1691e^{0}$ &        \\
			 1/4  & $4.1440e^{-03}$  &  2.982  & $ 5.4401e^{-01}$ & 1.962 \\
			 1/8  & $ 3.6195e^{-04}$  & 3.011   & $ 1.0989e^{-01}$ & 1.995 \\
			 1/16 & $3.9636e^{-05}$  &  3.019  & $2.5661e^{-02}$ & 2.030 \\
			 1/32 & $5.0149e^{-06}$  &  3.051  & $6.5821e^{-03}$ & 2.098 \\
			\hline

		\end{tabular} }


	\end{table}
    
     \noindent Tables \ref{tab:Table1}--\ref{tab:Table4} provide the computational errors and the order of convergence for discontinuous $P_0-P_1-P_0$ and $P_1-P_2-P_1$ elements for the fluid region, and discontinuous $P_1-P_0$ and $P_2-P_1$ elements for porous media, respectively. For both cases, parameters are taken as $C_{11}=10$ and $C_{11}=20$, respectively. The time step $\Delta t$ has been taken as $O(h)$ and $O(h^2)$ for $P_0-P_1-P_0$, $P_1-P_0$, and $P_1-P_2-P_1$, $P_2-P_1$, respectively. We notice that the computational results of Tables $1$-$4$ support our theoretically derived convergence rates.
    
 \noindent\textbf{Example 7.2} 
The computational domain and the interface condition are chosen similar to Example 7.1.
{\small\begin{align*}
     \bu &= [(exp(x+y)+y)exp(t),-(exp(x+y)+x)exp(t)],  \\
     p &= (cos(\pi x)cos(\pi y)+x)exp(t), \quad
     \phi = (exp(x+y) - cos(\pi x) + xy)exp(t).
\end{align*}}
\begin{table}[ht!]
		\centering
		\caption{Numerical errors and convergence rates for Example 7.2 with ${P}_0$--${P}_1$--${P}_0$ elements.}
		\label{tab:Table5}
		
\small{		\begin{tabular}{|c|c c| c c| c c|}
		
			\hline
			
			 $h$ & $ |(\e^i_{\bs},\e_{\bu}^i)|_{\mathcal{A}}$ & Rate & $\|\bu^i-\bU^i_h\|_{L^2(\Omega)}$ & Rate & $\|p^i-P_h^i\|_{L^2(\Omega_f)}$ & Rate\\[1ex]
			
	\hline   1/2  & $3.4711e^{0}$  &          & $1.1761e^{0}$ &        & $5.0668e^{0}$ &     \\
			 1/4  & $2.0074e^{0}$  &  0.981  & $2.9853e^{-01}$ & 1.939 & $2.4873e^{0}$ & 0.947\\
			 1/8  & $ 1.6772e^{0}$  & 1.014   & $6.4333e^{-02}$ & 1.978 & $9.7868e^{-01}$ & 0.971\\
			 1/16 & $1.1743e^{0}$  &   1.016   & $1.5886e^{-02}$ & 2.017 & $7.6099e^{-01}$ & 1.039 \\
			 1/32 & $5.7788e^{-01}$  & 1.018   & $4.1409e^{-03}$ & 2.021 & $3.8805e^{-01}$ & 1.045 \\
			 1/64 & $2.9260e^{-01}$  &  1.023  & $1.0195e^{-03}$ & 2.022 & $2.0119e^{-01}$ & 1.051 \\
			\hline

		\end{tabular}}
	\end{table}
\begin{table}[ht!]
		\centering
		\caption{Numerical errors and convergence rates for Example 7.2 with ${P}_1$--${P}_0$ elements.}
		\label{tab:Table6}
		
	\small{	\begin{tabular}{|c|c c| c c| c c|}
		
			\hline
			
			 $h$ & $\|\phi-\Phi^i_h\|_{L^2(\Omega_p)}$ & Rate & $|(\e^i_{\bomega},\e_{\phi}^i)|_{\mathcal{B}}$ & Rate\\[1ex]
			
	\hline   1/2  & $3.2114e^{-01}$  &          & $4.0484e^{0}$ &       \\     
			 1/4  & $6.3883e^{-02}$  &  1.974  & $ 1.8240e^{0}$ & 0.992 \\
			 1/8  & $ 1.3078e^{-02}$  &  1.982  & $ 7.5843e^{-01}$ & 1.007 \\
			 1/16 & $2.8530e^{-03}$  &   2.019  & $3.5484e^{-01}$ & 1.095 \\ 
			 1/32 & $7.2621e^{-04}$  & 2.028   & $1.7655e^{-01}$ & 1.015 \\ 
			 1/64 & $1.8377e^{-04}$  & 2.032   & $8.8759e^{-02}$ & 1.026\\ 
			\hline

		\end{tabular}}


	\end{table}
\begin{table}[ht!]
		\centering
		\caption{Numerical errors and convergence rates for Example 7.2 with ${P}_1$--${P}_2$--${P}_1$ elements.}
		\label{tab:Table7}
		
	\small{	\begin{tabular}{|c|c c| c c| c c|}
		
			\hline
			
			 $h$ & $ |(\e^i_{\bs},\e_{\bu}^i)|_{\mathcal{A}}$ & Rate & $\|\bu^i-\bU^i_h\|_{L^2(\Omega)}$ & Rate & $\|p^i-P_h^i\|_{L^2(\Omega_f)}$ & Rate\\[1ex]
			
	\hline   1/2  & $1.4977e^{0}$  &          & $5.4749e^{-02}$ &        & $1.5136e^{0}$ &     \\
			 1/4  & $3.4016e^{-01}$  &  1.963  & $5.8584e^{-03}$ & 3.017 & $8.9997e^{-01}$ & 2.083\\
			 1/8  & $ 7.2347e^{-02}$  &  2.013  & $6.4989e^{-04}$ & 3.022 & $1.1302e^{-01}$ & 2.085\\
			 1/16 & $1.6928e^{-02}$  &  2.023    & $7.8000e^{-05}$ & 3.047 & $1.4376e^{-02}$ & 2.097 \\
			 1/32 & $4.3415e^{-03}$  & 2.095   & $9.4340e^{-06}$ & 3.058 & $2.0123e^{-03}$ & 2.099 \\
			\hline

		\end{tabular}}
	\end{table}
\begin{table}[ht!]
		\centering
		\caption{Numerical errors and convergence rates for Example 7.2 with ${P}_2$--${P}_1$ elements.}
		\label{tab:Table8} 
		
	\small{	\begin{tabular}{|c|c c| c c| c c|}
		
			\hline
			
			 $h$ & $\|\phi-\Phi^i_h\|_{L^2(\Omega_p)}$ & Rate & $|(\e^i_{\bomega},\e_{\phi}^i)|_{\mathcal{B}}$ & Rate\\[1ex]
			
	\hline   1/2  & $4.2570e^{-02}$  &          & $7.4881e^{-01}$ &       \\     
			 1/4  & $5.9015e^{-03}$  &  2.850  & $ 2.3889e^{-01}$ & 1.648 \\
			 1/8  & $ 7.1122e^{-03}$  & 3.002   & $ 3.9794e^{-02}$ & 1.960 \\
			 1/16 & $8.8777e^{-05}$  &   3.010  & $9.9782e^{-03}$ & 1.995 \\ 
			 1/32 & $1.1018e^{-05}$  & 3.052   & $2.5640e^{-03}$ & 2.058 \\ 
			\hline

		\end{tabular}}


	\end{table}
	\noindent Tables \ref{tab:Table5}--\ref{tab:Table8} indicate the computational errors and the order of convergence for discontinuous $P_0-P_1-P_0$ and $P_1-P_2-P_1$ elements for the fluid region, and discontinuous $P_1-P_0$ and $P_2-P_1$ elements for porous media, respectively. The parameters are chosen as $C_{11}=10$ and $C_{11}=20$ for both cases. The time step has been chosen similar to Example 7.1.  
We observe that these results are consistent with the theoretically derived results. 

 \begin{figure}[ht!]
  \centering
  \begin{subfigure}[b]{0.48\textwidth}
    \centering
    \includegraphics[width=\textwidth]{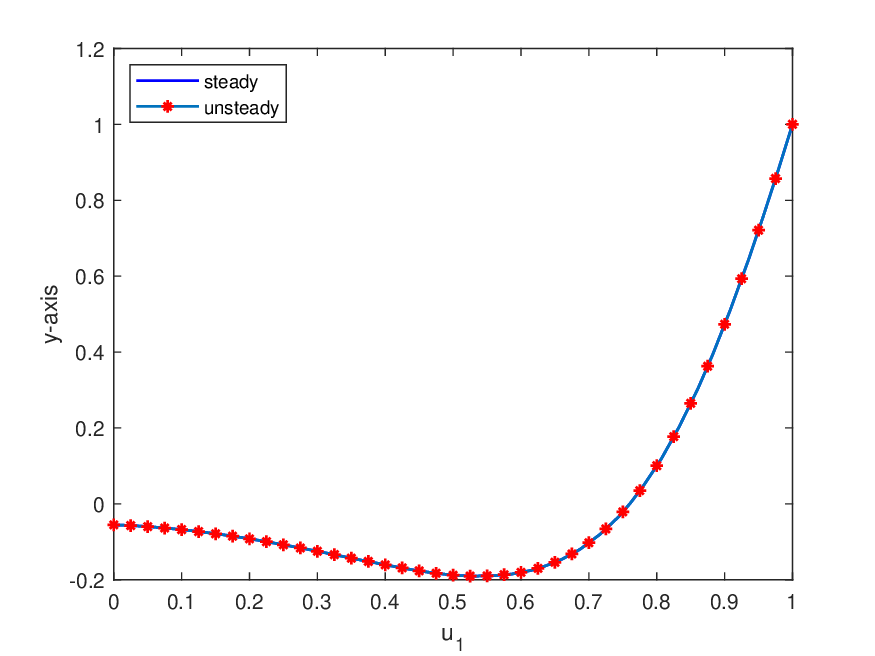}
    \caption{First velocity component along the line $x=0.5$.}
  \end{subfigure}
 \hfill
  \begin{subfigure}[b]{0.48\textwidth}
    \centering
    \includegraphics[width=\textwidth]{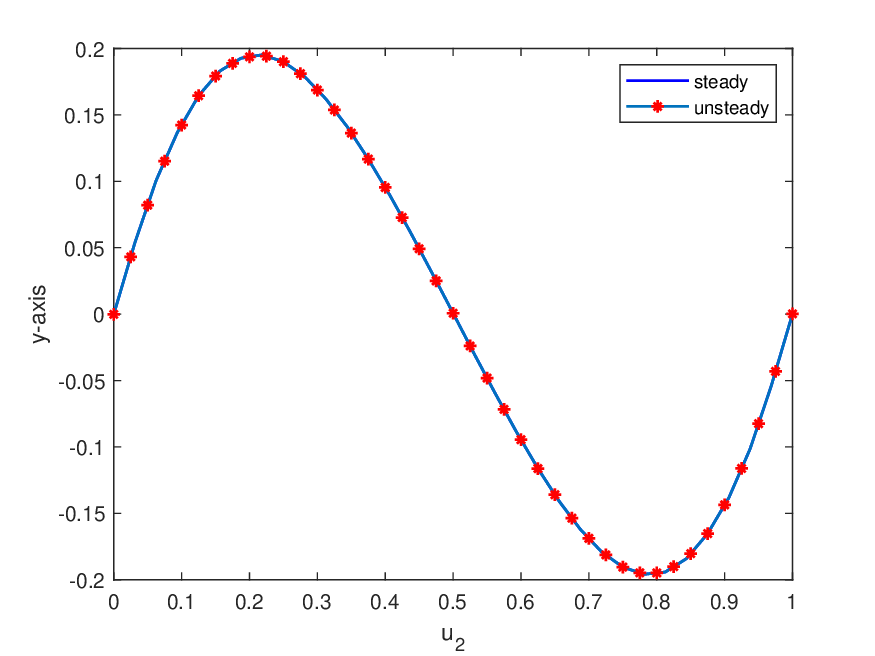}
    \caption{Second velocity component along the line $y=1.5$.}
    \label{fig:second}
  \end{subfigure}
  \centering
  \begin{subfigure}[b]{0.48\textwidth}
    \centering
    \includegraphics[width=\textwidth]{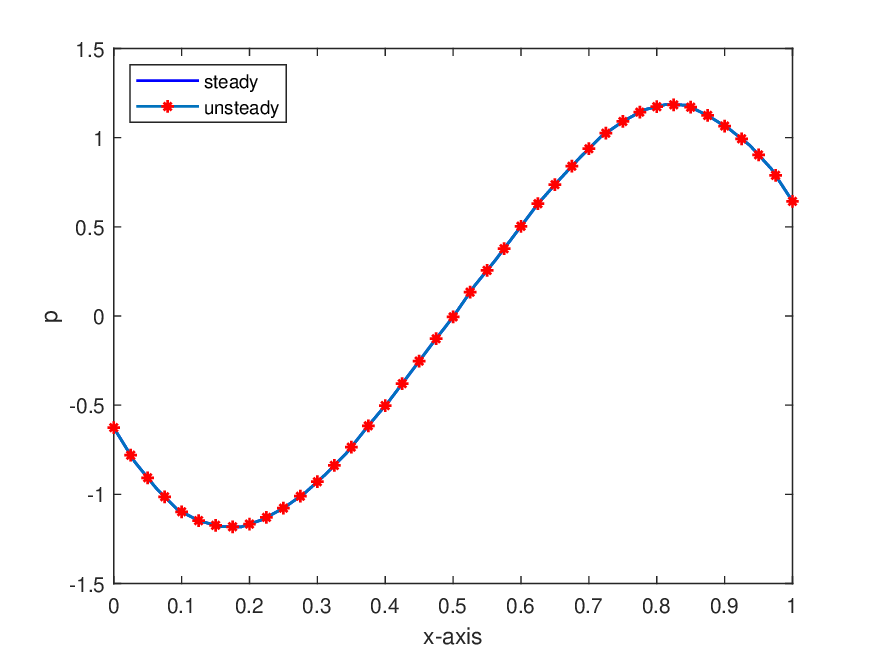}
    \caption{Pressure along the line $y=1.5$.}
  \end{subfigure}
 \hfill
  \begin{subfigure}[b]{0.48\textwidth}
    \centering
    \includegraphics[width=\textwidth]{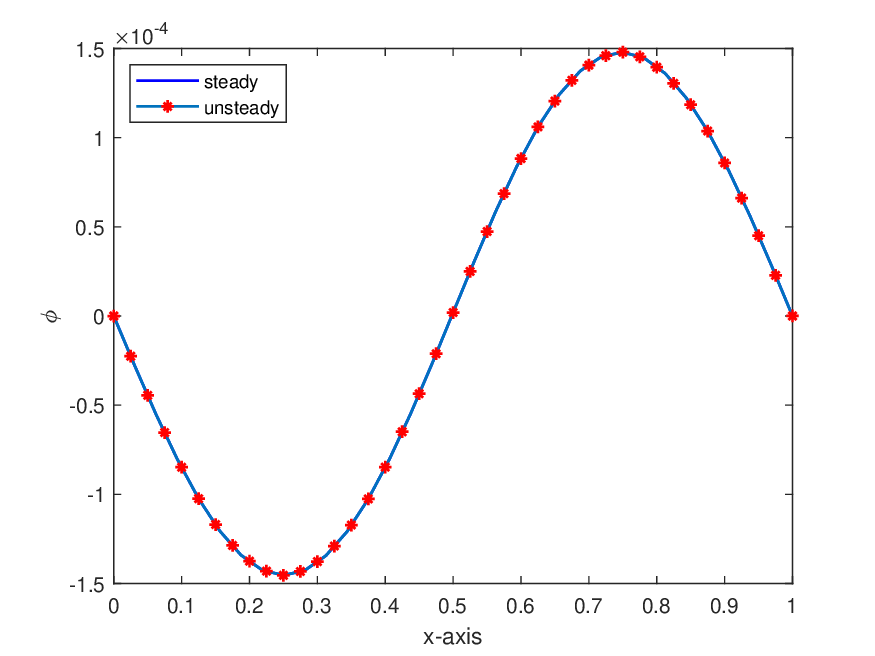}
    \caption{Piezometric head along the line $y=0.5$.}
    \label{fig:second1}
  \end{subfigure}
  \caption{Velocity components, pressure, and piezometric head for lid-driven cavity flow with $\mu=1$.}
  \label{fig:all}
\end{figure}
\begin{figure}[!ht]
  \centering
  \begin{subfigure}[b]{0.48\textwidth}
    \centering
\includegraphics[width=\textwidth]{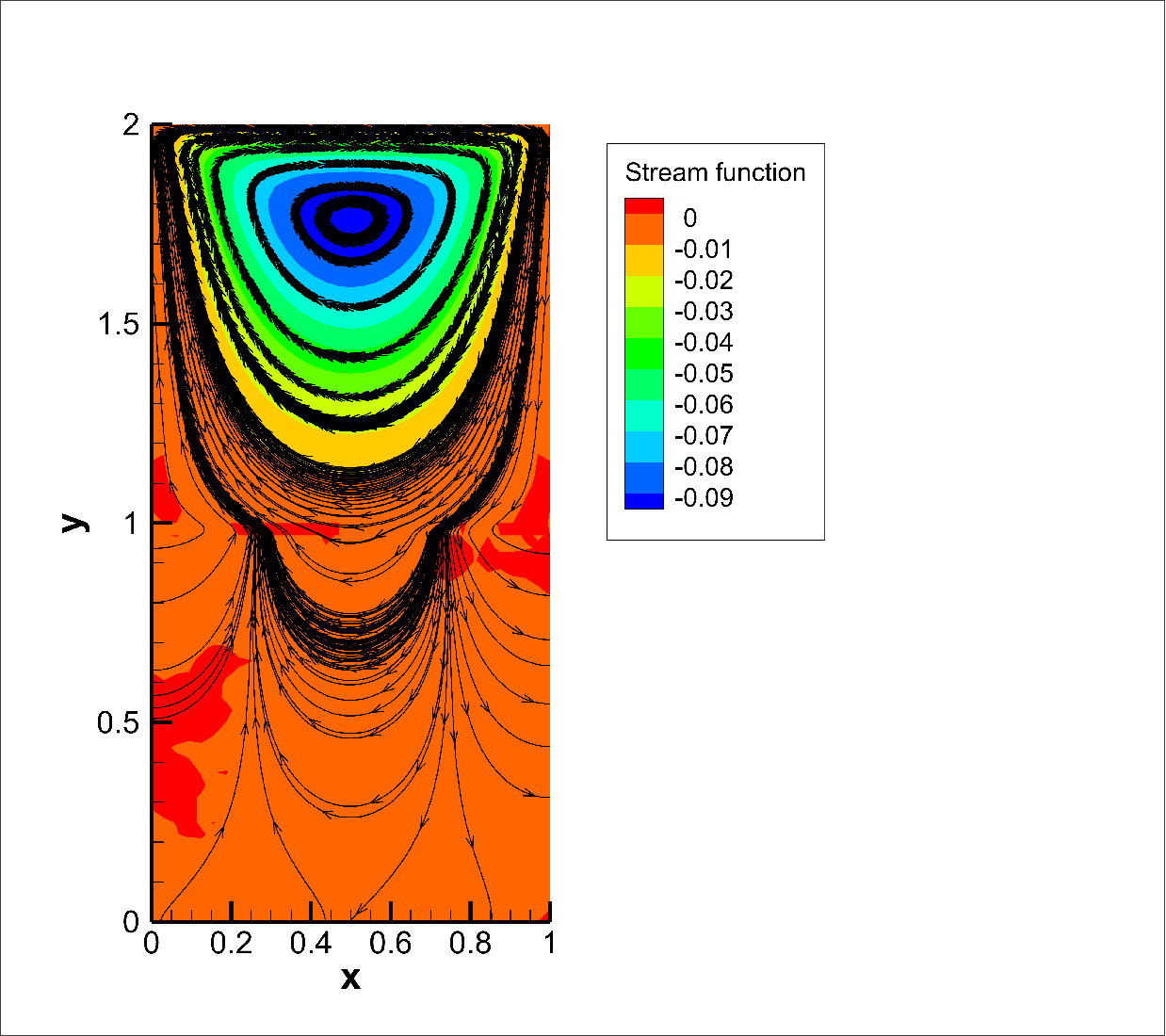}
\caption{$\mu = 1$, $\bK = 1$}\label{thelabel}
\end{subfigure}
\hfill
  \begin{subfigure}[b]{0.48\textwidth}
    \centering
\includegraphics[width=\textwidth]{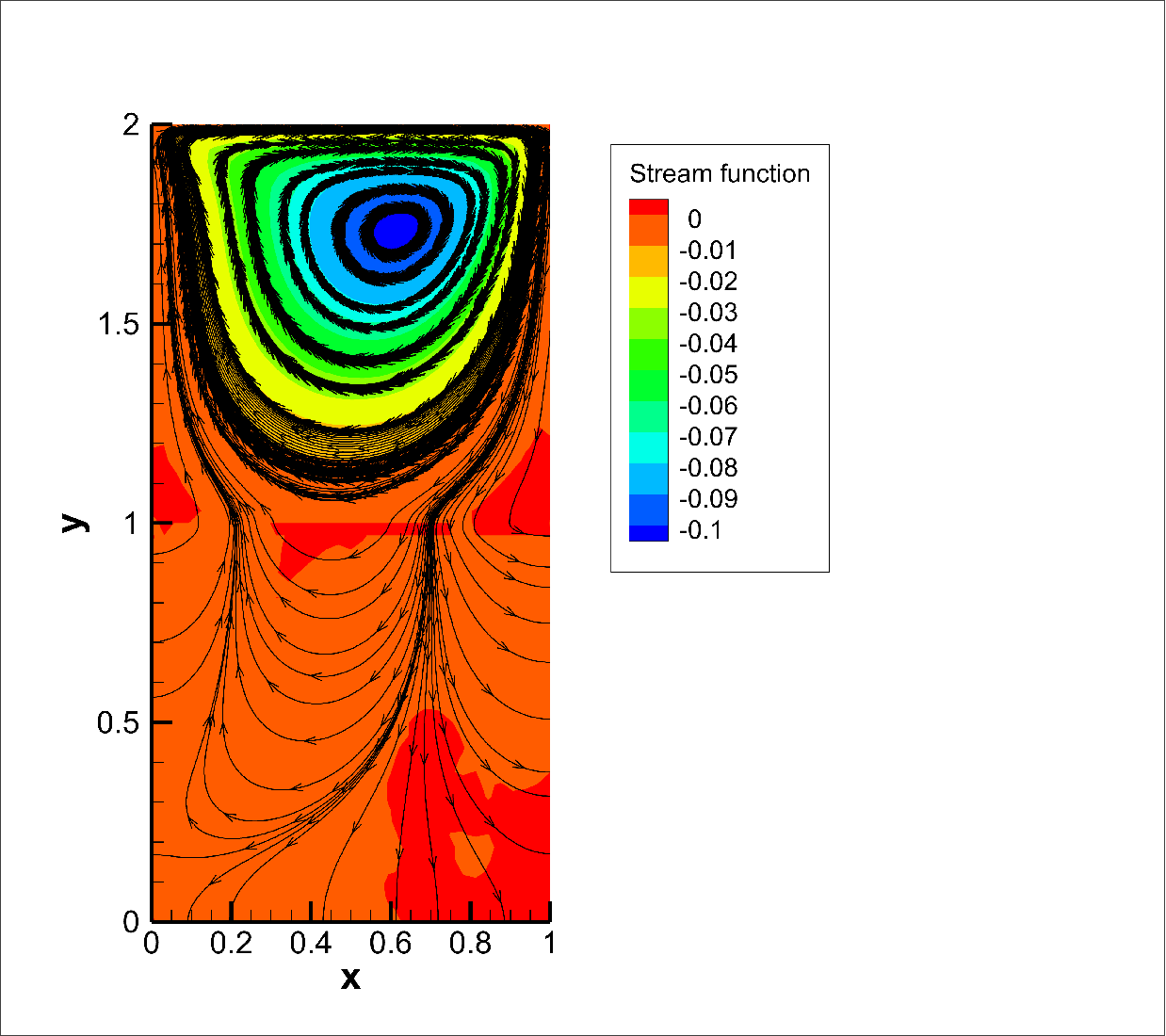} 
\caption{$\mu = 1/100$, $\bK = 1$}\label{thelabel1}
\end{subfigure}
\hfill
  \begin{subfigure}[b]{0.48\textwidth}
    \centering
\includegraphics[width=\textwidth]{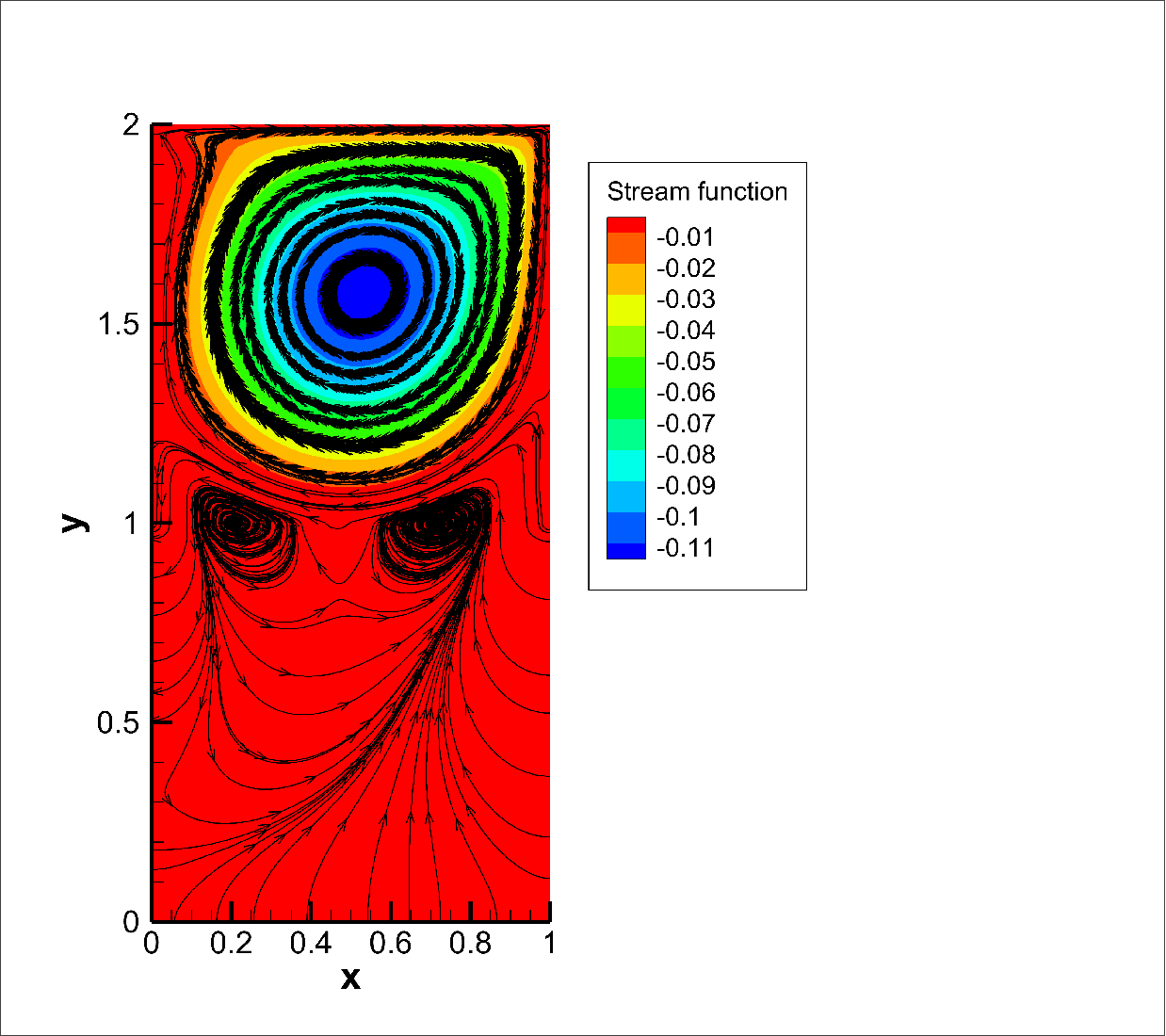} 
\caption{$\mu = 1/1000$, $\bK = 1$}\label{thelabel3}
\end{subfigure}
\hfill
  \begin{subfigure}[b]{0.48\textwidth}
    \centering
\includegraphics[width=\textwidth]{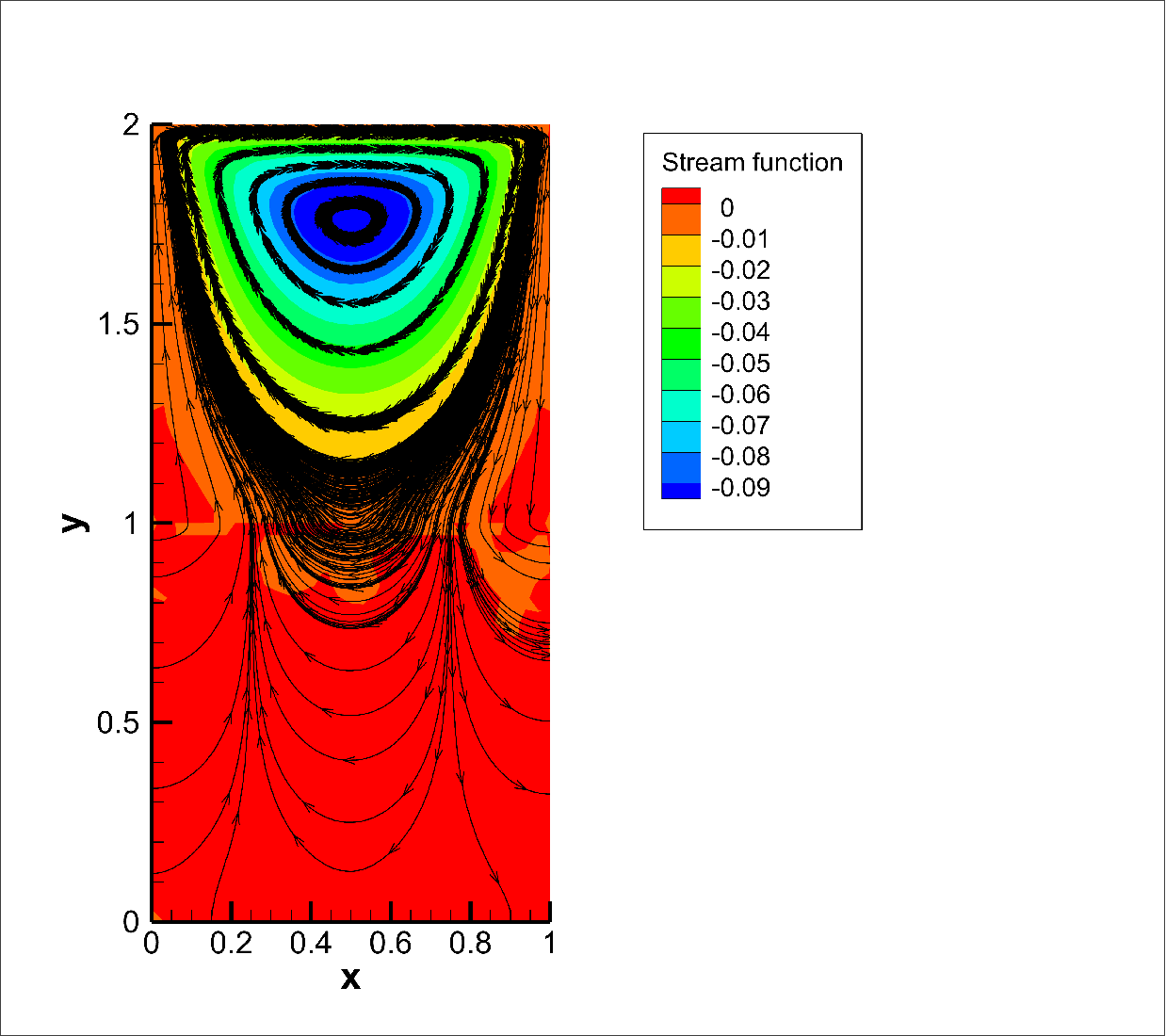} 
\caption{$\mu = 1$, $K = 1/1000$}\label{thelabel4}
\end{subfigure}
\caption{Stream lines for evolutionary NSD with $\mu=1,1/100,1/1000$ and $K=1,1/1000$ at final time $T=20$.}
  \label{fig:all1}
\end{figure}

\noindent \textbf{Example 7.3 (Benchmark Problem)} In this example, the lid-driven cavity flow test is considered on the computational regions $\Omega_f=(0,1)\times(1,2),$ $\Omega_p = (0,1)\times  (0,1) $, and $\Gamma=(0,1)\times \{1\} $ with $\f_f=(0,0)$ and $f_p=0$. The velocity on the upper boundary of the domain is $\bu=(1,0)$, and no-slip boundary conditions are considered on the other part of the domain. 
Below, the mentioned boundary conditions are presented:
{\scriptsize\begin{align*}
      \begin{cases}
     \bu = (1,0) & \text{on}\,\,(0,1)\times  \{2\},\\
     \bu =(0,0) & \text{on}\,\, \{0\}\times(1,2)\cup \{1\}\times(1,2), \\
     \phi=0 & \text{on}\,\, \partial\Omega_p\setminus\Gamma.
      \end{cases}
\end{align*}}
In order to perform the simulation, discontinuous $P_1-P_2-P_1$ elements for the fluid region and discontinuous $P_2-P_1$ elements for porous media have been chosen, respectively, with final time $T=20$, $\Delta t=0.01$, and $h=\frac{1}{32}$. Here, we have considered different $\mu=\{1,1/100,1/1000\}$, $\bK=\{1,1/1000\}$, $C_{11}=20$, and other parameters are $1$. In Figure \ref{fig:all}, we have shown a comparison along the lines $(0.5, y)$, $(x, 1.5)$, and $(x, 0.5)$ between unsteady and steady NSD problems. We observe that for large time, the unsteady NSD approximate solution converges to the steady NSD approximate solution. Figure \ref{fig:all1} yields the streamlines for the time-dependent problem for various values of $\mu$ and $\bK$. As values of $\mu$ decrease, we notice that the primary vortex moves toward the cavity’s center. Also, a second and third vortex may also develop near the interface. For the smaller value of the permeability, the streamlines of the numerical velocity show the stability of the method. From the figures, it can be observed that the numerical scheme can effectively simulate the fluid flow.

\vspace*{-5mm}
\section{Conclusion}\label{s7}
In this article, the LDG method for coupled unsteady NSD equations has been utilized, contributing novel insights. 
A locally conservative fully discrete numerical scheme for the coupled model has been proposed and analyzed along with its existence and uniqueness. Then, optimal error estimates in the energy norm have been derived with the help of standard projections. However, they are not adequate to achieve an optimal error estimate in the $L^\infty(\textbf{L}^2)$-norm. Therefore, modified Oseen operators have been introduced in the LDG setup, and their approximation properties have been established to conquer the difficulties. With the help of modified Oseen operators, optimal error estimates for velocity and piezometric head in the $L^\infty(\bL^2)$ and $L^\infty(L^2)$-norms have been derived, respectively. Also, an $L^\infty(L^2)$-norm error bound for pressure is obtained. For temporal discretization, the backward Euler method has been considered, and fully discrete error estimates have been established. Finally, numerical experiments have been performed to confirm our theoretical findings and prove the efficiency of the scheme.
\vspace*{-6mm}
\section*{Acknowledgements} 
\vspace*{-2mm}
The first author would like to express her heartfelt gratitude to the Department of Science and Technology (DST), Government of India, for the financial support (DST/INSPIRE Fellowship/IF210197). 
The second and third authors would like to express their gratitude to Anusandhan National Research Foundation (ANRF), Government of India, under the Core Research Grant project no. CRG/2023/006366.
\vspace*{-6mm}	
\section*{Declarations}
\vspace*{-2mm}
\noindent \textbf{Conflict of interest} The authors declare no competing interests.

\end{document}